\documentclass[a4paper,11pt]{article}
\usepackage[utf8]{inputenc}
\usepackage{fullpage}
\usepackage{amsmath, amssymb, amstext, amsfonts, amsthm, bbm, array, enumerate,scalerel}
\usepackage{mathrsfs, dsfont,mathabx}
\usepackage{nicefrac}
\usepackage[dvipsnames]{xcolor}
\usepackage{tikz}
\usepackage{stmaryrd}
\usepackage{tikz-3dplot}
\usepackage{float}
\usepackage{appendix}
\usepackage{longtable}
\usepackage[round]{natbib}
\usetikzlibrary{patterns}
\usetikzlibrary{plotmarks}
\usepackage{shuffle}
\usepackage{ stmaryrd }
 \usepackage{bbm}

\newtheoremstyle{boldremark}
    {\dimexpr\topsep/2\relax} 
    {\dimexpr\topsep/2\relax} 
    {}          
    {}          
    {\bfseries} 
    {.}         
    {.5em}      
    {}          

\RequirePackage[colorlinks,citecolor=blue, urlcolor=blue, linkcolor=blue ]{hyperref}

\theoremstyle{plain}
\newtheorem{theorem}{Theorem}[section]
\newtheorem{lemma}[theorem]{Lemma~}
\newtheorem{corollary}[theorem]{Corollary}
\newtheorem{proposition}[theorem]{Proposition}

\theoremstyle{definition}
\newtheorem{definition}[theorem]{Definition}
\newtheorem{example}[theorem]{Example}

\newtheorem*{assumption*}{\assumptionnumber}
\providecommand{\assumptionnumber}{}
\makeatletter

\theoremstyle{boldremark}
\newtheorem{remark}[theorem]{Remark}
\newtheorem{notation}[theorem]{Notation}

\newtheorem{assumpt}[theorem]{Assumption}

\numberwithin{equation}{section}

\def\Ex{\mathbb{E}}
\def\Rset{\mathbb{R}}
\def\Nset{\mathbb{N}}
\def\X{\mathbb{X}^{(2)}}

\newcommand{\cadlag}{c\`adl\`ag~}
\newcommand{\caglad}{c\`agl\`ad~}
\newcommand{\Cadlag}{C\`adl\`ag~}
\newcommand{\Ito}{It\^{o}~}
\newcommand{\R}{\mathbb R}
\renewcommand{\P}{\mathbb P}
\newcommand{\Q}{\mathbb Q}
\newcommand{\E}{\mathbb E}
\renewcommand{\L}{\mathbb L}
\newcommand{\Acal}{\mathcal A}
\newcommand{\e}{\epsilon}

\usepackage{scalerel,stackengine}
\stackMath
\newcommand\reallywidehat[1]{%
\savestack{\tmpbox}{\stretchto{%
  \scaleto{%
    \scalerel*[\widthof{\ensuremath{#1}}]{\kern.1pt\mathchar"0362\kern.1pt}%
    {\rule{0ex}{\textheight}}
  }{\textheight}%
}{2.4ex}}%
\stackon[-6.9pt]{#1}{\tmpbox}%
}
\title{Universal approximation theorems for continuous functions \\ of \cadlag  paths  and Lévy-type signature models}
\author{Christa Cuchiero \thanks{University of Vienna, Department of Statistics and Operations Research, Data Science @ Uni Vienna, Kolingasse 14-16, 1090 Wien, Austria, christa.cuchiero@univie.ac.at
} 
\and Francesca Primavera  \thanks{University of Vienna, Department of Statistics and Operations Research, Kolingasse 14-16, 1090 Wien, Austria, francesca.primavera@univie.ac.at}
\and Sara Svaluto-Ferro \thanks{University of Verona, Department of Economics, Via Cantarane 24, 37129 Verona, Italy, sara.svalutoferro@univr.it \newline
The authors gratefully acknowledge financial support 
through grant Y 1235 of the START-program.
}}
\date{}

\begin{document}
\maketitle

\begin{abstract}
We prove a universal approximation theorem that allows to approximate continuous functionals of \cadlag (rough) paths uniformly in time and on compact sets of paths via linear functionals of their time-extended signature. Our main motivation to treat this question comes from  signature-based models for finance that allow for the inclusion of jumps. Indeed, as an important application, we define a new class of universal signature models based on an augmented L\'evy process, which we call \emph{L\'evy-type signature models}. They extend continuous signature models for asset prices as proposed e.g.~by \cite{ASS:20} in several directions, while still preserving universality and tractability properties. To analyze this, we first show that the signature process of a generic multivariate L\'evy process is a polynomial process on the extended tensor algebra and then use this for pricing and hedging approaches within L\'evy-type signature models. 
\end{abstract}

\noindent\textbf{Keywords:} \cadlag rough paths, signature, universal approximation theorems, financial modeling with jumps\\
\noindent \textbf{MSC (2020) Classification:} 60L10, 60J76 
\tableofcontents

\section{Introduction}

Methods based on the signature of a path represent a non-parametric way for extracting characteristic features from time series data which is essential in machine learning tasks. This explains why these techniques are more and more applied in econometrics and mathematical finance, see e.g.~\cite{BHLPW:20, KLP:20, ASS:20, LNP:20, NSWLB:20, BHRS:21, CGGOT:21, MH:21, AGTZ:22}
and the references therein. Indeed, signature-based methods allow for data-driven  modeling approaches, while
stylized facts or 
first principles from mathematical finance can still be easily guaranteed. 

The notion of signature of a path goes back to \cite{CHEN:57, CHEN:77}  (see also \cite{M:54, F:81, S:88}) and plays a prominent role in the context of rough path theory, pioneered by \cite{L:98}. We also refer to the monographs by  \cite{FV:10} and \cite{FH:14}. It has initially been developed to deal with controlled differential equations driven by rough signals in a pathwise way. Indeed, signature can be understood as an enhancement of the input signal with iterated integrals, which in turn allows to express the solution of the rough differential equation as a continuous map of the input signal enhanced in this way.

More generally, the choice of the signature as feature map capturing the specific characteristics of the path, can be explained by a \emph{universal approximation theorem} (UAT) according to which continuous (with respect to  certain variation distances) functionals of \emph{continuous} paths can be approximated on compact sets of paths by linear functionals of the time-extended signature. This result is  however only proved for continuous paths and therefore leaves open the question of approximating continuous functionals of the more general set of \emph{\cadlag}paths, which have particular relevance when it comes to  financial modeling.  Based on advances on the signature of \cadlag paths by \cite{FS:17} and  using the Skorokhod-$J_1$-topology
on (Lie group valued) \cadlag paths, we here present a UAT that solves this question. Firstly, we prove that continuous path functionals over the time interval $[0,1]$ can be uniformly approximated on compact sets of paths by linear functionals of the time-extended signature evaluated at the final time $1$. Secondly, we extend this result by considering functionals of the stopped paths, stopped at all deterministic times $t \in [0,1] $, which can thus be associated with \emph{non-anticipative} path-functionals (see Remark \ref{rem:non anticip path funct}).  We then prove that the coefficients of the linear functionals do not depend on $t \in [0,1]$ but can be
chosen uniformly in time. This is subject of Theorem \ref{th: UAT non antici}, our main result in the first part of the paper.

Our principal motivation to treat this problem comes from  signature-based models for finance that allow for an inclusion of jumps.
Indeed, all signature models for asset prices that have been proposed so far, (see \cite{ASS:20, CGS:22, CGMS:23}), build on stochastic processes with continuous trajectories. In order to extend them to processes with \cadlag trajectories and to show universality properties among classical jump-diffusion models in finance, the UAT proved in Section \ref{sec:UAT} is important, as will become apparent from our definition of \emph{L\'evy-type signature models} below.

The approach that we follow consists in parameterizing the model itself or its characteristics as linear functionals of the signature of an augmented Lévy process $X$ that we call market's primary underlying process. More precisely,
$\mathbb{X}$ denotes the signature of the $N+2$ dimensional Lévy process (specified under some risk neutral measure $\mathbb{Q}$)
\begin{equation*}
	X_t:t\mapsto \left(t,W_t,\int_0^t\int_\Rset x\ (\mu-\nu)(ds,dx), \int_0^t\int_\Rset x^2\ \mu(ds,dx),\dots,\int_0^t\int_\Rset x^N\ \mu(ds,dx)\right),
\end{equation*} where $\mu$ denotes a homogeneous Poisson random measure with compensator $\nu$ and $W$ a one-dimensional Brownian motion.
As introduced in \cite{FS:17} for general semimartingales, we define the signature $\mathbb X$ of  the L\'evy process $X$, as solution of the Marcus stochastic differential equation (see \cite{KP:95} for more details) in the extended tensor algebra $T((\R^d))$ over $\Rset^d$,
\begin{align*}
	d\mathbb{X}=\sum_{i=1}^{d}\mathbb{X}\otimes \diamond \ dX^i,\qquad
	\mathbb{X}_0=(1,0,0,\dots)\in G((\Rset^{d})). 	
\end{align*}
We refer to Section \ref{sec:preliminaries} for the precise definitions of the involved quantities.
We then introduce \emph{Lévy-type signature models}  under some risk neutral measure $\mathbb{Q}$ via
\begin{equation}\label{eq:levymodelintro}
	S(\boldsymbol{\ell})_t=S_0+ \int_0^t\left(\sum_{|J|\leq n}\ell^J_W\langle \epsilon_{J},\mathbb{X}_{s^-}\rangle \right)dW_s+\int_0^t\int_\Rset \left(\sum_{|J|\leq n}\ell^J_\nu\langle \epsilon_{J},\mathbb{X}_{s^-}\rangle \right)x  \ (\mu-\nu)(ds,dx),
\end{equation}
for some parameters $ \ell^J_W,\ell^J_\nu\in \Rset$. Here $\langle \e_J,\mathbb X_t\rangle$ denotes the component $\mathbb X_t$ corresponding to the $J$-th basis element in $T((\R^d))$.

As already indicated above, due to the linear characteristics in $\mathbb{X}$, these models can be considered as \emph{universal}, since by our \cadlag versions of the UAT the  characteristics of a L\'evy driven SDE  can be approximated arbitrarily well, if they are continuous 
(with respect to the Skorokhod-$J_1$-topology) path functionals of the (enhanced) L\'evy process $X$. This applies in particular to  jump-diffusions as considered e.g.~in \cite[Section III.2c]{JS:87}.

To analyze tractability properties of models of the form \eqref{eq:levymodelintro}, we first show that the truncated signature process of a generic multivariate Lévy process  is a \emph{polynomial process on the truncated tensor algebra} and derive its expected value via  the so-called moment formula, thus solving a linear ODE (see \cite{CKT:12, FL:20}). Specifically, we prove that under moment assumptions on the Lévy measure, the expected signature of a Lévy process assumes the Lévy-Kintchine form described in \cite{FS:17}.

This in turn can be used to get  explicit formulas of the expected signature of $S$ in terms of the expected signature of $X$, which makes pricing of so-called sig-payoffs, which are nothing else than linear functionals of the signature of $(t, S_t)_{t \geq 0}$ (see \cite{LNP:20}), very efficient. A key step to obtain this result is to 
rewrite $S(\boldsymbol{\ell})_t$ itself as linear functionals of $\mathbb{X}_t$, a reformulation that we shall call \emph{sig-model representation} see Corollary \ref{cor1}).
Not only pricing of sig-payoffs (see Section \ref{sec:pricing}) be traced back to the signature of $X$, but also the quadratic hedging problem can be solved by computing the explicit form of the mean variance hedging strategy in terms of $\mathbb{X}$ (see Section \ref{sec:hedging}).

Note that by our UAT sig-payoffs approximate all $J_1$-continuous (possibly path-dependent) payoffs arbitrary well on compact set of paths, so that their prices and hedges can again be used as proxies. The corresponding coefficients of the sig-payoffs can be obtained by linear regression methods using typical price trajectories, e.g.~generated by several appropriate models with jumps. For calibration purposes this procedure is usually not accurate enough (comparable with  polynomial approximations, see also the Appendix in \cite{CGS:22}), but the obtained prices and hedges can for instance be used for variance reduction techniques. Let us remark here that one crucial advantage of signature based models, in particular the current L\'evy-type signature model, is that calibration via Monte Carlo methods is highly efficient. Indeed, L\'evy-type signature models share the same advantage as the continuous model of \cite{CGS:22}, since the
calibration task can be split into an offline sampling and a standard optimization. This is again due to the sig-model representation
which allows to precompute all signature samples of $\mathbb{X}$ offline and thus to avoid an (Euler)
simulation of trajectories in each optimization step. This adds further importance to the sig-model representation and suggests an interpretation of the compensator $\nu$ as hyperparameter (similarly to optimization tasks in the context of machine learning).
Note that standard Monte Carlo methods can then further be enhanced by variance reduction using sig-payoffs as explained above. 

Finally, in Section \ref{sec:change} we discuss equivalent measure changes that allow to keep the tractable structure of L\'evy-type signature models also under the physical measure $\mathbb{P}$.

Our results build on the theory of \cadlag rough paths as commenced in \cite{FS:17}, but related concepts like semimartingale signature (see \cite{T:21}) could be pursued as well. 
Let us finally remark that another strand of research related to finance where (càdlàg) rough paths play a significant role is robust and model-free finance, see for example \cite{PP:16, ACLP:21} and in particular \cite{ALP:21} for \cadlag rough paths foundations for robust finance. 
Indeed, in this context the theory of rough integration allows one to go beyond classical It\^o-integrals and to
overcome issues associated with null sets that are inherent in stochastic models, when e.g.~dealing with volatility uncertainty. This in turn opens the door to a worst case analysis and a
robustification of financial models and strategies. Note that we here do not pursue this direction but interpret all integrals in our financial model setup as It\^o-integrals.

The remainder of the paper is structured as follows. In Section \ref{sec:preliminaries}, we introduce the algebraic setting proper of rough path theory and recall the notion of \emph{weakly geometric \cadlag rough paths} and their \emph{Marcus signature}. Particular attention is also given to the Marcus signature of \cadlag semimartingales. Section \ref{sec:UAT} is dedicated to the universal approximation theorem. In  Section \ref{sec:Levy} we introduce L\'evy-type signature models, discuss their universality properties and derive pricing and hedging formulas based on polynomial technology. All the technical proofs are given in the Appendix. Finally, in order to enhance the accessibility of the paper, in Section \ref{sec:auxiliary} we include some auxiliary results on the Marcus signature of weakly geometric \cadlag rough paths.

\section{Preliminaries} \label{sec:preliminaries}
\subsection{Algebraic setting} \label{sec:auxiliary}

Fix $d\in \mathbb N$ and let $\Rset^d$  be the Euclidean space. The extended tensor algebra over $\Rset^d$ is defined by
	\[ T((\Rset^d)):=\prod_{n=0}^\infty(\Rset^d)^{\otimes n},
	\]
	where $(\Rset^d)^{\otimes n}$ denotes the $n$-tensor space of $\Rset^d$ with the convention $(\Rset^d)^{\otimes 0} := \Rset$. We equip
	$T((\Rset^d))$ with the standard addition $+$, tensor multiplication $\otimes$ and scalar multiplication.
	For $N\in\Nset$, the truncated tensor algebra is defined by
	\[
	T^N(\Rset^d):=\bigoplus_{n=0}^N(\Rset^d)^{\otimes n}.
	\]
	Elements of the extended tensor algebra will be denoted in blackboard bold face, e.g.~$\mathbbm{a}=(\mathbbm{a}^{(n)})_{n=0}^\infty\in T((\Rset^d)).$ Elements of the truncated tensor algebra instead, in bold face, e.g.~$\mathbf{a}=(\mathbbm{a}^{(n)})_{n=0}^N\in T^N(\Rset^d)$. For any $\mathbf{a}\in T^N(\Rset^d)$, we set 
 \begin{align}\label{eq:normTN}
 |\mathbf{a}|_{T^N(\Rset^d)}:=\max_{n=0,\dots,N}|\mathbbm{a}^{(n)}|_{(\Rset^d)^{\otimes n}}
 \end{align}
 and denote by $\rho$ the relative induced distance.
 Let $\pi_n:T((\Rset^d))\rightarrow(\Rset^d)^{\otimes n}$ be the map such that for $\mathbbm{a}\in T((\Rset^d)),$ $\pi_n(\mathbbm{a})=\mathbbm{a}^{(n)}$, and $\pi_{\leq N}:T((\Rset^d))\rightarrow T^N(\Rset^d)$ be such that for $\mathbbm{a}\in T((\Rset^d)),$ $\pi_{\leq N}(\mathbbm{a})=\mathbf{a}=(\mathbbm{a}^{(n)})_{n=0}^N$. For $c\in \Rset$, set 
$$T^N_c(\Rset^d):=\{\mathbf{a}\in T^N(\Rset^d)\colon \mathbbm{a}^{(0)}=c\}.$$
 The space $T_1^N(\Rset^d)$ is a Lie group under the tensor multiplication $\otimes$, truncated beyond level $N$. The neutral element with respect to $\otimes$ is $\mathbf{1}:=(1,0,\dots,0)\in T_1^N(\Rset^d)$. Moreover, for any $\mathbf{a}=(\mathbf{1}+\mathbf{b})\in T_1^N(\Rset^d)$, with $\mathbf{b}\in T_0^N(\Rset^d)$, its inverse is given by
 \begin{align}\label{eq:inv}
     \mathbf{a}^{-1}=\sum_{k=0}^N(-1)^k\mathbf{b}^{\otimes k}.
 \end{align}

The exponential and logarithm maps are defined as follows:
\begin{equation}\label{eqn3}
\begin{aligned}
			\exp^{(N)}:T&_{0}^N(\Rset^d)\rightarrow T_{1}^N(\Rset^d)\qquad\quad&	\log^{(N)}:T&_{1}^N(\Rset^d)\rightarrow T_{0}^N(\Rset^d)\\
		&\mathbf{b}\mapsto\mathbf{1}+\sum_{k=1}^N\frac{\mathbf{b}^{\otimes k}}{k!},\qquad&
		&\mathbf{1}+\mathbf{\mathbf{b}}\mapsto\sum_{k=1}^N(-1)^{k+1}\frac{\mathbf{b}^{\otimes k}}{k},
	\end{aligned}\end{equation}
	where the tensor multiplication is again always truncated beyond level $N$.
 
	Let $\mathfrak{g}^N(\Rset^d)$ be the \emph{free step-$N$ nilpotent Lie algebra} over $\Rset^d$, i.e.~
\begin{align*}
\mathfrak{g}^N(\Rset^d):=\{0\} \oplus \Rset^d\oplus[\Rset^d,\Rset^d]\oplus\dots\oplus \underbrace{[\Rset^d,[\Rset^d,\dots [\Rset^d,\Rset^d]]]}_{(N-1)\text{ brackets}}\subseteq
T_0^N(\Rset^d),
\end{align*}
where, for $\mathbf{a}\in T_0^{M}(\Rset^d)$, $1\leq M\leq N-1$, $\mathbbm{b}^{(1)}\in \Rset^d$, $[\mathbbm{b}^{(1)},\mathbf{a}]:=\mathbbm{b}^{(1)}\otimes \mathbf{a}-\mathbf{a}\otimes \mathbbm{b}^{(1)}$. 
Its image through the exponential map is a subgroup of $T_1^N(\Rset^d)$ with respect to $\otimes$. It is called \emph{free step-$N$ nilpotent Lie group} and is denoted by
\begin{align}\label{eq:GN}
    G^N(\Rset^d):=\exp^{(N)}(\mathfrak{g}^N(\Rset^d)).
\end{align}
We equip it with the so-called Carnot-Caratheodory (CC) norm $\|\cdot\|_{CC}$ and the induced (left-invariant) metric, denoted by $d_{CC}$.
We refer to Chapter 7 in \cite{FV:10} for more details. 

 Let $I = (i_1,\dots, i_n)$ be a multi-index with entries in $\{1,\dots,d\}$. Denoting by $\epsilon_1,\ldots,\epsilon_d$ the canonical basis of $\mathbb R^d$, we use the notation $|I|:=n$,
	\begin{equation*}
	S(I):=i_1+i_2+\dots+i_n,
	\qquad\text{and}\qquad
	\epsilon_I:=\epsilon_{i_1}\otimes \epsilon_{i_2}\otimes\dots\otimes  \epsilon_{i_n}.
	\end{equation*}
	Observe that $(\epsilon_I)_I$ is the canonical orthonormal basis of  $(\Rset^d)^{\otimes n}$.
	Furthermore, we denote by $\e_\emptyset$ the basis element of $(\R^d)^{\otimes 0}$ and set $|\emptyset|:=0$.  We also set $I':=(i_1,\dots, i_{n-1})$ for $n>1$, ${I'}=\emptyset$ for $n=1$,   $I^{''}:=(I')'$ for $n>1$, and use the convention $\e_{I''}=0$ for $n=1$.

 For two multi-indices $I\in \{1,\dots,d\}^n$, $J\in \{1,\dots,d\}^m$, and two indices of length $1$, $a,b\in\{1,\dots,d\}$, the shuffle product $\shuffle$ is defined recursively by
\begin{align*}
	&I\shuffle \emptyset=\emptyset\shuffle I=I,\\
	&(I,a)\shuffle (J,b)=((I\shuffle (J,b)),a)+(((I,a)\shuffle J),b),
\end{align*}
where $(I, a)$ denotes the concatenation of multi-indices. 

Given $\mathbbm{a}\in T((\Rset^d)),$ we write $\mathbbm{a}_{I}:=\langle \epsilon_{I},\mathbbm{a}\rangle$. We then define the following set
	\[
	\mathscr{L}_{\Rset^d}:=\text{span}\{\mathbbm a\mapsto \mathbbm a_{I}\colon  |I|\geq 0  \},
	\]
and call elements of $\mathscr{L}_{\Rset^d}$ \textit{linear functionals} on $T((\Rset^d))$.

The set of group-like elements is defined as follows
\[
G((\Rset^d)):=\{\mathbbm{a}\in T((\Rset^d)) \ | \ \pi_{\leq N }(\mathbbm{a})\in G^N(\Rset^d) \text{ for all }N \}.
\] 
Let $\mathbbm{a}\in G((\Rset^d))$ be a group-like element and $I\in \{1,\dots,d\}^n$, $J\in \{1,\dots,d\}^m$ two multi-indices. Then, we have that 
\begin{align}\label{eq:shuffle}
    \langle \epsilon_{I},\mathbbm{a}\rangle \langle \epsilon_{J},\mathbbm{a}\rangle =\langle \epsilon_{I}\shuffle \epsilon_{J},\mathbbm{a}\rangle
\end{align}
where $\epsilon_{I}\shuffle \epsilon_{J}:=\sum_{k=1}^K\epsilon_{I_k}$ where $K, I_k$ for $k=1, \ldots, K$ are determined via $I\shuffle J=\sum_{k=1}^K I_k$.

\subsection{\Cadlag rough paths and Marcus signature}\label{sec23}
In this section, we recall the definition of weakly geometric (w.g.)~\cadlag rough paths and Marcus signature, as originally formulated in \cite{FS:17}, to which we refer for a deeper discussion. We adopt here the Lie group-valued point of view and refer to Appendix \ref{sec:equi} for an equivalent definition.

Let 
$C([0, 1],G^N(\Rset^d))$ and $D([0, 1],G^N(\Rset^d))$ be the space of continuous and \cadlag maps
(paths), respectively, from the interval $[0, 1]$\footnote{The choice of any other interval $[a,b]$ does not effect any of the results presented in the rest of the paper.} into the metric space $(G^N(\Rset^d), d_{CC})$, both endowed with the  $J_1$-topology (see Definition~\ref{def:J1} and Remark~\ref{rem:J1cont}). For $\mathbf{X}\in D([0,1],G^N(\Rset^d))$, $s,t\in [0,1]$, $s\leq t$, we denote the path increments (with respect to the group structure) via $\mathbf{X}_{s,t}:=\mathbf{X}_{s}^{-1}\otimes \mathbf{X}_t$ and the path jumps as $\Delta\mathbf{X}_{t}:= \lim_{s \nnearrow t} \mathbf{X}_{s,t}$. 

For $p > 0$, we denote by $[p]$ the entire part of $p$. For a partition $\mathcal{D}=\{0=t_0<t_1<\dots<t_k=1\}$ of $[0,1]$ we write $\sum_{t_i\in \mathcal{D}}$ for the summation over all points in $\mathcal{D}$. The \emph{$p$-variation} of $\mathbf{X}\in D([0,1],G^N(\Rset^d))$ is then computed via
\begin{align*}
	\|\mathbf{X}\|_{p-var}:=\sup_{\mathcal{D}\subset[0,1]}\left(\sum_{t_i\in\mathcal{D}}d_{CC}(\mathbf{X}_{t_i,t_{i+1}})^p\right)^{\frac{1}{p}}.
\end{align*}
We denote by $C^p([0,1],G^N(\Rset^d))$  and $D^p([0,1],G^N(\Rset^d))$ the subspaces of $C([0,1],G^N(\Rset^d))$  and $D([0,1],G^N(\Rset^d))$, respectively, consisting of paths of finite $p$-variation

We then define the homogeneous variation distance on $D^p([0,1],G^N(\Rset^d))$ via
\begin{align}\label{eq:hommetric}
d_p(\mathbf{X},\mathbf{Y}):=\sup_{\mathcal{D}\subset[0,1]}\left(\sum_{t_i\in\mathcal{D}}d_{CC}(\mathbf{X}_{t_i,t_{i+1}},\mathbf{Y}_{t_i,t_{i+1}})^p\right)^{\frac{1}{p}}, \quad \text{for }\mathbf{X},\mathbf{Y}\in D^p([0,1],G^N(\Rset^d)).
\end{align}
\begin{definition}\label{def: roughpathsGROUP}
	Let $p\in[1,3)$. A \cadlag path $\mathbf{X}:[0,1]\rightarrow G^{[p]}(\Rset^d)$ is said to be a \emph{w.g.~\cadlag $p$-rough path over $\mathbb R^d$} if $\|\mathbf{X}\|_{p-var}<\infty$. For $p\geq 2$, a w.g.~\cadlag $p$-rough path is said to be \emph{Marcus-like} if for all $t\in [0,1]$,
 \begin{align*}
     \log^{(2)}(\Delta\mathbf{X}_t)\in \{0\}\oplus \Rset^d\oplus \{0\}.
 \end{align*}
 \begin{assumpt}\label{assumption}
     Throughout the paper, we assume that all the w.g.~\cadlag rough paths start at $\mathbf{1}\in G^{[p]}(\Rset^d)$, i.e. $\mathbf{X}_0=\mathbf{1}$.
 \end{assumpt}

 \begin{remark}
 \begin{enumerate}
     \item  Note that \emph{Marcus-like} means that jumps only occur at level $1$. Therefore, denoting $\pi_1(\mathbf{X})$ by $X$ (with values in $\mathbb{R}^{d}$), it holds that
    \[
  \log^{(2)}( \Delta \mathbf{X}_{t})=\lim_{s \nnearrow t} \log^{(2)}(\mathbf{X}_{s,t})=
    (0,  \Delta X_{t}, 0),
    \]     
    where the last $0$ stands here for the zero matrix in $[\mathbb{R}^{d}, \mathbb{R}^d]$.
    \item  The notion of w.g.~rough paths naturally extends to arbitrary low regularity. We restricted the presentation only to paths of finite $p$-variation for $p\in [1,3)$, as this is what matters the most in our stochastic analysis setting.
     
 \end{enumerate}

 \end{remark}
\end{definition}
By Lyons' extension theorem (see e.g.~Theorem~9.5 in \cite{FV:10}), every w.g.~continuous $p$-rough paths admits a unique extension to a path of finite $p$-variation with values in $G^N(\Rset^d)$, where $N>[p]$. 
A similar result holds also for w.g.~\cadlag $p$-rough paths. However, in the \cadlag setting, uniqueness is obtained provided that an additional algebraic condition is satisfied.
\begin{theorem}\label{minimal jump extension}
	Let $p\in[1,3)$ and $\Nset \ni N>[p]$. A w.g.~\cadlag $p$-rough path $\mathbf{X}:[0,1]\rightarrow G^{[p]}(\Rset^d)$ admits a unique extension to a \cadlag path $\mathbb{X}^{N}:[0,1]\rightarrow G^{N}(\Rset^d)$, such that $\mathbb{X}^{N}$ starts from $\mathbf{1}\in G^N(\Rset^d)$, is of finite $p$-variation with respect to the $CC$ metric $d_{CC}$ on $G^N(\Rset^d)$, and satisfies the following  condition:
	\begin{equation}\label{eqn11}
	\log^{(N)}(\Delta \mathbb{X}_t^N)=\log^{([p])}(\Delta \mathbf{X}_t) \quad  \text{for all }\ t\in[0,1].
	\end{equation}
	$\mathbb{X}^{N}$ is called \emph{minimal jump extension} of $\mathbf{X}$ in $G^N(\Rset^d)$.
\end{theorem}
\begin{remark}
    Equation \eqref{eqn11} 
has to be interpreted as follows: $\langle \e_I,\log^{(N)}(\Delta \mathbb{X}_t^N)\rangle=0$ for all $t\in[0,1]$, $|I|>1$ if $\mathbf{X}$ is Marcus-like, $|I|>[p]$ otherwise.
\end{remark}
Finally, the Marcus signature of a w.g.~\cadlag rough path is defined as the unique path extension provided by Theorem \ref{minimal jump extension}.
\begin{definition}
    Let $p\in[1,3)$ and $\mathbf{X}:[0,1]\rightarrow G^{[p]}(\Rset^d)$ be a w.g.~\cadlag $p$-rough path. The \textit{Marcus signature} of $\mathbf{X}$ is defined as the unique path
    \begin{align*}
        \mathbb{X}:[0,1]\rightarrow G((\Rset^d)),
    \end{align*}
    such that for all $\Nset \ni N>[p]$, $\pi_{\leq N}(\mathbb{X})=\mathbb{X}^N$, where $\mathbb{X}^N$ denotes the minimal jump extension of $\mathbf{X}$ in $G^N(\Rset^d)$.
\end{definition}
\begin{remark}
    If, in particular, $\mathbf{X}:[0,1]\rightarrow G^{[p]}(\Rset^d)$ is a continuous w.g.~rough path, $\mathbb{X}$ is nothing else than the unique path extension specified by Lyons' extension theorem.
\end{remark}
\begin{notation}
    From now on, we refer to $\mathbb{X}$ as \textit{signature} of $\mathbf{X}$ and to $\mathbb{X}^N$ as \textit{truncated signature} of order $N$ of $\mathbf{X}$.
\end{notation}

\vspace{0.2cm}

\subsubsection*{Marcus signature in a nutshell} 
Here, we report briefly the techniques employed in the proof of Theorem \ref{minimal jump extension} to see how the extension is constructed. It mostly relies on Marcus' idea of turning a \cadlag $\Rset^d$-valued path into a continuous one, see \cite{M:78,M:81} (which explains the choice of the name \textit{Marcus signature}). This is done by introducing an additional time interval at each jump time and replacing the jumps by a straight line which connects the states before and after the jump. In the current \cadlag rough path setting, where the state space is $G^{[p]}(\Rset^d)$, the jumps of $\mathbf{X}$ are connected by means of path-functions  mapping a pair $(\mathbf{x},\mathbf{y})\in(G^{[p]}(\Rset^d),G^{[p]}(\Rset^d))$ to a continuous path which links them (see Section 2.3 in \cite{CF:19}). In particular, for computing the minimal jump extension $\mathbb{X}^{N}$, the \textit{log-linear path-function} $\phi$ is used,  on which we will focus here. It is defined as follows:
\begin{align}\label{log linear path}
	\phi:G^{[p]}(\Rset^d)\times G^{[p]}(\Rset^d)&\rightarrow C([0,1],G^{[p]}(\Rset^d))\nonumber \\
	(\mathbf{x},\mathbf{y})&\mapsto\Big(s\mapsto\mathbf{x}\otimes \exp^{([p])}(s\log^{({[p]})}(\mathbf{x}^{-1}\otimes\mathbf{y} ))\Big).
	\end{align}
The outcome of this construction is a w.g.~continuous $p$-rough path $\mathbf{X}^\phi\in C^p([0,1],G^{[p]}(\Rset^d))$. Thus, by Lyons' extension theorem, it admits a unique extension $\mathbb{X}^{\phi,N}\in C^p([0,1],G^N(\Rset^d))$. The \cadlag extension is then obtained via a time-change of $\mathbb{X}^{\phi,N}$. That is, $\mathbb{X}^{N}:=\mathbb{X}^{\phi,N}_{\tau}$, for a pre-specified \cadlag function $\tau:[0,1]\rightarrow[0,1]$ such that the additional time intervals are removed. As a result of this construction, $\mathbb{X}^N$ can be explicitly computed by solving a so-called \emph{Marcus-type rough differential equation} (RDE) for w.g.~\cadlag rough paths (see \cite{CF:19}).

Roughly speaking, a solution of a Marcus-type RDE\footnote{Solutions of Marcus-type RDEs are also called \textit{geometric} solutions. In the presence of jumps, this notion has to be distinguished with the concept of \textit{forward} solution, introduced in \cite{FZ:17}.} driven by $\mathbf{X}$ is defined as the $\tau$-time changed solution of the corresponding continuous RDE driven by  $\mathbf{X}^{\phi}$, that is the continuous path built via $\phi$. 
For some given suitable vector field $V$, such a RDE is then denoted by
\begin{equation}\label{eq: generalMarcus}
    dY_t=V(Y_t) \diamond d \mathbf{X}_t.
\end{equation}
Observe that the same type of equation has also been considered in a semimartingale context by \cite{KP:95}. 

Recall now that by construction the minimal jump extension of Theorem \ref{minimal jump extension} satisfies $\mathbb{X}^{N}=\mathbb{X}^{\phi,N}_{\tau}$, and that $\mathbb{X}^{\phi,N}$ solves as continuous Lyons' extension of $\mathbf{X}^\phi$ a \emph{linear RDE} (see e.g.~Lemma~2.10 in \cite{L:07}, Exercise 4.6 in \cite{FH:14}). Hence, by the definition of the solution concept for Marcus-type RDEs, $\mathbb{X}^{N}$ solves \eqref{eq: generalMarcus} with $V$ being the identity. We summarize this whole reasoning in the following corollary.

\begin{corollary}\label{corominimal}
	Let $p\in[1,3)$, $\Nset \ni N>[p]$, and $\mathbf{X}:[0,1]\rightarrow G^{[p]}(\Rset^d)$ be a w.g.~\cadlag $p$-rough path. The minimal jump extension $\mathbb{X}^{N}$ with values in $G^N(\Rset^d)$ satisfies the Marcus-type RDE\footnote{We refer to Section \ref{sec:auxiliary} for more details on the Marcus-type RDE \eqref{Rdemarcus}.}
	\begin{align}\label{Rdemarcus}
		&d\mathbb{X}^{N}=\mathbb{X}^{N}\otimes \diamond d\mathbf{X},\qquad \mathbb{X}^{N}_0=\mathbf{1}\in G^N(\Rset^d), 
	\end{align}
whose explicit form can be written using \eqref{eqn3} and the notation of Section~\ref{sec23}  as 
\begin{align}\label{Rdemarcusexplicit}
	\mathbb{X}^{N}_t=1+\int_0^t\mathbb{X}^{N}_{s^-}\otimes d\mathbf{X}_s+\sum_{0<s\leq t}\mathbb{X}^{N}_{s^-}\otimes \{\exp^{(N)}(\log^{([p])}(\Delta \mathbf{X}_s))-\Delta \mathbf{X}_s\}.
\end{align}
The integral in \eqref{Rdemarcusexplicit} is understood as a rough integral and the summation term is well defined as an absolutely summable series.
\end{corollary}
\begin{remark}
	\label{rem:solution}  The minimal jump extension $\mathbb{X}^{N}$ is the \emph{unique} solution of equation \eqref{Rdemarcus}. Indeed, since equation \eqref{Rdemarcus} is an example of Marcus-type RDE 
with linear vector fields,  uniqueness follows from the uniqueness results for continuous linear RDEs (see Section 10.7 in \cite{FV:10}). Thus, given a w.g.~\cadlag $p$-rough path $\mathbf{X}$, its signature can also be defined without ambiguity as the unique solution of the following Marcus-type RDE in the extended tensor algebra:
	\begin{align*}
		d\mathbb{X}=\mathbb{X}\otimes \diamond d\mathbf{X},\qquad
		\mathbb{X}_0=(1,0,0,\dots)\in G((\Rset^d)).
	\end{align*}
\end{remark}

\subsection{Signature of \cadlag semimartingales}\label{sec:semimartroughpath}
\Cadlag semimartingales fit well into the theory of \cadlag rough paths. Indeed, every semimartingale admits a canonical lift which is a.s.~a Marcus-like \cadlag $p$-rough path for any $p > 2$. For the proof of the subsequent proposition, we refer to \cite{FS:17} and \cite{CF:19}.
\begin{proposition}\label{prop2}
	Let $p\in (2,3)$, $X:=(X_t)_{t\in[0,1]}$ be a $\Rset^d$-valued \cadlag semimartingale on a filtered probability space $(\Omega, \mathcal{F},(\mathcal{F}_t)_{t\in[0,T]},\mathbb{P})$, and $([X,X]^c_t)_{t\in[0,1]}$ its $(\Rset^d)^{\otimes 2}$-valued continuous quadratic variation. Set 
	\[
	\mathbb{X}^{(2)}_{0,t}:=\int_0^t X_{0,r^-}\otimes dX_r+\frac{1}{2}[X,X]^c_{t}+\frac{1}{2}\sum_{0< u \leq t}\Delta X_u\otimes\Delta X_u, \qquad t\in[0,1],
	\]
where the integral is understood in It\^o's sense. Define for all $t\in [0,1]$, $\mathbf{X}_t:=(1,X_{t}-X_0,\mathbb{X}^{(2)}_{0,t})$.
 Then $\mathbf{X}(\omega)\in D^p([0,1],G^2(\Rset^d))$  a.s. We call $\mathbf{X}$ Marcus lift of $X$.
\end{proposition}

The 
signature of semimartingales can be computed by adopting stochastic integration methods instead of resorting to rough integrals.
 Proposition~4.16 in \cite{CF:19} accounts for this. For completeness, we report an adaptation of this statement to our specific context.
 
\begin{proposition}\label{prop:sigsemimaritngale}
	Let $X$ be an $\Rset^d$-valued semimartingale and $\mathbf{X}$ its Marcus lift. It holds that the Marcus-type RDE
	\begin{align*}
		d\mathbb{X}=\mathbb{X}\otimes \diamond d\mathbf{X},\qquad
		\mathbb{X}_0=(1,0,0,\dots)\in G((\Rset^d))
	\end{align*}
coincides a.s.~with the Marcus SDE (see \cite{KP:95})
\begin{align}\label{MarcusSDE1}
	d\mathbb{X}=\mathbb{X}\otimes \diamond dX,\qquad 
	\mathbb{X}_0=(1,0,0,\dots)\in G((\Rset^d)).  
\end{align}
\end{proposition}

The explicit form of equation \eqref{MarcusSDE1}\footnote{Existence and uniqueness of the solution of the Marcus SDE \eqref{MarcusSDE1} are provided by Theorem 3.2 in \cite{KP:95}.} is given by
\begin{align*}
	\mathbb{X}_t=&1+\int_{0}^{t}\mathbb{X}_{s^-}\otimes dX_s+\frac{1}{2}\int_{0}^{t}\mathbb{X}_{s^-} \otimes  d[X,X]^c_s+\sum_{0<s\leq t} \mathbb{X}_{s^-}\otimes \{\exp(\Delta X_s)-\Delta X_s -1\},
\end{align*}
where
$\exp(\Delta X_s)=\sum_{k=0}^\infty \frac 1 {k! }(\Delta X_s)^{\otimes k}$.
This means that for each multi-index $I$ with entries in $\{1,\dots,d\}$, we get
\begin{align*}
		\langle\epsilon_{I},\mathbb{X}_t\rangle =&\int_{0}^t \langle\epsilon_{I'},\mathbb{X}_{s^-}\rangle dX_s^{i_{|I|}}
		+\frac{1}{2}\int_{0}^t\langle\epsilon_{I''},\mathbb{X}_{s^-}\rangle\  d[X^{i_{{|I|}-1}},X^{i_{|I|}}]^c_s \\
		&\qquad+\sum_{0<s\leq t}\sum_{\e_{I_1}\otimes \e_{I_2}=\e_I}\frac{1}{|I_2|!}\langle\epsilon_{I_1},\mathbb{X}_{s^-}\rangle \langle \e_{I_2}, (\Delta X_s)^{\otimes|I_2|}\rangle1_{\{|I_2|> 1\}}.
	\end{align*}
Observe that here we consider It\^o type integrals.
\begin{remark}
    \begin{enumerate}
        \item If $X$ is a $\Rset$-valued \cadlag semimartingale, by \Ito's formula we get that
        \begin{align*}
            \mathbb{X}_t=(1,X_t-X_0,\frac{(X_t-X_0)^2}{2},\dots,\frac{(X_t-X_0)^n}{n!},\dots)\in G((\Rset)).
        \end{align*}
        \item If $X$ is a continuous $\Rset^d$-valued semimartingale, $\mathbb{X}$ coincides with the Stratonovich signature, see e.g. \cite{FV:10}.
    \end{enumerate}
\end{remark}

\section{Universal approximation theorem}\label{sec:UAT}
In this section, we present a universal approximation theorems  (UAT) for continuous functionals of w.g.~\cadlag rough paths. We start by specifying the topology and the subspace of paths that will be considered.
\subsection{The subspace of time-extended paths and the topology} \label{sec:time-extended}
\begin{definition}\label{def:J1}
Let $(E,d)$ be a metric space and $D([0,1],E)$ the space of \cadlag paths on it. Denote by $\Lambda$ the set of all strictly increasing bijections of $[0,1]$ to itself. The \textit{Skorokhod $J_1$-metric} on $D([0,1],E)$ is defined via
\begin{equation*}
	\sigma_\infty(X,Y):=\inf_{\lambda\in \Lambda}\{|\lambda|\vee\sup_{s\in[0,1]}d(X_{\lambda(s)},Y_s)\}, \qquad X,Y\in D([0,1],E),
\end{equation*}
where,  $|\lambda|:=\sup_{s\in[0,1]}|\lambda(s)-s|$. 
\end{definition}
\begin{remark}\label{rem:J1cont}
    We endow the space of \cadlag paths considered throughout the paper with the $J_1$-topology, induced by the metric $\sigma_\infty$. In particular, if $(E,d)=(G^{N}(\Rset^{d}),d_{CC})$, $N\in \Nset$, the Skorokhod $J_1$-metric on $D([0,1],G^{N}(\Rset^{d}))$ reads
\begin{align}\label{eq:skordCC}
    \sigma_\infty(\mathbf{X},\mathbf{Y}):=\inf_{\lambda\in \Lambda}\{|\lambda|\vee \sup_{s\in[0,1 ]}d_{CC}(\mathbf{X}_{\lambda(s)},\mathbf{Y_s})\}, \qquad \mathbf{X},\mathbf{Y}\in D([0,1],G^{N}(\Rset^{d})).
\end{align}
\end{remark}
Next, we introduce a particular subspace of w.g.~\cadlag rough paths.
\begin{definition}\label{def:time-paths}
Let $p\in [1,3)$. The subspace of \textit{time-extended w.g.~\cadlag $p$-rough paths} is defined as follows:\footnote{For the explanation of the index $-1$ used here we refer to Notation \ref{notation} below.}
     \begin{enumerate}
         \item If $p\in[1,2)$, $$\widehat{D}^p([0,1],G^{1}(\Rset^{d+1})):=\Big\{\widehat{\mathbf{X}}\in \ D^p([0,1],G^{1}(\Rset^{d+1}))\ | \text{ for all }t\in [0,1], \ \langle \e_{-1},\widehat{\mathbf{X}}_t\rangle :=t\Big\}.$$
         \item If $p\in[2,3)$,
     \end{enumerate}
\begin{align}\label{eq:Wset}
\widehat{D}^p([0,1],G^{2}(\Rset^{d+1})):=\Big\{\widehat{\mathbf{X}}\in \ &D^p([0,1],G^{2}(\Rset^{d+1}))\ | \text{ for all }t\in [0,1], \ \langle \e_{-1},\widehat{\mathbf{X}}_t\rangle :=t, \\ \nonumber 
& \text{ for all }i\in \{-1,1,\dots,d\}, \ \langle \e_{(i,-1)},\widehat{\mathbf{X}}_t\rangle :=\int_0^t\langle \e_{i},\widehat{\mathbf{X}}_{s-}\rangle ds\Big\}.
\end{align}
\end{definition}
\begin{notation}\label{notation}
In Definition \ref{def:time-paths} we use the index $-1$ to denote the time component of a time-extended w.g.~\cadlag rough paths. This notation will be used throughout the paper.
\end{notation}

\begin{remark} \phantomsection \label{rem: time-extended paths}
\begin{enumerate} 
    \item By Assumption \ref{assumption}, for all $\widehat{\mathbf{X}}\in \widehat{D}^p([0,1],G^{[p]}(\Rset^{d+1}))$, it holds that $\widehat{\mathbf{X}}_0=\mathbf{1}\in G^{[p]}(\Rset^{d+1})$.
    \item \label{rem: time-extended paths ii}
    Observe that due to the smoothness of the identity path $X_t:=t$ for $t\in[0,1]$, given a w.g.~\cadlag rough path $\mathbf{X}\in D^p([0,1],G^{[p]}(\Rset^d))$, it is always possible to build a time-extended one $\widehat{\mathbf{X}}$ with values in $G^{[p]}(\Rset^{d+1})$, i.e.  $\widehat{\mathbf{X}}\in \widehat{D}^p([0,1],G^{[p]}(\Rset^{d+1}))$. Trivially, the reverse is also true.

    \item 
    Let $p\in[2,3)$. If $\mathbf{X}, \tilde{\mathbf{X}}\in D^p([0,1],G^{2}(\Rset^{d+1})) $ are such that $\pi_1(\mathbf{X})=\pi_1(\tilde{\mathbf{X}})$, then one can show that there exists a \cadlag path of finite $\frac{p}{2}$ variation $S:[0,1]\rightarrow[\Rset^{d+1},\Rset^{d+1}]$,  such that for all $t\in[0,1]$,$$\tilde{\mathbf{X}}_t=\mathbf{X}_t+(1,0,S_t-S_0).$$
 This in particular implies that for any $\mathbf{X}\in D^p([0,1],G^{2}(\Rset^{d+1})) $, there exists infinitely many $\tilde{\mathbf{X}}\in D^p([0,1],G^{2}(\Rset^{d+1})) $ with the same $\pi_1(\mathbf{X})$.
 However, when some components of $\pi_1(\mathbf{X})$ are sufficiently regular, the corresponding components of the second level $\pi_2(\mathbf{X})$ are frequently fixed to be the
paths specified via the Young integral (\cite{YOUNG}), considered as a canonical choice. In the definition of the set \eqref{eq:Wset}, we follow this line. The set of time-extended w.g.~\cadlag rough paths is in fact the subset of w.g.~\cadlag rough paths such that the first component (indexed by $-1$) of the underlying path is the time and all the components of the second level path involving the index $-1$ are defined via the Young integral. Observe that by the 
shuffle property \eqref{eq:shuffle}, for all $\widehat{\mathbf{X}}\in \widehat{D}^p([0,1],G^{2}(\Rset^{d+1}))$, $i\in \{-1,1,\dots,d\}$, $t\in [0,1]$, it holds that   $$\langle \e_{(-1,i)},\widehat{\mathbf{X}}_t\rangle=t \ \langle \e_{i},\widehat{\mathbf{X}}_t\rangle-\langle \e_{(i,-1)},\widehat{\mathbf{X}}_t\rangle.$$
This specific choice is needed in order to recover consistency with the classical Young integral (see Lemma \ref{Rough=Young}), which is key for showing condition \ref{itiii} in the proof of Proposition \ref{UAT1} (see also Corollary \ref{coro:uniqueness}).
\end{enumerate}
\end{remark}

\subsection{UAT for continuous path functionals}
We start by proving that continuous path functionals can be uniformly approximated on compact sets by linear functionals of the time-extended signature evaluated at the final time.  
\begin{proposition}\label{UAT1}
Let $p\in[1,3)$, $K\subset \widehat{D}^p([0,1],G^{[p]}(\Rset^{d+1}))$ a compact subset, bounded with respect to the $p$-variation norm, and $F:\widehat{D}^p([0,1],G^{[p]}(\Rset^{d+1}))\rightarrow \Rset$ a continuous functional. For each $\widehat{\mathbf{X}}\in K$, denote by $\mathbb{\widehat{X}}$ its signature. Then, for every $\varepsilon>0$ there exists a linear functional $l\in \mathscr{L}_{\Rset^{d+1}}$ such that\[
	\sup_{\widehat{\mathbf{X}}\in K}|F(\widehat{\mathbf{X}})-l(\mathbb{\widehat{X}}_1)|\leq\varepsilon.
	\]
\end{proposition}
The proof of Proposition \ref{UAT1} can be found in Appendix~\ref{appendixA:UAT1}.

\begin{remark}
\begin{enumerate}
\item On the subset of continuous paths the $J_1$-metric corresponds to the sup one defined via
$$
    d_{\infty}(\mathbf{X},\mathbf{Y}):=\sup_{s\in[0,1 ]}d_{CC}(\mathbf{X}_{s},\mathbf{Y_s}), \qquad \mathbf{X},\mathbf{Y}\in C([0,1],G^{[p]}(\Rset^{d+1})).
$$
Therefore, Theorem~\ref{UAT1} states in particular that $d_{\infty}$-continuous functionals of \textit{continuous paths} can be approximated, uniformly on compact sets, by linear functionals of the time-extended signature evaluated at the final point. For different versions of approximations theorems for continuous functionals of \textit{continuous} paths we refer to \cite{LLN:13,KO:19,LNP:20}.
	\item Requiring the set $K$ to be bounded with respect to the $p$-variation norm is not redundant. For instance, on the subset of continuous paths the sup distance $d_{\infty}$ is dominated by the $p$-variation one. This implies that not every $J_1$-compact set is bounded with respect to the $p$-variation norm. 
\item Propositon \ref{UAT1} can be reformulated with respect to the rough paths variant of the Skorokhod $M_1$-topology introduced in Section 3.2 in \cite{CF:19}. Indeed, by Theorem 3.13 in \cite{CF:19} one can see that the set A given by \eqref{setA} is a linear subspace of the continuous functions from $K$ to $\Rset$, when on $K$ the Skorokhod $M_1$-topology is considered. Notice however that neither version of these approximation results implies the other, as in general, the two Skorokhod topologies are not comparable. See Proposition 3.10 and Remark 3.6 in \cite{CF:19} for more details.
\end{enumerate}

\end{remark}

Next, for illustrative purposes, we also provide a simple adaptation of Proposition \ref{UAT1} to the semimartingales setting. Rember that for a semimartingale $X$ its Marcus lift ${\mathbf{X}}$ is computed as in Proposition \ref{prop2} and its signature ${\mathbb{X}}$ is the solution of the Marcus SDE \eqref{MarcusSDE1}. 
\begin{corollary}\label{coro:UAT1 semimart}
    Under the hypothesis of Proposition \ref{UAT1}, for $p\in(2,3)$. Let $X=(X_t)_{t\in [0,1]}$ be an $\Rset^{d}$-valued semimartingale, set $\widehat{X}:=(t,X_t)_{t\in[0,1]}$ and denote by $\widehat{\mathbf{X}}$ its Marcus lift and by $\widehat{\mathbb{X}}$ its signature. Then, for every $\varepsilon>0$ there exists a linear functional $l\in \mathscr{L}_{\Rset^{d+1}}$ such that\[
	\sup_{\widehat{\mathbf{X}}(\omega)\in K}|F(\widehat{\mathbf{X}}(\omega))-l(\mathbb{\widehat{X}}_1(\omega))|\leq\varepsilon
	\]
 almost surely.
\end{corollary}
\begin{proof}
     Since the \Ito integral with respect to semimartingales of finite variation coincides a.s. with the Young integral, it holds that $\widehat{\mathbf{X}}(\omega)\in \widehat{D}^p([0,1],G^{2}(\Rset^{d+1}))$ a.s.. The rest of the proof follows as in Proposition \ref{UAT1}.
\end{proof}

\begin{remark}
Note that the continuous functional $F$ in Corollary \ref{coro:UAT1 semimart}   can of  course be a functional that only involves $X$ and not its Marcus lift.
\end{remark}

As a direct result of the proof of Proposition \ref{UAT1}, the following statement concerning \textit{uniqueness of the signature of time-extended w.g.~\cadlag rough paths} holds true.
\begin{corollary}\label{coro:uniqueness}
    Let $\widehat{\mathbf{X}}\in \widehat{D}^p([0,1],G^{[p]}(\Rset^{d+1}))$ and denote by $\widehat{\mathbb{X}}_1$ its signature evaluated at the final point $1$. Then, $\widehat{\mathbb{X}}_1$ uniquely determines $\widehat{\mathbf{X}}$. 
\end{corollary}
\begin{proof}
     This result follows directly from the proof of condition \ref{itiii} in Proposition \ref{UAT1}. The reasoning therein does not depend on the set $K$, nor on the metric $\sigma_{\infty}$, showing that a time-extended w.g.~\cadlag $p$-rough path is completely characterized by its signature evaluated at the final time $1$.
\end{proof}

\begin{remark}
    Let $\widehat{X}:=(t,X_t)_{t\in[0,1]}$ be a time-extended \cadlag semimartingale with $\widehat{X}_0=0$. By Corollary \ref{coro:uniqueness} we get that a.s. $\widehat{\mathbb{X}}_1(\omega)$ uniquely determines $X(\omega)$.
\end{remark}
 We now present the main result of the first part of the paper: an extended version of Proposition \ref{UAT1} which leads to an approximation of paths functionals not only on compact sets of paths but also  uniformly in time. For its formulation, we first need to introduce the notion of a \textit{stopped w.g.~\cadlag rough path}.
 \begin{definition}\label{def:stopped paths}
Let $p\in[1,3)$.  Given ${{\mathbf{X}}}\in {D}^p([0,1],G^{[p]}(\Rset^{d})))$ and $t\in [0,1]$, we define the \textit{stopped w.g.~\cadlag rough path of ${\mathbf{X}}$} stopped at time $t$ as the \cadlag path ${\mathbf{X}}^t:[0,1]\rightarrow G^{[p]}(\Rset^{d})$  given by
$${\mathbf{X}}^t_u:=
         {\mathbf{X}}_u1_{\{u<t\}}+{\mathbf{X}}_t1_{\{u\geq t\}}.$$
 \end{definition}
The proof of the following theorem can be found in Appendix \ref{appendixA:UAT non antici}.

\begin{theorem}\label{th: UAT non antici}
	Let $p\in[1,3)$, $K\subset \widehat{D}^p([0,1],G^{[p]}(\Rset^{d+1}))$ a compact subset, bounded with respect to the $p$-variation norm, and $F:{D}^p([0,1],G^{[p]}(\Rset^{d+1}))\rightarrow \Rset$ a continuous path functional. For each $\widehat{\mathbf{X}}\in \widehat{D}^p([0,1],G^{[p]}(\Rset^{d+1}))$, denote by $\widehat{\mathbb{X}}$ its signature. Then, for every $\varepsilon>0$ there exists a linear functional ${l}\in \mathscr{L}_{\Rset^{d+1}}$ such that\[
	\sup_{t\in[0,1]}\sup_{{\widehat{\mathbf{X}}}\in K}|F({\widehat{\mathbf{X}}^t})-{l}(\widehat{\mathbb{{X}}}_t)|\leq\varepsilon.
	\]
\end{theorem}

\begin{remark}\phantomsection\label{rem:non anticip path funct}
\begin{enumerate}
    \item  Observe that by Definition~\ref{def:stopped paths}, for each ${\mathbf{X}}\in {D}^p([0,1],G^{[p]}(\Rset^{d}))$  
    the map ${\bf  X}\mapsto F({\mathbf{X}}^t)$
   just depends on ${\bf X}$ through its values $({\bf X}_s)_{s\in [0,t]}$ on $[0,t]$. Therefore, for a given path functional $F:{D}^p([0,1],G^{[p]}(\Rset^{d}))\rightarrow \R$, it is always possible to associate a so-called non-anticipative path functional $\tilde F:[0,1]\times D^p([0,1], G^{[p]}(\mathbb{R}^d))$ (see e.g.~\cite{DUP:09,CONT:10,CONTFOUR:13, CONTANNA:17}) such that for all $t\in [0,1]$, ${\bf X}\in D^p([0,1], G^{[p]}(\mathbb{R}^d))$,
   \begin{align*}
      \tilde{F}(t, {\bf X}):=F({\bf X}^t).
   \end{align*}
    \item A similar result in the setting of \textit{continuous} rough paths has been formulated in \cite{KLP:20}, where, however, different techniques are used. Indeed, in contrast to our approach, \cite{KLP:20} consider non-anticipative functionals on the space of continuous rough paths
  equipped with an appropriate topology (fully clarified in \cite{BHRS:21}).
\end{enumerate}
\end{remark}

Similarly as in Corollary \ref{coro:UAT1 semimart}, in the semimartingales setting, Theorem \ref{th: UAT non antici} reads as follows.
\begin{corollary}
    Under the hypothesis of Theorem \ref{th: UAT non antici}, for $p\in (2,3)$. Let $X=(X_t)_{t\in [0,1]}$ be an $\Rset^{d}$-valued semimartingale, set $\widehat{X}:=(t,X_t)_{t\in[0,1]}$ and denote by $\widehat{\mathbf{X}}$ its Marcus lift and by $\widehat{\mathbb{X}}$ its signature. Then, for every $\varepsilon>0$ there exists a linear functional ${l}\in \mathscr{L}_{\Rset^{d+1}}$ such that\[
	\sup_{t\in[0,1]}\sup_{{\widehat{\mathbf{X}}(\omega)}\in K}|F({\widehat{\mathbf{X}}^t}(\omega))-{l}(\widehat{\mathbb{{X}}}_t(\omega))|\leq\varepsilon
	\]
 almost surely.
\end{corollary}

Finally, we conclude this section with some illustrative examples of path functionals for which Theorem \ref{th: UAT non antici} holds true.

\begin{example}\label{example functionals UAT}
(i) For $|I|\leq 2$, consider the path functional $F:{D}([0,1],G^2(\Rset^{d+1}))\to\R$ given by 
$${F({{\mathbf{X}}}):= \sup_{s\leq 1}|\langle \e_I,{{\mathbf{X}}}_s\rangle|} 
$$
and note that for all $t\in[0,1]$,
$
			{F({{\mathbf{X}}}^t)=\sup_{s\leq t}|\langle \e_I,{{\mathbf{X}}}_s\rangle|} 
$.
Furthermore, by Proposition VI.2.4 in \cite{JS:87}, $F$ is continuous with respect to the Skorokhod $J_1$-topology at all paths $\mathbf{X}$ that do not admit a jump at time $1$.  
  Let $K$ be as in Theorem~\ref{th: UAT non antici} and assume that all elements in $K$ do not jump at $1$.  By Theorem~\ref{th: UAT non antici}, for every $\varepsilon>0$ there is a map $ l\in \mathscr L_{\Rset^{d+1}}$ such that
		\begin{align*}
			\sup_{t\in[0,1]}\sup_{{\widehat{\mathbf{X}}}\in K}|F(\widehat{\mathbf{X}}^t)-{l}(\widehat{\mathbb{{X}}}_t)|=\sup_{t\in[0,1]}\sup_{{\widehat{\mathbf{X}}}\in K}|\sup_{s\leq t}|\langle \e_I,\widehat{\mathbf{X}}_s\rangle|
			-{l}(\widehat{\mathbb{{X}}}_t)|\leq \varepsilon.
		\end{align*}
(ii) Let $Y$ be the solution of the Marcus-type RDE
\begin{align}\label{eq:example MarcusRDE}
		dY_s=V(Y_s)\diamond d{\mathbf{X}}_s, \quad Y_0=y_0,
		\end{align}
		for some suitable vector field $V$ and let ${\mathbf{X}}$ be a w.g.~\cadlag $p$-rough paths, for $p\in[2,3)$. Consider the path functional $F:{D}^p([0,1],G^2(\Rset^{d+1}))\to\R$ given by
  $$ {F({{\mathbf{X}}})}:=Y_1.
		$$
    For each $t\in [0,1]$, denote by $\widehat{Y}^t_1$ the solution map evaluated at time $1$ of the Marcus-type RDE \eqref{eq:example MarcusRDE} driven by the stopped time-extended path ${\widehat{\mathbf{X}}}^t\in {D}^p([0,1],G^2(\Rset^{d+1}))$. Observe that $F({\widehat{\mathbf{X}}}^t)=\widehat{Y}^t_1=\widehat{Y}_t$.
Moreover, by Proposition 3.18 in 
\cite{CF:19} and a similar argument as in the proof of Proposition \ref{UAT1}, $F$
		is continuous with respect to the Skorokhod $J_1$-topology on sets of bounded $p$-variation norm. 
  Thus, for each $\varepsilon>0$ and $K$ satisfying the conditions of Theorem \ref{th: UAT non antici},   there exists ${l}\in\mathscr{L}_{\Rset^{d+1}}$ such that
		\begin{align*}
\sup_{t\in[0,1]}\sup_{{\widehat{\mathbf{X}}}\in K}|F(\widehat{\mathbf{X}}^t)-{l}(\widehat{\mathbb{{X}}}_t)|=\sup_{t\in[0,1]}\sup_{{\widehat{\mathbf{X}}}\in K}|\widehat{Y}_t
	-{l}(\widehat{\mathbb{{X}}}_t)|\leq \varepsilon.
		\end{align*}
Therefore, solutions of Marcus-type RDEs can be approximated uniformly in time by linear functionals of the signature of the driving signal\footnote{By Theorem 5.3 in \cite{FZ:17}, the same type of argument applies to solutions of RDEs of \textit{forward-type}. Notice in particular that this concept of RDE is in line with the \Ito theory of  \cadlag semimartingales (see e.g.~Proposition 6.9 in \cite{FZ:17}).}.
\end{example}

\section{Lévy-type signature models}\label{sec:Levy}
\subsection{The market's primary process and the model}
In this section we introduce the class of Lévy-type signature models. The key object is a multidimensional process, which will have the role of encoding all the randomness of the model. Because of this, it will be called \emph{market's primary process}.

For its definition, we consider a standard Brownian motion $W$ and a homogeneous Poisson random measure $\mu$ on a filtered probability space $(\Omega, \mathcal{F},(\mathcal{F}_t)_{t\in[0,T]},\mathbb{Q})$, whose filtration $(\mathcal{F}_t)_{t\in[0,T]}$, is the augmented one generated by $W$ and $\mu$. We suppose that $\mu$ and $W$ are independent and that $\mu$ has intensity measure $\nu(dt,dx)=dt\times F(dx)$. Here, $dt$ denotes the Lebesgue measure on $([0,T],\mathcal{B}([0,T]))$ and $F(dx)$ a $\sigma$-finite measure on $(\Rset,\mathcal{B}(\Rset))$ satisfying $F(\{0\})=0$ and 
\begin{align}\label{eq:momentLevy}
     \int _{\R} |x|^kF(dx) <\infty \qquad \text{for all } k\geq 2.
\end{align}
\begin{definition}\label{def2}
Fix  $N\geq 2$. The \emph{market's primary process} is  the $\Rset^{N+2}$-valued process $X$ given by
$$
	X_t:= \left(t,W_t,\int_0^t\int_\Rset x\ (\mu-\nu)(ds,dx), \int_0^t\int_\Rset x^2\ \mu(ds,dx),\dots,\int_0^t\int_\Rset x^N\ \mu(ds,dx)\right).
$$
\end{definition}
Observe that we do not use the $\widehat{\cdot}$ notation, even though the first component of $X$ is $t$. Moreover, from now on we will use the following numbering
for the components of the $\Rset^{N+2}$-valued market's primary process $X$
 \begin{align*}\label{components X}
 	&X^{-1}_t=t,\\
 	&X^{0}_t=W_t,\nonumber \\
 	&X^{1}_t=\int_0^t\int_\Rset x\ (\mu-\nu)(ds,dx),\nonumber \\
 	&X^{k}_t=\int_0^t\int_\Rset x^k\ \mu(ds,dx),\qquad  \text{ for }2\leq k\leq N.\nonumber 
 \end{align*}

The next proposition states that the market's primary process is a Lévy process.
\begin{proposition}\label{Xislevy}
	The $\Rset^{N+2}$-valued market's primary process   $X$  introduced in Definition~\ref{def2} is a Lévy process and its triplet $(b^{X},C^{X}, F^{X})$, associated with the ``truncation function"\footnote{See Definition II.2.3 and Corollary II.2.38 in \cite{JS:87}.} $\chi (y)=y$, is given by 
	\begin{align*}
		&b^{X}=(1,0,0,\int_\Rset x^2 F(dx),\dots, \int_\Rset x^N F(dx)),\\ \\
		&C^{X}_{i,j}=0 \ \text{if}\  i\neq 0\text{ or }j\neq 0, \  \ C^X_{0,0}=1,\\ \\
		&F^{X}(dy)=h_*F(dx),
	\end{align*}
	where $h(x):=(0,0,x,x^2,x^3,\dots,x^N)$.
\end{proposition}
\begin{proof}
	Observe that $X$ is the sum of the two independent $\Rset^{N+2}$-valued processes
$Y$ and $Z$ given by
	$
	Y_t=\left(t,W_t,0,\dots,0\right)
	$
	and
	\[
	Z_t=\left(0,0,\int_{0}^{t}\int_{\Rset}x\  ({\mu}-\nu)(ds,dx),\int_0^t\int_{\Rset}x^2 \mu(ds,dx),\dots,\int_0^t\int_{\Rset}x^N \mu(ds,dx)\right).
	\]
	Using the same numbering for the components, we get that $Y$ is a Lévy process with triplet $(b^Y,C^Y,0)$, where 
	$
	b^Y=(1,0,\dots,0)$, $C^Y_{i,j}=0$ for $i\neq 0$ or $j\neq 0$ and $C^Y_{0,0}=1$.
	
	Next, we verify that $Z$ is a well defined Lévy process. Let $M$ be the random measure on $\Rset_+\times \Rset^{N+2}$, given by $
	M:=H_*\mu$ for $H(t,x):=(t,h(x))$. Since $\mu$ is a Poisson random measure on $\Rset_+\times \Rset$, and $H$ is a measurable map, $M$ is a Poisson random measure on $\Rset_+\times \Rset^{N+2}$ with intensity $$\nu^Z(dt,dy)=H_*(dt\times F(dx))=dt\times (h_*F(dx)).$$
	 Moreover, since
	$
	(h_*F)(\{0\})=F(\{0\})=0,
	$ and 
	\[
	\int_{\Rset^{N+2}}(1\wedge\|y\|^2)h_*F(dy)=\int_{\Rset}(1\wedge\|h(x)\|^2)F(dx)<\infty,
	\]
	 we also have that $h_*F$ is a Lévy measure on $\Rset^{N+2}$. Using that
	 $$Z=b^Zt+	\int_{0}^t\int_{\Rset^{N+2}}y \ (M-\nu^Z)(ds,dy),$$
	for $b^Z=\left(0,0,0,\int_{\Rset}x^2 F(dx),\dots,\int_{\Rset}x^N F(dx)\right)$, we can conclude that $Z$ is a well defined Lévy process, whose triplet with respect to $\chi(y)=y$ is given by $\left(b^Z,0,h_*F)\right)$. The claim follows.	
\end{proof}
\begin{remark}\label{rem:finitmomentX}
    Observe that by condition \eqref{eq:momentLevy}, the $\Rset^{N+2}$-valued Lévy process $X$ introduced in Definition~\ref{def2} has finite $k$-moments, for all $k\in \Nset$.
\end{remark}
 Since $X$ is a Lévy process, following the discussion of Section \ref{sec:semimartroughpath}, we can compute its signature $\mathbb{X}$ by solving the $T((\Rset^{N+2}))$-valued Marcus SDE given by
\begin{align*}
	d\mathbb{X}&=\mathbb{X}\otimes \diamond \ dX,\qquad
	\mathbb{X}_0=(1,0,0,\dots)\in T((\Rset^{N+2})), \nonumber	
\end{align*}
whose component-wise version is given by 	
\begin{equation}\label{component-wise}
\begin{aligned}
		\langle\epsilon_{I},\mathbb{X}_t\rangle =&\int_{0}^t \langle\epsilon_{I'},\mathbb{X}_{s^-}\rangle dX_s^{i_{|I|}}
		+1_{\{i_{{|I|}-1}=i_{|I|}=0\}}\frac{1}{2}\int_{0}^t\langle\epsilon_{I''},\mathbb{X}_{s^-}\rangle ds \\
		&\qquad+\sum_{0<s\leq t}\sum_{\e_{I_1}\otimes \e_{I_2}=\e_I}\frac{1}{|I_2|!}1_{\{|I_2|> 1\}}\langle\epsilon_{I_1},\mathbb{X}_{s^-}\rangle \Delta X_s^{S(I_2)}1_{\{I_2\in\{1,\dots,N\}^{|I_2|}\}},
	\end{aligned}
	\end{equation}
	for a each multi-index $I$ with entries in $\{-1,0,1,\dots,N\}$.
	
	We are now interested in computing the expected signature of $\mathbb X$, i.e.~
$\mathbb E[\langle \epsilon_I,\mathbb X\rangle]$
for each multi-index $I$. One can prove (see e.g.~\cite{FS:17}),
that the expected signature of  an $\Rset^d$-valued Lévy process is finite whenever all the moments of the Lévy measure are finite.
Moreover, it can be computed analytically using the method in \cite{FS:17} or by recognizing that ${\mathbb{X}}^N$ is a polynomial process on $T^N(\Rset^d)$ for each $N\in\Nset$, and applying the so-called moment formula (see \cite{CKT:12, FL:20}).
As in \eqref{eqn3}, for each $\mathbbm a\in T((\mathbb R^d))$ such that $\pi_0(\mathbbm{a})=0$, we use the notation
\begin{align}\label{eq:exp}
    \exp({\mathbbm a}):=1+\sum_{k=1}^\infty\frac{\mathbbm{a}^{\otimes k}}{k!}.
\end{align}
The proof of the following result via polynomial technology can be found in Appendix~\ref{appendixB:expsiglevy}.
\begin{proposition}\label{expsiglevy}
	Fix $d\in \Nset$. Let $b^L\in\Rset^d$, $C^L\in (\R^d)^{\otimes 2}$ a positive definite $d\times d$ symmetric matrix and $F^L$ a $\sigma$-finite measure on $(\Rset^d,\mathcal{B}(\Rset^d))$ such that $F^L(\{0\})=0$ and 
	\begin{equation}\label{eqn1}
	\int_{\R^d}\|x\|^kF^L(dx)<\infty\quad \text{for all}\ k\geq 2.
	\end{equation}
	Let $L$ be an $\Rset^d$-valued Lévy process with triplet $(b^L,C^L,F^L)$ with respect to the ``truncation function'' $\chi(x)=x$. For each $N\in \mathbb N$, the truncated signature $\mathbb{L}^N$ of $L$ is a $T^N(\Rset^d)$-valued polynomial process. This in particular implies that all the moments of $\mathbb{L}^N$ exist and are finite.
	Moreover,
		\begin{equation}\label{eqn4}
	\Ex[\mathbb{L}_t]=	\exp\left(\left(b^L+\frac{C^L}{2}+\int_{\Rset^d}(\exp(x)-1-x) F^L(dx)\right)t\right).
	\end{equation}
\end{proposition}
\begin{remark}\label{rem2}
Let us restate expression \eqref{eqn4}  in a more explicit form. Observe that since
\[
	\int_{\Rset^d}(\exp(x)-1-x) F^L(dx)
	=\left(0,0,\frac{1}{2} \int_{\Rset^d}x^{\otimes 2}F^L(dx),
	\dots,\frac{1}{k!}\int_{\Rset^d}x^{\otimes k}F^L(dx),\dots \right),
	\]
\eqref{eqn4} can be written as $\Ex[\mathbb{L}_t]=	\exp(Qt)$ for
\begin{equation*}
    Q:=\left(0,b^L,\frac{1}{2} \left(C^L+\int_{\Rset^d}x^{\otimes 2}F^L(dx)\right),
	\dots,\frac{1}{k!}\int_{\Rset^d}x^{\otimes k}F^L(dx),\dots \right).
\end{equation*}
	Applying once again the definition of the exponential function \eqref{eq:exp}, we obtain that for $|I|>0$ it holds 
	\begin{equation*}
	    \Ex[\langle \epsilon_I,\mathbb L_t\rangle]=\sum_{k=1}^n\frac {t^k}{k! }\sum_{\epsilon_{I_1}\otimes\cdots\otimes \epsilon_{I_k}=\epsilon_I} \langle\epsilon_{I_1},Q\rangle\cdots\langle \epsilon_{I_k},Q\rangle,
	\end{equation*}
	where 
	$\langle \epsilon_\emptyset,Q\rangle=0$,\ 
$	\langle \epsilon_{j_1},Q\rangle=b_{j_1}^L$,
$	\langle \epsilon_{j_1,j_2},Q\rangle=\frac 1 2 (C_{j_1,j_2}^L+\int_{\mathbb R^d}x_{j_1}x_{j_2}F^L(dx))$, and 
$$\langle \epsilon_J,Q\rangle=\frac 1 {k! }\int_{\mathbb R^d}x_{j_1}\cdots x_{j_m}F^L(dx)$$ for each $J=(j_1,\ldots ,j_m)$ for $m\geq 3$. For example, setting $I=(1,3)$ we get that
$$\Ex[\langle \epsilon_{(1,3)},\mathbb L_t\rangle]
=t\langle\epsilon_{(1,3)},Q\rangle
+\frac {t^2}2
\langle\epsilon_{1},Q\rangle\langle \epsilon_{3},Q\rangle
=\frac t2\left(C_{1,3}^L+\int_{\mathbb R^d}x_1x_3 F^L(dx)\right)
+\frac {t^2}2b_1^Lb_3^L.$$
    \end{remark}
    \begin{remark}\label{rem3}
        Observe that following the argument of Lemma~A.2 in \cite{CGS:22}, we can use the i.i.d.  increments of a Lévy process to get
        $$\mathbb{E}[\langle e_{I}, {\mathbb{L}}_{s+t}\rangle|\mathcal F_s]
	=\sum_{e_{I_1}\otimes e_{I_2}=e_I}\langle e_{I_1}, {\mathbb{L}}_{s}\rangle
	\mathbb{E}[
	\langle e_{I_2}, {\mathbb{L}}_{t}\rangle].$$
        Using this representation, we can see that the formula obtained in Proposition~\ref{expsiglevy} also provides an explicit representation of the conditional moments of $\mathbb L$.
    \end{remark}
We are now ready to define a Lévy-type signature model, which will represent the price process of an arbitrary financial asset under an equivalent martingale measure $\mathbb{Q}$. 
\begin{definition}
	Let $X$ be the market's primary process introduced in Definition~\ref{def2} and $\mathbb{X}$ its signature. We define a \emph{Lévy-type signature model} as follows 
	\begin{equation}\label{Levymodel}
		S(\boldsymbol{\ell})_t=S_0+ \int_0^t\bigg(\sum_{|J|\leq n}\ell^J_W\langle \epsilon_{J},\mathbb{X}_{s^-}\rangle \bigg)dW_s+\int_0^t\int_\Rset \bigg(\sum_{|J|\leq n}\ell^J_\nu\langle \epsilon_{J},\mathbb{X}_{s^-}\rangle \bigg)x  \ (\mu-\nu)(ds,dx),
	\end{equation}
	where $ \ell^J_W,\ell^J_\nu\in \Rset$ and $J$ are multi-indices of length at most $n\in\Nset$ whose indices are in $\{-1,0,1,\dots,d\}$ with $d\leq N$. 
\end{definition}
From now on we always assume that
\begin{equation*}
nd+1\leq N.
\end{equation*}
This technical condition is needed in order to achieve an equivalent representation of the model, called sig-model representation, that will be derived in the next section. We refer to Remark~\ref{gooddefinitiontilda} for a precise discussion on this.

\begin{remark}\label{Levytyperemark}\phantomsection
	\begin{enumerate}
 \item The motivation for this type of model comes from Theorem \ref{th: UAT non antici} proved in Section \ref{sec:UAT}. Indeed, note that any continuous path functional of the enhanced market's primary underlying process can be approximated by linear functionals of its signature. The characteristics of $S(\boldsymbol{\ell})$ can thus be interpreted as approximations of continuous path functionals $F$ of $\mathbf{X}$, with $\mathbf{X}$ denoting the Marcus lift of $X$ as computed in Proposition \ref{prop2}.
		\item The no-arbitrage condition is satisfied. In fact, the process $S(\boldsymbol{\ell})$ is a square integrable martingale. This is a direct consequence of Proposition~\ref{Xislevy}, Remark~\ref{rem:finitmomentX}, Proposition~\ref{expsiglevy}, and the shuffle product formula according to which every polynomial on the signature admits a linear representation.
		\item The price process $S(\boldsymbol{\ell})$ is a \cadlag semimartingale whose characteristics, with respect to the ``truncation function'' $\chi(x)=x$, are given by
		\begin{align}\label{eq:characteristics}
		   B^S_t=0, \quad C^S_t=\int_0^t\bigg(\sum_{|J|\leq n}\ell^J_W\langle \epsilon_{J},\mathbb{X}_{s^-}\rangle \bigg)^2ds, \quad \nu^S(dt\times dx)=K_t(dx)\times dt,
		\end{align}
		where
		\begin{align*}
		    &K_t(A)=\int_\Rset 1_{A\setminus\{0\}}\bigg(\sum_{|J|\leq n}\ell^J_\nu\langle \epsilon_{J},\mathbb{X}_{t^-}\rangle \ x\bigg)F(dx).
		\end{align*}
	\end{enumerate}
\end{remark}

\subsection{Sig-model representation}
The purpose of this section is to achieve an alternative representation of a Lévy-type signature model, namely to express $S(\boldsymbol{\ell})$ itself as linear functional of the signature of the market's primary process. We will see that this new representation, that will be called sig-model representation, is more suitable for addressing the problem of pricing, and model calibration to both time-series data and option prices. Indeed, the sig-model representation is fully in line with the representation of \cite{CGS:22} in the simpler continuous case, so that all the arguments from there concerning tractability for pricing and calibration apply. We start by representing It\^o integrals as linear functionals of the signature of $X$.
\begin{definition}\label{def1}
For each multi-index $I$ and each $j\in\{-1,\ldots,d\}$
we consider the following multi-index transformation
	\begin{align*}
	&(\e_I;\e_j)^{\thicksim}
	:=\e_I\otimes \e_j-\frac{1}{2}\e_I\otimes\e_{-1}1_{\{i_{|I|}=j=0\}}
	+\sum_{\e_{I_1}\otimes \e_{I_{2}}=\e_I}
\alpha(|I_2|)(\e_{I_1}\otimes \e_{S(I_2)+j})1_{\{j>0\}}1_{\{I_2\in\{1,\ldots,d\}^{|I_2|}\}},
\end{align*}
where 
\begin{align*}
\alpha(|I_2|):=\sum_{k=1}^{|I_2|}(-1)^k\alpha(|I_2|,k),
\end{align*}
and 
denoting by $\mathcal I(r,k)$ the set of all multi-indices $J\in \Nset^k$ such that $S(J)=r$,  
$$\alpha(r,k):=\sum_{J\in \mathcal I(r,k)}\prod_{i=1}^k\frac 1 {(j_i+1)!}.$$ 
    \end{definition}

The proof of the next proposition can be found in Appendix~\ref{appendixB:index-transf}.
\begin{proposition}\label{index-transf}
	Fix $n\in\Nset$,  $I\in \{-1,0,1,\dots,d\}^n$, and $j\in\{-1,0,1,\dots,d\}$. Suppose that $N\geq 2$ and $nd+j\leq N$. Then,
\begin{align*}
		\int_{0}^{t}\langle \epsilon_{I},\mathbb{X}_{s^-}\rangle \  dX^{j}_s=& \ \langle (\epsilon_{I};\epsilon_{j})^{\thicksim},\mathbb{X}_t\rangle.
\end{align*}
\end{proposition}

   \begin{remark}\label{gooddefinitiontilda}
Fix $j\geq 1$ and observe that by Proposition~\ref{index-transf} it holds
	\begin{align*}
		\int_{0}^t\langle \epsilon_{I},\mathbb{X}_{s^-}\rangle dX^{j}_s
		&=
		\langle (\epsilon_{I};\epsilon_{j})^{\thicksim},\mathbb{X}_t\rangle\\
		&=\langle \epsilon_{I}\otimes \epsilon_{j},\mathbb{X}_{t}\rangle
		+\sum_{\e_{I_1}\otimes \e_{I_2}=\e_I}\alpha(|I_2|)\langle \epsilon_{I_1}\otimes\epsilon_{S(I_2)+j},\mathbb{X}_t\rangle1_{\{I_2\in\{1,\ldots,d\}^{|I_2|}\}},
\end{align*}
	for some coefficients $\alpha(|I_2|)\in \Rset$. Since 
\[\max_{|I_2|\leq n} S(I_2)+j=nd+j,
	\]	
we can conclude that $nd+j\leq N$ is a necessary condition for being able to write  $\int_0^t\langle \e_I,\mathbb X_{s-}\rangle dX_s^j$
as a linear functional of the signature of $X$ for each multi-index $|I|\leq n$.
\end{remark}

The sig-model representation is a direct application of Proposition~\ref{index-transf}. We thus omit the proof of the next corollary.
\begin{corollary}\label{cor1}
	The Lévy-type signature model defined in \eqref{Levymodel} admits the following sig-model representation 
	\begin{equation} \label{eq:sigrep}
S(\boldsymbol{\ell})_t=S_0+ \sum_{|J|\leq n}\big(\ell^J_W\langle (\epsilon_{J};\epsilon_0)^\thicksim,\mathbb{X}_{t}\rangle +\ell^J_\nu\langle (\epsilon_{J};\epsilon_{1})^\thicksim,\mathbb{X}_{t}\rangle\big).
	\end{equation}
\end{corollary}

\begin{remark}
\begin{enumerate}\phantomsection
\item Observe that thanks to the sig-model representation, $S(\boldsymbol{\ell})_t$ itself can be thought as an approximation of some continuous path functional of the enhanced market's underlying process. 

\item Moreover, it
allows to reduce the calibration to time-series data to a linear regression problem in the coefficients $\ell$ similarly as in Section 4.1 of \cite{CGS:22}.
Indeed, suppose that a trajectory of the market's primary process is observable, which  means that we can extract from time-series asset price data the trajectory of the ``market Brownian motion'' $W$ and the ``market Poisson random measure'' $\mu$, together with its compensator. Note here, that $X$ is specified under $\mathbb{Q}$, which means that the $\mathbb{Q}$-Brownian motion and the $\mathbb{Q}$-Poisson random measure already account for the market price of risk. In other words we extract by means of  quadratic variation estimators similarly as in Example~4.1 of \cite{CGS:22} a $\mathbb{P}$-Brownian with drift (corresponding to the $\mathbb{Q}$-Brownian motion)  and a $\mathbb{P}$-Poisson random measure  with potentially a different intensity ($\mathbb{Q}$-Poisson random measure).
Once $X$ is extracted, its signature  can be computed along the observed trajectory. Using a simple linear regression, the coefficients $\ell^J_W$
can then be matched with the market's continuous quadratic variation since the continuous quadratic variation of $S(\boldsymbol{\ell})$ is given by \eqref{eq:characteristics}. By means of \eqref{eq:sigrep},  the remaining parameters $\ell^J_{\nu}$ can then be estimated by performing a linear regression directly on the price data.

\item For pricing and calibrating to option data the sig-model representation also has a crucial advantage. Indeed, Monte-Carlo prices of any (possibly path-dependent) option on $S(\boldsymbol{\ell})$ can be easily obtained by computing all signature samples of $\mathbb{X}$ offline. Calibration to option data then reduces to a simple optimization task to infer $\ell^J_W$ and $\ell^J_{\nu}$ from market data. In contrast to classical models or when directly working with \eqref{Levymodel}, the sig-model represenation does not involve any Euler discretization in each optimization step, which is a crucial advantage. We refer to \cite{CGS:22} and \cite{CGMS:23} where this advantage has also been exploited.

\end{enumerate}

\end{remark}
\begin{example}\label{jumpsize1}
In Remark~\ref{gooddefinitiontilda} we have explained why in order to achieve the sig-model representation one has to include some moments of the jump-measure $\mu$ in the definition of the
primary process $X$ (see Definition~\ref{def2}). 
We now want to show that this necessary condition is no longer needed if the last component of $X$ consists of a standard Poisson process.

In this context $\nu$ satisfies
$$\nu(dt,dx)=dt\times \lambda \delta_1(dx),$$
for some $\lambda >0$ where $\delta_1$ denotes the Dirac measure.  Set then
	\begin{equation}\label{Xjump1}
		X_t:t\mapsto \left(t,W_t,\int_0^t\int_\Rset x\ \mu(ds,dx)\right),
\end{equation}
and observe that  for all $k\geq2$ and  $t\in[0,1]$ 
	\[
	\int_0^t\int_\Rset x^k\mu(ds,dx)
	=\int_0^t\int_\Rset x\ \mu(ds,dx)=\langle \e_{1},X_t\rangle
	.
	\]
Notice that to simplify the notation we do not compensate the sum of jumps in the jump-component of $X$.
As a consequence, for $j=1$ the multi-index transformation introduced in Definition~\ref{def1} simplifies to
$$
\begin{aligned}
	&(\e_I;\e_{1})^{\thicksim}:=\e_I\otimes \e_1
	+\sum_{\e_{I_1}\otimes \e_{I_{2}}=\e_I}
\alpha(|I_2|)(\e_{I_1}\otimes \e_{1})1_{\{I_2\in\{1,\ldots,1\}^{|I_2|}\}}. 
\end{aligned}
$$
Observe in particular that the linear combination describing $(\e_I;\e_{1})^{\thicksim}$ just involves indices with values in $\{-1,0,1\}$. The pairing with $\mathbb X_t$, where $X$ is as in \eqref{Xjump1}, is thus well defined. By Proposition~\ref{index-transf} we can conclude that
the Lévy-type signature model 
\begin{equation}\label{Sjump1}
		S(\boldsymbol{\ell})_t=S_0+ \int_0^t\left(\sum_{|J|\leq n}\ell^J_W\langle \epsilon_{J},\mathbb{X}_{s^-}\rangle \right)dW_s+\int_0^t\int_\Rset \left(\sum_{|J|\leq n}\ell^J_\nu\langle \epsilon_{J},\mathbb{X}_{s^-}\rangle \right)x  \ (\mu-\nu)(ds,dx),
\end{equation}
with $ \ell^J_W,\ell^J_\nu\in \Rset$ and $J$ multi-index with values in $\{-1,0,1\}$ can be rewritten as
\begin{align*} 
S(\boldsymbol{\ell})_t=S_0+ \sum_{|J|\leq n}&\ell^J_W\langle (\epsilon_{J};\epsilon_0)^{\thicksim},\mathbb{X}_{t}\rangle 
+ \sum_{|J|\leq n}\ell^J_\nu\langle (\epsilon_{J};\epsilon_{1})^{\thicksim}-\lambda(\epsilon_{J};\epsilon_{-1})^{\thicksim},\mathbb{X}_{t}\rangle.
\end{align*}

It is important to notice that the improvement of tractability achieved in this setting comes with a loss in generality. In particular, if a jump occurs at time $t$ its size has to be 
$$\sum_{|J|\leq n}\ell^J_\nu\langle \epsilon_{J},\mathbb{X}_{t^-}\rangle.$$
and is thus known in $t^-$.
\end{example}

\subsection{Pricing of sig-payoffs}\label{sec:pricing}
In this section we tackle the problem of pricing a contingent claim in Lévy-type signature models.  As an application of Proposition \ref{UAT1}, we consider contingent claims whose payoff is represented as a linear functional of the signature of the time-extended price process (extended price process henceforth), i.e.~so-called sig-payoffs, already analyzed in \cite{LNP:20, ASS:20, CGS:22} and the references therein. In order to give a clear treatment of the problem, we start by computing the signature of the time-extended price process
\begin{equation}\label{extendedpriceprocess}
	\widehat{S}(\boldsymbol{\ell})_t:=(t,S(\boldsymbol{\ell})_t),
\end{equation} whose components will be denoted by
\begin{align*}
	\widehat{S}(\boldsymbol{\ell})_t^{-1}=t,\qquad \widehat{S}(\boldsymbol{\ell})_t^1=S(\boldsymbol{\ell})_t.
\end{align*}
As $S(\boldsymbol{\ell})$ is given by \eqref{Levymodel},  the extended price process $\widehat{S}(\boldsymbol{\ell})$ is an $\Rset^2$-valued \cadlag semimartingale. Hence, its signature $\widehat{\mathbb{S}}(\boldsymbol{\ell})$ can be computed following Proposition~\ref{MarcusSDE1}, and thus
	\begin{equation}
\begin{aligned}\label{sigofS}
		\langle\epsilon_{I},\mathbb{\widehat{S}(\boldsymbol{\ell})}_t\rangle 
		=&\int_{0}^t \langle\epsilon_{I'},\mathbb{\widehat{S}(\boldsymbol{\ell})}_{s^-}\rangle d{\widehat{S}(\boldsymbol{\ell})}_s^{i_{|I|}}
		+1_{\{i_{|I|}=i_{{|I|}-1}=1\}}\frac{1}{2}\int_{0}^t\langle\epsilon_{I''},\mathbb{\widehat{S}(\boldsymbol{\ell})}_{s^-}\rangle\ d[{S}(\boldsymbol{\ell}),{S}(\boldsymbol{\ell})]_s^c \\
		&\qquad+\sum_{0<s\leq t}\sum_{\e_{I_1}\otimes \e_{I_2}=\e_I}\frac{1}{|I_2|!}1_{\{|I_2|> 1\}}\langle\epsilon_{I_1},\mathbb{\widehat{S}(\boldsymbol{\ell})}_{s^-}\rangle \Delta {{S}(\boldsymbol{\ell})}_s^{|I_2|}1_{\{I_2=(1,\ldots,1)\}},
	\end{aligned}
	\end{equation}
	for each $I$ with entries in $\{-1,1\}$.

Now, we show how to express the signature of $\widehat{S}(\boldsymbol{\ell})$ in terms of the signature of the underlying process $X$. This result is again based on the multi-index transformation introduced in Definition~\ref{def1}. Before making  the argument precise, let us make the following observations and introduce some useful notation.

Consider the price process defined as in equation \eqref{Levymodel}. Let $M:=\sum_{i=0}^n (d+2)^i$ be the number of possible multi-indices with length at most $n$ and indices in $\{-1,0,1,\dots,d\}$. Observe that we are considering also the empty index $\emptyset$. Define 
\[
\{\epsilon_{J^h}; 1\leq h\leq M\},
\]
to be an ordered version of $\{\epsilon_J, |J|\leq n\}$. For instance, we may consider the order with respect to the length first and alphabetically next. Then,
\begin{equation*}
	S(\boldsymbol{\ell})_t=S_0+\sum_{h=1}^M\ell^{W}_h\langle (\epsilon_{J^h};\epsilon_0)^{\thicksim},\mathbb{X}_{t}\rangle + \sum_{h=1}^M\ell^{\mu}_h\langle (\epsilon_{J^h};\epsilon_1)^{\thicksim},\mathbb{X}_{t}\rangle.
\end{equation*}
For $\alpha=(\alpha^1,\alpha^2,\dots,\alpha^M)\in \Nset^M$ we write
\begin{align*}
	\epsilon^{\shuffle \alpha}:=\frac{1}{\alpha!}(\shuffle^{\alpha^1}\epsilon_{J^1})\shuffle (\shuffle^{\alpha^2}\epsilon_{J^2})\shuffle \dots \shuffle (\shuffle^{\alpha^M}\epsilon_{J^M}),
\end{align*}
and for 
$
\ell_{1,\dots,M}:=(\ell_{1},\ell_2,\dots ,\ell_M)\in \Rset^M,
$ we also set
$\ell_{1,\dots,M}^\alpha:=\ell_1^{\alpha^1}\cdots\ell_M^{\alpha^M}.$

We are now ready to state the announced result. It will become clear later that this is the key result to solve analytically the pricing problem for sig-payoffs. The proof can be found in Appendix~\ref{appendixB:fromStoX}, where a discussion on the necessity of the lower bound on $N$ is also provided.
\begin{theorem}\label{fromStoX}
	Let $\widehat{S}(\boldsymbol{\ell})$ be the time-extended price process defined in \eqref{extendedpriceprocess} and $\widehat{\mathbb{S}}(\boldsymbol{\ell})$ its signature. Let $I$ be a multi-index with indices in $\{-1,1\}$. If $N\geq |I|(dn+1)$, then there exists a linear combination of multi-indices $U_I(\boldsymbol{\ell})$ with indices in $\{-1,0,1,\dots,N\}$, such that the following equality holds
	\begin{equation}\label{multSX}
	    \langle\e_I,\widehat{\mathbb{S}}(\boldsymbol{\ell})_t\rangle=\langle \epsilon_{U_I(\boldsymbol{\ell})},\mathbb{X}_t\rangle.
	\end{equation}
	Furthermore, 
	$
	\langle \epsilon_{U_I(\boldsymbol{\ell})},\mathbb{X}_t\rangle
	$
	is a polynomial of degree $|I|$ in $(\ell^W_{{1,\dots,M}},\ell^\nu_{{1,\dots,M}})\in \Rset^{2M}$, and  $U_I(\boldsymbol{\ell})$ can be computed recursively on the length $|I|$ of $I$ setting $\e_{U_\emptyset(\boldsymbol{\ell})}=\e_{\emptyset}$ and
	\begin{align*}
		{\e_{U_I(\boldsymbol{\ell})}}=&(\e_{{U_{I'}(\boldsymbol{\ell})}};{\e_{-1}})^{\thicksim}1_{\{i_{|I|}=-1\}}+\sum_{\alpha\in \Nset^M, \ S(\alpha)=1}(\ell^W_{{1,\dots,M}})^{\alpha} (\e_{U_{I'}(\boldsymbol{\ell})}\shuffle {\epsilon^{\shuffle \alpha}};\e_0)^{\thicksim}1_{\{i_{|I|}=1\}}\\ 
		&+\sum_{\alpha\in \Nset^M, \ S(\alpha)=1}(\ell^\nu_{{1,\dots,M}})^{\alpha} (\e_{U_{I'}(\boldsymbol{\ell})}\shuffle {\epsilon^{\shuffle \alpha}};\e_1)^{\thicksim}1_{\{i_{|I|}=1\}}\\
		&+\sum_{\alpha\in \Nset^M, \ S(\alpha)=2}(\ell^W_{{1,\dots,M}})^{\alpha} (\e_{U_{I''}(\boldsymbol{\ell})}\shuffle \epsilon^{\shuffle \alpha};\e_{-1})^{\thicksim}1_{\{|I|>1\}}1_{\{(i_{|I|-1},i_{|I|})=(1,1)\}} \\ 
		&+\sum_{\e_{I_1}\otimes \e_{I_2}=\e_I} \sum_{\alpha\in \Nset^M, \ S(\alpha)=|I_2|}(\ell^\nu_{{1,\dots,M}})^{\alpha} (\e_{U_{I_1}(\boldsymbol{\ell})}\shuffle \epsilon^{\shuffle \alpha};\e_{|I_2|})^{\thicksim}1_{\{|I_2|>1\}}1_{\{I_2=(1,\dots,1)\}}.
	\end{align*}
\end{theorem}
\begin{remark}
    Consider the setting of Example \ref{jumpsize1}. From an inspection of the proof of Theorem~\ref{fromStoX}, it can be deduced that even in this case we do not need to consider higher moments of the jump measure in the definition of the market's primary process. This means that for every multi-index $I$ with indices in $\{-1,1\}$, there always exists a linear combination of multi-indices $U_I(\boldsymbol{\ell})$ with entries in $\{-1,0,1\}$, such that 
	\begin{equation}\label{fromSjump1toXjump1}
	    \langle\e_I,\widehat{\mathbb{S}}(\boldsymbol{\ell})_t\rangle=\langle \epsilon_{U_I(\boldsymbol{\ell})},\mathbb{X}_t\rangle.
	\end{equation}
	Here $S(\boldsymbol{\ell})$ and $X$ denote the processes defined in equation \eqref{Sjump1} and  \eqref{Xjump1} respectively. 
	More precisely, the representation in \eqref{fromSjump1toXjump1} holds with $(\e_I;\e_1)^{\thicksim}$ replaced by
	$(\e_I;\e_1)^{\thicksim}-\lambda (\e_I;\e_{-1})^{\thicksim}$ and $(\e_I;\e_k)^{\thicksim}$ replaced by $(\e_I;\e_1)^{\thicksim}$ for each $k\geq 2$.
\end{remark}

Finally, we recall the notion of sig-payoffs as introduced in \cite{LNP:20} and address the problem of computing their arbitrage free prices.
\begin{definition}\label{sigpayoffdef}
	Let $\widehat{S}(\boldsymbol{\ell})$ be the time-extended price process defined in \eqref{extendedpriceprocess} and $\widehat{\mathbb{S}}(\boldsymbol{\ell})$ its signature. We define a \emph{sig-payoff} with maturity $T>0$ as the payoff $C^S$ that pays to the holder of the claim an amount equal to 
	\begin{align} \label{eq:HS}
	C^S=\sum_{|I|\leq m}h^I\langle \epsilon_{I},\widehat{\mathbb{S}}(\boldsymbol{\ell})_T\rangle,
	\end{align}
	where, $h^I\in \Rset$ and $I$ are multi-indices of length at most $m\in \Nset$ with indices in $\{-1,1\}$.
\end{definition}

Observe that by Theorem~\ref{fromStoX}, Proposition~\ref{Xislevy}, and Proposition~\ref{expsiglevy}  the random variable $C^S$ given in \eqref{eq:HS} is  square integrable.

Applying Theorem~\ref{fromStoX}, we are now ready to provide a pricing formula for sig-payoffs.

\begin{corollary}
	Let 
	\[
	C^S=\sum_{|I|\leq m}h^I\langle \epsilon_{I},\widehat{\mathbb{S}}(\boldsymbol{\ell})_T\rangle
	\]
	be a sig-payoff with maturity $T>0$. If $N\geq m(nd+1)$, then, its price can be expressed by
	\begin{equation}\label{price}
		\sum_{|I|\leq m}h^I\Ex_{\mathbb{Q}}[\langle \epsilon_{I},\widehat{\mathbb{S}}(\boldsymbol{\ell})_T\rangle]=\sum_{|I|\leq m}h^I\Ex_{\mathbb{Q}}[\langle \epsilon_{U_I(\boldsymbol{\ell})},\mathbb{X}_T\rangle],
	\end{equation}
	where $U_I(\boldsymbol{\ell})$ is as in Theorem~\ref{fromStoX}. This in particular implies that 
	\eqref{price} is a polynomial of degree at most $m$ in  $(\ell^W_{{1,\dots,M}},\ell^\nu_{{1,\dots,M}})\in \Rset^{2M}$. 
\end{corollary}
\begin{proof}
	The results follows from Theorem~\ref{fromStoX}.
\end{proof}
\begin{remark}
	Observe that, thanks to Proposition~\ref{expsiglevy}, the expectations $\Ex_{\mathbb{Q}}[\langle \epsilon_{U_I(\boldsymbol{\ell})},\mathbb{X}_T\rangle]$, for $|I|\leq m$, can be computed analytically. Hence, after their computation, the calibration of the model to option prices corresponding to sig-payoffs reduces to a polynomial optimization problem in $(\ell^W_{{1,\dots,M}},\ell^\nu_{{1,\dots,M}})\in \Rset^{2M}$. For usual calibration to call and put options the above described Monte-Carlo approach can be combined with variance reduction techniques  based on \eqref{price}. For more details about these topics in the continuous setting see e.g.~\cite{LNP:20} and \cite{CGS:22}.
\end{remark}
\subsection{Hedging of sig-payoffs}\label{sec:hedging}
In the Lévy-type signature model the market is not complete. 
One risk management criterion is 
to minimize the hedging error in a mean square sense, i.e.~to use a quadratic loss function.
In this context losses and gains are treated in a symmetric manner and the criterion to be minimized can be either the squared hedging error at maturity or measured locally in time (see e.g.~\cite{SCHWEIZER:98,SC:90, SCW:95, S:99}, \cite{Pham:00} and the references therein for an extensive analysis  of these problems).
In what follows, we focus on the first approach which leads to the construction of the so-called mean variance hedging strategy. It specifically consists in solving the optimization problem
			\begin{align}\label{hproblem}
			    	\Ex\left[\left(C-v-\int_0^T\theta_s\  dS_s\right)^2\right].
			 \end{align}
over all initial endowments $v\in\Rset$ and all hedging strategies $\theta$.
Here, $C$ denotes the payoff at time $T$ of a European contingent claim and $S$ the price process of the considered risky asset.
When the price process $S$ is a martingale (thus working under a risk neutral measure $\mathbb{Q}$), this problem has been explicitly solved in \cite{S:85}. 
In particular, the authors showed that, given any contingent claim with maturity $T>0$ and square integrable payoff $C$, there exists a unique solution $\varphi^*=(v^*,\theta^*)$ to the optimization problem \eqref{hproblem} which is computed as follows. Let $V$ denote the martingale generated by $C$, that is $V_t:=\Ex[C|\mathcal{F}_t]$, $t\leq T$. Then, the optimal initial endowment and the optimal hedging strategy, 
are given by
\begin{equation}\label{generalstrategy}
v^*=V_0,	 \qquad  \qquad \theta^*=\frac{d [ V,S]^{pred}}{d [S]^{pred}},
\end{equation}
where $[\,\cdot\,]^{pred}$ and  $[\cdot,\cdot]^{pred}$ stand for predictable quadratic variation and covariation, respectively.

We now adapt these  results to our setting, working under $\mathbb{Q}$ where the price process is a martingale. For a discussion why measuring the risk under $\mathbb{Q}$ can be reasonable we refer to \cite{A:10}. The proof of the subsequent theorem can be found in Appendix~\ref{appendixB:localriskminimization}.
\begin{theorem}\label{localriskminimization}
Fix a multi-index $I$ and let 
$
C^S=\langle \epsilon_{I},\widehat{\mathbb{S}}(\boldsymbol{\ell})_T\rangle
$
be a sig-payoff with maturity $T>0$. The unique solution $\varphi^*=(v^*,\theta^*)$ to the optimization problem \eqref{hproblem} with respect to $C^S$ is given by
\begin{align*}
v^*&=\Ex[\langle \epsilon_{I},\widehat{\mathbb{S}}(\boldsymbol{\ell})_T\rangle]
   \\
  \theta_t^*&=\frac
  {
  \sum_{\e_{J_1}\otimes \e_{J_2}\otimes \e_{J_3}=\e_{U_I(\boldsymbol \ell)}}
    \sum_{|H|\leq n}
      \langle\epsilon_{J_1\shuffle H},\mathbb{X}_{t^-}\rangle
      \gamma_2(J_2,H,\boldsymbol{\ell})
      \Ex[\langle \e_{J_3},\mathbb X_{T-t}\rangle] }
      {
      \Big(\sum_{|H|\leq n}\ell^J_W\langle \epsilon_{H},\mathbb{X}_{t^-}\rangle \Big)^2+\Big(\sum_{|H|\leq n}\ell^H_\nu\langle \epsilon_{H},\mathbb{X}_{t^-}\rangle \Big)^2\int_\Rset x^2F(dx)
      },
\end{align*}
where $\gamma_2(J_2,H,\boldsymbol\ell)=\gamma_W(J_2)\ell^H_W
      +
    \gamma_\nu(J_2)\ell^H_\nu\int_{\mathbb R} x^{S(J_2)+1} F(dx)$ for
$$\gamma_W(J_2)=1_{\{J_2=(0)\}}\quad\text{and}\quad\gamma_\nu(J_2)=\frac{1}{|J_2|!}1_{\{|J_2|\geq 1\}}1_{\{J_2\in\{1,\dots,N\}^{|J_2|}\}}.$$

\end{theorem}

\begin{remark}\label{rem4} We report here some further observations about Theorem~\ref{localriskminimization}.
\begin{enumerate}
\item By linearity the results of Theorem~\ref{localriskminimization} extend to $C^S=\sum_{|I|\leq m}h^I\langle \epsilon_{I},\widehat{\mathbb{S}}(\boldsymbol{\ell})_T\rangle$.
\item  Remark~\ref{rem2} provides explicit formulas for $\Ex[\langle \e_{J},\mathbb X_{T}\rangle]$ and, as a consequence of Theorem~\ref{fromStoX}, also for $v^*=\Ex[\langle \e_{U_I(\boldsymbol \ell)},\mathbb{{X}}_T\rangle]$.

\item
The above result gives a fully explicit formula for the mean variance strategy that hinges solely on the option and the model's parameters as well as the signature of the market's primary process $X$. 
\item\label{itiv} If $i_{|I|-1}=-1$ and $i_{|I|}=1$ by \eqref{sigofS} we have that
$\langle\epsilon_{I},\mathbb{\widehat{S}(\boldsymbol{\ell})}_t\rangle 
		=\int_{0}^t \langle\epsilon_{I'},\mathbb{\widehat{S}(\boldsymbol{\ell})}_{s^-}\rangle d{{S}(\boldsymbol{\ell})}_s,$
		which is a true martingale. This yields $V_t=\langle\epsilon_{I},\mathbb{\widehat{S}(\boldsymbol{\ell})}_t\rangle $ and hence
		$$
			\theta^*_t=\frac{d[V,{{S}(\boldsymbol{\ell})}]^{pred}_t}{d [ {{S}(\boldsymbol{\ell})}]^{pred}_t} =\frac{\langle\epsilon_{I'},\mathbb{\widehat{S}(\boldsymbol{\ell})}_{t^{-}}\rangle d [S(\boldsymbol{\ell})]^{pred}_t}{d [S(\boldsymbol{\ell})]^{pred}_t}
			=\langle\epsilon_{I'},\mathbb{\widehat{S}(\boldsymbol{\ell})}_{t^{-}}\rangle.
			$$
\item A contingent claim $H^S$ is attainable if it almost surely admits the representation
$$
	H^S=\Ex[H^S]+\int_0^T \xi^*_s dS(\boldsymbol{\ell})_s.
$$
In the case of a sig-payoff
$
H^S=\sum_{|I|\leq m}h^I\langle \epsilon_{I},\widehat{\mathbb{S}}(\boldsymbol{\ell})_T\rangle
$
this is for instance the case if for each $I$ such that $h^I\neq0$ it holds $i_{|I|-1}=-1$ and $i_{|I|}=1$. In this case we indeed have that
\begin{align*}
	H^S=&h^{\emptyset}+\sum_{1\leq |I|\leq m } h^I\int_0^T\langle \epsilon_{I'},\widehat{\mathbb{S}}(\boldsymbol{\ell})_{s^-}\rangle dS(\boldsymbol{\ell})_s,
\end{align*}
and thus in particular $\Ex[H^S]=h_\emptyset$ and $\xi^*_s=\sum_{1\leq |I|\leq m } h^I\langle \epsilon_{I'},\widehat{\mathbb{S}}(\boldsymbol{\ell})_{s^-}\rangle$.
\end{enumerate}
\end{remark}

\subsection{Measure transformations for L\'evy-type signature models }\label{sec:change}

In the previous sections, we worked with Lévy-type signature models under some martingale measure $\mathbb{Q}$. Precisely, we started with a market's primary process $X$ and its signature, and then  defined a Lévy-type signature model as the process obtained by integrating, in \Ito sense, linear functionals of the signature of $X$ with respect to its martingale components.  

The goal of this section is to specify the dynamics of such a Lévy-type signature model under the real world measure $\mathbb{P}$. In particular, we will investigate conditions guaranteeing that a Lévy-type signature model under $\mathbb Q$ remains a Lévy-type signature model also under $\mathbb P$, meaning that its characteristics are expressed in terms of the signature of a primary $\mathbb{P}$-Lévy process whose components include the $\mathbb{P}$-driving signals. As a result, explicit formulas for the computation of the moments of the price process also under the measure $\mathbb{P}$ can be provided. This can then be exploited for instance in view of generalized methods of moments for parameter estimation under $\mathbb{P}$.

First of all, recall that the equivalence of the measures $\mathbb{Q}$ and $\mathbb{P}$ implies the existence of a positive $\mathbb{Q}$-martingale $Z$, with $\Ex_\mathbb{Q}[Z]=1$, such that 
\begin{equation}\label{eqn17}
 d\mathbb{P}=Z\, d\mathbb{Q}.   
\end{equation}
Moreover, since we are dealing with the augmented filtration generated by $W$ and $\mu$, the martingale $Z$ is the solution of the following SDE
\begin{align}\label{exponentialmartingale}
	&dZ_t=Z_{t^-}\left(f(t)dW^\mathbb{Q}_t+\int_\Rset (e^{g(t,x)}-1)(\mu-\nu^\mathbb{Q})(dt,dx)\right),\qquad Z_0=1, 
\end{align}
where $f(t)$ and $g(t,y)$ are predictable processes such that $Z$ is a true martingale (see e.g.~Theorem~2.1 in \cite{K:04}). We here indicate the $\mathbb{Q}$ dependence of the driving signals $W^{\mathbb{Q}}$ and $\mu-\nu^\mathbb{Q}$.  In such a case, we have that (see e.g.~Theorem~2.3 in \cite{K:04})
\begin{equation*}
	W^\mathbb{P}=W^\mathbb{Q}-\int_{0}^{\cdot}f(s)ds
\end{equation*}
is a standard $\mathbb{P}$-Brownian motion and the compensator of the jump measure $\mu$ is
\begin{equation*}
	\nu^\mathbb{P}(dt,dx)=e^{g(t,x)}\nu^\mathbb{Q}(dt,dx)=e^{g (t,x)}F(dx)dt.
\end{equation*}
Moreover, the dynamics under $\P$ of the $\Q$-Lévy model 
defined in \eqref{Levymodel} are given by
\begin{align}\label{SunderP}
dS(\boldsymbol{\ell})_t=&\bigg(\sum_{|J|\leq n}\ell^J_W\langle \epsilon_{J},\mathbb{X}_{t^-}\rangle f(t)+\sum_{|J|\leq n}\ell^J_\nu\langle \epsilon_{J},\mathbb{X}_{t^-}\rangle  \int_\Rset y(e^{g(t,x)}-1)F(dx)\bigg)dt\nonumber \\
+& \bigg(\sum_{|J|\leq n}\ell^J_W\langle \epsilon_{J},\mathbb{X}_{t^-}\rangle \bigg)dW^\mathbb{P}_t\\
+& \int_\Rset \bigg(\sum_{|J|\leq n}\ell^J_\nu\langle \epsilon_{J},\mathbb{X}_{t^-}\rangle \bigg) x  \ (\mu-\nu^\mathbb{P})(dt,dx)\nonumber. 
\end{align}
Therefore, its characteristics are expressed in terms of linear functionals of the signature of the process $X$, which is in general not a $\mathbb{P}$-Lévy process and whose components do not include the $\mathbb{P}$-driving signals. 

In order to recover some tractability also under $\P$, we are interested in Lévy-type signature models which satisfy \eqref{multSX} also for some $\P$-Lévy process $Y$ which fulfills the hypothesis of Proposition~\ref{expsiglevy}. 
\begin{definition}
A process $S(\boldsymbol{\ell})$ is called $\P$-$\Q$-Lévy-type signature model if it satisfies \eqref{Levymodel} and admits the representation
\begin{align*}
    		S(\boldsymbol{\ell})_t
    		=S_0&+ 
    		\int_0^t\bigg(\sum_{|J|\leq n}\bar\ell^J_t\langle \epsilon_{J},\mathbb{Y}_{s^-}\rangle \bigg)ds+
    		\int_0^t\bigg(\sum_{|J|\leq n}\bar\ell^J_W\langle \epsilon_{J},\mathbb{Y}_{s^-}\rangle \bigg)d\overline W_s\\
    		&+
    		\int_0^t\int_\Rset \bigg(\sum_{|J|\leq n}\bar\ell^J_\nu\langle \epsilon_{J},\mathbb{Y}_{s^-}\rangle \bigg)x  \ (\overline\mu-\overline\nu)(ds,dx),
\end{align*}
for some $\bar\ell^J_t, \bar\ell^J_W,\bar\ell^J_\nu\in \Rset$ and some Lévy process $Y$ of the form
$$
    Y_t:=  \left(t,\overline W_t,\int_0^t\int_\Rset x\ (\overline\mu-\overline\nu)(ds,dx), \int_0^t\int_\Rset x^2\ \overline\mu(ds,dx),\dots,\int_0^t\int_\Rset x^N\ \overline\mu(ds,dx)\right),
$$
for some $N$ large enough, $\P$-Brownian motion $\overline W$, and homogeneous $\P$-Poisson random measure $\overline \mu$ with compensator $\overline \nu$, such that $Y$
 has finite moments of all orders.
\end{definition}

The proof of the following proposition can be found in Appendix \ref{proof: prop1}

\begin{proposition}\label{prop1}
    Let $S(\boldsymbol{\ell})$ be a Lévy-type signature model under $\Q$ and fix $\P$ according to \eqref{eqn17} and \eqref{exponentialmartingale}, for $g(t,x)=g(x)$ a deterministic function and
    \begin{equation}\label{eqn9}
		f(t)=\sum_{|I|\leq p}f^{I}\langle \epsilon_{I},\mathbb{X}_{t^-}\rangle,
	\end{equation}
	with $f^I\in \Rset$, $p\in \Nset$, $I$ multi-indices with indices in the set $\{-1,0,1,2,\dots,N\}\setminus\{0,1\}$. If 
	\begin{align}\label{conditionmart}
	    \int_{\Rset}(e^{g(x)}-1)^2F(dx)<\infty \quad \text{ and }\quad \Ex_\Q\bigg[\exp\Big(\frac{1}{2}\int_0^Tf(s)^2ds\Big)\bigg]<\infty,
	\end{align}
	then $S(\boldsymbol{\ell})$ is a $\P$-$\Q$-Lévy-type signature model with respect to the process $Y$ given by
\begin{equation}\label{eq:Y}
Y_t:=  \left(t,W_t^\P,\int_0^t\int_\Rset x\ (\mu-\nu^\P)(ds,dx), \int_0^t\int_\Rset x^2\ \mu(ds,dx),\dots,\int_0^t\int_\Rset x^N\ \mu(ds,dx)\right).
\end{equation}
 \end{proposition}
    \begin{remark}
         Since tractability under $\Q$ is a key feature for the results of the paper, the process $S(\boldsymbol{\ell})$ in Proposition~\ref{prop1} is first specified through its dynamics under $\Q$. Applying the above  measure change procedure in the reverse direction one can however also start with a process first specified under $\P$.
    \end{remark}

\appendix

\section{Proofs of Section \ref{sec:UAT}} 
For notational convenience, in the proofs of Section \ref{sec:UAT} we focus only on the set of time-extended w.g.~\cadlag $p$-rough paths, for $p\in[2,3)$. However, as pointed out in Remark \ref{rem:proof UAT1} \ref{rem:proof UAT1 i} below, the exact same reasoning can also be applied when working with time-extended w.g.~\cadlag $p$-rough paths, for $p\in[1,2)$.
\subsection{Proof of Proposition \ref{UAT1}}\label{appendixA:UAT1}
	This result follows by the Stone-Weierstrass theorem applied to the set $A$ given by
	\begin{align}\label{setA}
		A:=\text{span}\{l:K\rightarrow\Rset \ ; \ \widehat{\mathbf{X}}\mapsto\langle \epsilon_I,\widehat{\mathbb{X}}_1\rangle\colon |I|\geq0\}.
	\end{align}
	Therefore, we must prove that $A$ satisfies the following conditions: \begin{enumerate}
	    \item\label{iti} it is a linear subspace of the space of continuous functions from $K$ to $\Rset$,
	    \item\label{itii}it is a sub-algebra and contains a
	non-zero constant function,
	    \item\label{itiii}it separates points.
	\end{enumerate}
	\ref{iti}: Let $N\geq3$ and denote by $\widehat{\mathbb{X}}^{N}$ the solution of the Marcus-type RDE given by \eqref{Rdemarcus}. By Proposition~3.18(ii) in \cite{CF:19}, the solution map of equation \eqref{Rdemarcus} 
	\begin{align}\label{solutionmap}
(K,\sigma_\infty)&\rightarrow (D([0,1],T^N(\Rset^{d+1})),\sigma_\infty)\\
	\widehat{\mathbf{X}}&\mapsto\widehat{\mathbb{X}}^N\nonumber 
	\end{align}
	is continuous on $K$ for every $N\geq 3$, where the $J_1$-metric on $D([0,1],T^N(\Rset^{d+1}))$ is computed via the distance $\rho$ induced by the norm \eqref{eq:normTN}. As already explained in  Remark~\ref{rem:solution}, equation \eqref{Rdemarcus} is in fact an example of a Marcus-type RDE \eqref{eq: generalMarcus} with linear vector fields. Therefore, by applying the result stated in Theorem~10.53 of \cite{FV:10}, one can see that all the assumptions of Proposition~3.18 in \cite{CF:19} are satisfied. 
	
	Next, observe that for a generic metric space $(E,d)$, the evaluation map
	\begin{align*}
		(D([0,1],E), &\ \sigma_\infty)\rightarrow (E,d)\nonumber \\
		&X\mapsto X_1
	\end{align*}
	is continuous. More precisely, let $X\in D([0,1],E)$ and let $(X^n)_n\subset D([0,1],E)$ be a sequence such that $X^n\rightarrow X$ as $n\rightarrow\infty$ with respect to the $J_1$-topology. Then, there exists $(\lambda^n)_n\subset \Lambda$ such that $|\lambda^n|\rightarrow0$ and $\sup_{s\in[0,1]}d(X^n_{\lambda^n(s)},X_s)\rightarrow0$ as $n\rightarrow\infty$ (see e.g.~Chapter 3 in \cite{billing}), and in particular  $d(X^n_{\lambda^n(1)},X_1)=d(X^n_1,X_1)\rightarrow0$, as $n\rightarrow\infty$.
	
	In the present setting this result yields that  continuity of \eqref{solutionmap} implies continuity of
	\begin{align*}
	 (K,\sigma_\infty)&\rightarrow (T^N(\Rset^{d+1}),\rho)\\
	\widehat{\mathbf{X}}&\mapsto\widehat{\mathbb{X}}_1^N.\nonumber 
	\end{align*}
 Since the map $\widehat{\mathbb{X}}_1^N\mapsto\langle \epsilon_I,\widehat{\mathbb{X}}_1^N\rangle$
	is  continuous for each $I$, the claim follows.

\noindent\ref{itii}: Since $\widehat{\mathbb{X}}_1$ is a group-like element, the shuffle property holds. Therefore, $A$ is a sub-algebra. Moreover, it also contains the 
	non-zero constant function
	$
	\widehat{\mathbf{X}}\mapsto\langle \epsilon_{\emptyset},\widehat{\mathbb{X}}_1\rangle\equiv1.
	$

	\noindent\ref{itiii}: As a final step, we need to show that the set $A$ separates points.
	Let us consider $\widehat{\mathbf{X}},\widehat{\mathbf{Y}}\in  K$, with  $\widehat{\mathbf{X}}\neq \widehat{\mathbf{Y}}$. We prove that there exists $n\in \Nset$ and $i,k,j\in\{-1,1,\dots,d\}$ such that 
	\begin{align*}
	    &\langle (\epsilon_i\shuffle \epsilon_{-1}^{\otimes n})\otimes \epsilon_{-1}+((\epsilon_j\otimes \epsilon_k)\shuffle \epsilon_{-1}^{\otimes n})\otimes \epsilon_{-1},\widehat{\mathbb{X}}_1\rangle\\
	    &\qquad\neq\langle (\epsilon_i\shuffle \epsilon_{-1}^{\otimes n})\otimes \epsilon_{-1}+((\epsilon_j\otimes \epsilon_k)\shuffle \epsilon_{-1}^{\otimes n})\otimes \epsilon_{-1},\widehat{\mathbb{Y}}_1\rangle.
	\end{align*}
	We reason by contradiction. Assume that for all $n\in \Nset$ and all $i,k,j\in\{-1,1,\dots,d\}$ 
	\begin{align}\label{finaleq}
	    &\langle (\epsilon_i\shuffle \epsilon_{-1}^{\otimes n})\otimes \epsilon_{-1}+((\epsilon_j\otimes \epsilon_k)\shuffle \epsilon_{-1}^{\otimes n})\otimes \epsilon_{-1},\widehat{\mathbb{X}}_1\rangle\\
	    &\qquad=\langle (\epsilon_i\shuffle\nonumber \epsilon_{-1}^{\otimes n})\otimes \epsilon_{-1}+((\epsilon_j\otimes \epsilon_k)\shuffle \epsilon_{-1}^{\otimes n})\otimes \epsilon_{-1},\widehat{\mathbb{Y}}_1\rangle.
	\end{align}
	In particular, for $k=j=-1$, by Proposition~\ref{prop:rough=Young}, for all $n\in\Nset$ and $i\in\{-1,1,\dots,d\}$ we get
	\begin{align*}
		\langle (\epsilon_i\shuffle \epsilon_{-1}^{\otimes n})\otimes \epsilon_{-1},\widehat{\mathbb{X}}_1\rangle=\langle (\epsilon_i\shuffle \epsilon_{-1}^{\otimes n})\otimes \epsilon_{-1},\widehat{\mathbb{Y}}_1\rangle.
	\end{align*}
	By Proposition~\ref{prop:rough=Young},
	$
	\langle (\epsilon_i\shuffle \epsilon_{-1}^{\otimes n})\otimes \epsilon_{-1},\widehat{\mathbb{X}}_1\rangle=\int_0^1\langle \e_{i},\widehat{\mathbf{X}}_{s-}\rangle \frac{s^{n}}{n!}ds,
	$
	and
	$
	\langle (\epsilon_i\shuffle \epsilon_{-1}^{\otimes n})\otimes \epsilon_{-1},\widehat{\mathbb{Y}}_1\rangle=\int_0^1\langle \e_{i},\widehat{\mathbf{Y}}_{s-}\rangle \frac{s^{n}}{n!}ds
	.$ By Corollary 4.24 in \cite{Brezis} and the \caglad property of the paths $$s\mapsto \langle \e_{i},\widehat{\mathbf{X}}_{s-}\rangle, \qquad s\mapsto \langle \e_{i},\widehat{\mathbf{Y}}_{s-}\rangle,$$ both starting in $0$ by Assumption \ref{assumption}, we then deduce that $\langle \e_{i},\widehat{\mathbf{X}}\rangle =\langle \e_{i},\widehat{\mathbf{Y}}\rangle $.
	
Similarly, choosing $i=-1$, equation \eqref{finaleq} yields that
	\begin{align*}
		\langle (\epsilon_k\otimes \epsilon_j)\shuffle \epsilon_{-1}^{\otimes n})\otimes \epsilon_{-1},\widehat{\mathbb{X}}_1\rangle=\langle (\epsilon_k\otimes \epsilon_j)\shuffle\nonumber  \epsilon_{-1}^{\otimes n})\otimes \epsilon_{-1},\widehat{\mathbb{Y}}_1\rangle
	\end{align*}
	for all $n\in \Nset$ and all $j,k\in\{-1,1,\dots,d\}$. Once again, by Proposition~\ref{prop:rough=Young}, 
	\[
	\langle (\epsilon_k\otimes \epsilon_j)\shuffle \epsilon_{-1}^{\otimes n})\otimes \epsilon_{-1},\widehat{\mathbb{X}}_1\rangle=\int_0^1  \langle\e_{(k,j)},\widehat{\mathbf{X}}_{s^-}\rangle\frac{s^{n}}{n!}ds,
	\]
	and
	\[
	\langle (\epsilon_k\otimes \epsilon_j)\shuffle \epsilon_{-1}^{\otimes n})\otimes \epsilon_{-1},\widehat{\mathbb{Y}}_1\rangle=\int_0^1\langle\e_{(k,j)},\widehat{\mathbf{Y}}_{s^-}\rangle\frac{s^{n}}{n!}ds.
	\]
The same reasoning as above yields then that $\langle\e_{(k,j)},\widehat{\mathbf{X}}\rangle=\langle\e_{(k,j)},\widehat{\mathbf{Y}}\rangle$ and thus $\widehat{\mathbf{X}}= \widehat{\mathbf{Y}}$.

\begin{remark}\phantomsection \label{rem:proof UAT1}
\begin{enumerate}
    \item \label{rem:proof UAT1 i} We presented the proof only for $p\in [2,3)$. However, an analogous and simpler result can be given for continuous functionals of w.g.~\cadlag $p$-rough path, with $p\in[1,2)$. Indeed, the only difference in this setting concerns the interpretation of the Marcus RDE \eqref{Rdemarcus} as the differential equation \eqref{eq:MarcusYoung}.
    \item On the subset of continuous paths the continuity of the solution map \eqref{solutionmap} can be deduced by standard theorems in rough paths theory. We refer to  e.g.~Theorem~10.26 and Corollary 10.28 in \cite{FV:10}, where continuity with respect to the variation metric $d_p$ \eqref{eq:hommetric} is also provided.
\end{enumerate}
   \end{remark}

\subsection{Proof of Theorem \ref{th: UAT non antici}}\label{appendixA:UAT non antici}
 We start by summarizing the reasoning used in the proof of Theorem \ref{th: UAT non antici}. Recall that we analyze here only the case of $p \in [2,3)$ but as explained in Remark \ref{rem:proof UAT1} \ref{rem:proof UAT1 i}, the proof can be easily adapted to $p\in [1,2)$. 
 
Consider the 
set of stopped w.g. rough paths given by 
\begin{align*}
			H:= \left\{\widehat{\mathbf{X}}^t\colon t\in [0,1],\widehat{\mathbf{X}}\in K\right\}.
\end{align*}
We need to show that continuous path functionals can be uniformly approximated on $H$ by linear functionals of the signature evaluated at the deterministic time of stopping $t\in[0,1]$. To this end, we rely on Proposition~\ref{UAT1}. 
In particular, in order to recover all the conditions needed to apply Proposition \ref{UAT1}, we consider the auxiliary  subset of $\widehat{D}^p([0,1],G^{2}(\Rset^{d+2}))$ given by
\begin{align*}
			 \widehat  H:= \left\{{\reallywidehat{\widehat{\mathbf{X}}^t}}\colon t\in [0,1],{\widehat{\mathbf{X}}}\in K\right\},
		\end{align*}
  consisting of all 
   w.g.~\cadlag rough paths with values in $G^{2}(\Rset^{d+2})$ 
  obtained as
  time extension of paths in $ H $. 
  Observe that by Remark \ref{rem: time-extended paths} \ref{rem: time-extended paths ii} $\widehat H$ is a well-defined subset of $\widehat{D}^p([0,1], G^2(\Rset^{d+2}))$.

Next, we first prove that the set $\widehat H$
is relatively compact in $D([0,1], G^2(\Rset^{d+2}))$ (Step 1). Then we show that its closure is bounded with respect to the $p$-variation norm and that it is a subset of the time-extended w.g. \cadlag rough paths
  $\widehat{D}^p([0,1],G^{2}(\Rset^{d+2}))$ (Step 2).  After introducing an auxiliary path functional and an application of Proposition~\ref{UAT1} (Step 3), we then just need to conclude that the obtained linear functional is of the desired form (Step 4).
\begin{remark}\label{rem: G2metric}
	Recall that $ \widehat{D}^p([0,1],G^{2}(\Rset^{d+2}))\subset D([0,1], G^2(\Rset^{d+2}))$ and that the metric space $(G^2(\Rset^{d+2}), d_{CC})$ is a Polish space where closed bounded sets are compact (see Corollary 7.50 in \cite{FV:10}). Therefore, a set $B\subset G^2(\Rset^{d+2}) $ such that $\sup_{{\bf x}\in B}\|{\bf x}\|_{CC}<\infty $ has compact closure.
 
 Furthermore, let $\rho$ denote the distance on $G^2(\Rset^{d+2})$ induced by the norm \eqref{eq:normTN}. By Proposition 7.49 in \cite{FV:10}, the map
	\begin{align*}
		id:(G^2(\Rset^{d+2}), d_{CC})\rightarrow (G^2(\Rset^{d+2}), \rho)
	\end{align*} is Lipschitz on bounded sets.
\end{remark}
\begin{notation}\label{notation proof UAT non anticipativ}
    Throughout the proof, for all $t\in[0,1]$, $\widehat{\mathbf{X}}\in K$, we set $\widecheck{\mathbf{X}}^t:=\reallywidehat{\widehat{\mathbf{X}}^t}$ and denote by $\widecheck{\mathbb{X}}^t$ its signature.
\end{notation}

\paragraph{Step 1:} We show that the set 
		$\widehat H$
		is relatively compact in $D([0,1], G^2(\Rset^{d+2}))$. By Theorem~12.3 in \cite{billing}, a set $A\subset D([0,1], G^2(\Rset^{d+2}))$ is relatively compact if and only if 
		\begin{align}\label{eq:setYrelcomp}
			\{\mathbf{Z}_s \ : \ s\in [0,1],\ \mathbf{Z}\in A\}
		\end{align}
		is relatively compact in $(G^2(\Rset^{d+2}), d_{CC})$  and
		\begin{align}\label{eq:2compact}
			\lim_{\delta\rightarrow0}\sup_{\mathbf{Z} \in A}\Big(\inf_{\{t_i\}_\delta}\max_{1\leq i\leq \nu}\sup_{s,u\in [t_{i-1},t_{i})}d_{CC}(\mathbf{Z}_s,\mathbf{Z}_u)\Big)=0
		\end{align}
		where $\{t_i\}_\delta$ denotes a partition  $0=t_0<t_1<\dots <t_\nu=1$ of $[0,1]$, satisfying $\min_{1\leq i\leq \nu}(t_i-t_{i-1})>\delta$. 
  First of all, notice that since $K$ is compact in $\widehat{D}^p([0,1],G^{2}(\Rset^{d+1}))$ and each element of $K$ is already time-extended, the set\footnote{By Notation \ref{notation proof UAT non anticipativ} $  {\widehat H}^1:=\left\{\reallywidehat{\widehat{\mathbf{X}}}\colon {\widehat{\mathbf{X}}}\in K\right\}.$}
  \begin{align*}
    {\widehat H}^1:=\left\{\widecheck{\mathbf{X}}^1\colon {\widehat{\mathbf{X}}}\in K\right\}\subset \widehat H,
  \end{align*}
   is compact in ${D}([0,1],G^{2}(\Rset^{d+2}))$.
		By Remark \ref{rem: G2metric} it thus holds that
  $$\sup_{{\widehat{\mathbf{X}}}\in  K}\sup_{s\in [0,1]}\|\widecheck{\mathbf{X}}^1_s\|_{CC}<\infty.$$
		By definition of $\widehat H$ and since $\|\widecheck{\mathbf{X}}^t_{t,s}\|_{CC}=(s-t)$ for each $s\in [t,1]$,  we can then compute  
		\begin{align*}
			\sup_{\widehat{\mathbf{Y}}\in \widehat H}\sup_{s\in [0,1]}\|\widehat{\mathbf{Y}}_s\|_{CC}
   &= \sup_{t\in [0,1],{\widehat{\mathbf{X}}}\in  K}\sup_{s\in [0,1]}\|\widecheck{\mathbf{X}}^t_s\|_{CC}\\
   &\leq \sup_{t\in [0,1],{\widehat{\mathbf{X}}}\in  K}\max\Big(
   \sup_{s\in [0,t]}\|\widecheck{\mathbf{X}}^1_s\|_{CC},
   \sup_{s\in [t,1]}(\|\widecheck{\mathbf{X}}^1_t\|_{CC}+\|\widecheck{{\mathbf{X}}}^t_{t,s}\|_{CC})\Big)<\infty.
		\end{align*}
  This shows that \eqref{eq:setYrelcomp} is satisfied by $\widehat H$.
		Similarly, observe that
$d_{CC}(\widecheck{{\mathbf{X}}}^t_s,\widecheck{{\mathbf{X}}}^t_u)=(u-s)$ for each $t\leq s\leq u$ and 
$$d_{CC}(\widecheck{{\mathbf{X}}}^t_s,\widecheck{{\mathbf{X}}}^t_u)
\leq 
d_{CC}(\widecheck{\mathbf{X}}^1_s,\widecheck{\mathbf{X}}^1_t)
+(u-t)
$$
for each $s\leq t\leq u$. We can then compute
		\begin{align*}
			&\lim_{\delta\rightarrow0}\sup_{\widehat{\mathbf{Y}}\in \widehat  H}\Big(\inf_{\{t_i\}_\delta}\max_{1\leq i\leq \nu}\sup_{s,u\in [t_{i-1},t_{i})}d_{CC}(\widehat{\mathbf{Y}}_s,\widehat{\mathbf{Y}}_u)\Big)\\
			&\qquad =\lim_{\delta\rightarrow0}\sup_{t\in [0,1],{\widehat{\mathbf{X}}}\in  K}\Big(\inf_{\{t_i\}_\delta}\max_{1\leq i\leq \nu}\sup_{s,u\in [t_{i-1},t_{i})}d_{CC}(\widecheck{{\mathbf{X}}}^t_s,\widecheck{{\mathbf{X}}}^t_u)\Big)=0,
		\end{align*}
  proving that $\widehat H$ satisfies \eqref{eq:2compact} too and is thus relatively compact.
		
		\paragraph{Step 2:} Let  ${cl(\widehat H)}$ be the closure of $\widehat H$ in $D([0,1], G^2(\Rset^{d+2}))$. We now show that $cl(\widehat H)$ is bounded with respect to the $p$-variation norm and that it is a subset of the time-extended w.g. rough paths over $\Rset^{d+2}$, i.e.
  $cl(\widehat H)\subset \widehat{D}^p([0,1],G^{2}(\Rset^{d+2})).$

 Let $(\widehat{\mathbf{Y}}^n)_{n}\subseteq \widehat H$ be a sequence and ${\mathbf{Y}}\in D([0,1], G^2(\Rset^{d+2}))$ be such that $\widehat{\mathbf{Y}}^n\rightarrow \mathbf{Y}$ in the $J_1$-topology.  By definition, there exists a sequence of time-reparametrization $(\lambda_n)_n$ such that 
 \begin{equation}\label{eq:lambdan}
				\lim_{n\rightarrow\infty}\sup_{s\in [0,1]}|\lambda_n(s)-s|=0,\quad\text{and}\quad
				\lim_{n\rightarrow\infty}\sup_{s\in [0,1]}d_{CC}(\widehat{\mathbf{Y}}^n_{\lambda_n(s)},\mathbf{Y}_s)=0. 
   \end{equation}
	 Let then $\mathcal D$ be  a fixed partition  of $[0,1]$. 
   Observe that since by assumption $K$ is bounded with respect to $\|\,\cdot\,\|_{p-var}$, $\widehat{H}^1$ and thus $\widehat H$ are bounded as well.
   We can thus compute
			\begin{align*}
				\sum_{t_i\in \mathcal D}d_{CC}(\mathbf{Y}_{t_i},\mathbf{Y}_{t_{i+1}})^p
    =&\lim_{n\rightarrow\infty } \sum_{t_i\in \mathcal D} d_{CC}\left (\widehat{\mathbf{Y}}^n_{\lambda_n(t_i)},\widehat{\mathbf{Y}}^n_{\lambda_n(t_{i+1})}\right)^p
				\leq \limsup_{n\rightarrow\infty }\|\widehat{\mathbf{Y}}^n_{\lambda_n(\cdot)}\|^p_{p-var}\\
				=&\limsup_{n\rightarrow\infty }\|\widehat{\mathbf{Y}}^n\|^p_{p-var}
    \leq \sup_{n\in \Nset} \|\widehat{\mathbf{Y}}^n\|^p_{p-var},
			\end{align*}
			where the second equality follows by the invariance of the $p$-variation norm with respect to time-reparametrization. Taking the sup over all the partition we get that $cl(\widehat{H})$ is bounded with respect to the $p$-variation norm.

			Next, note that $\lambda_n (0)=0$ for all $n\in \Nset$. Therefore, it holds that $\widehat{\mathbf{Y}}^n_{\lambda_n(0)}=\widehat{\mathbf{Y}}^n_0\rightarrow \mathbf{Y}_0={\bf 1}\in G^2(\Rset^{d+2})$. Moreover, by
   $$\lim_{n\to\infty}\sup_{s\in[0,1]}(\big|\|\widehat{\mathbf{Y}}^n_{\lambda_n(s)}\|_{CC}-\|{\mathbf{Y}}_s\|_{CC}\Big|)
   \leq \lim_{n\rightarrow\infty}\sup_{s\in [0,1]}d_{CC}(\widehat{\mathbf{Y}}^n_{\lambda_n(s)},\mathbf{Y}_s)=0,$$
   and the fact that $\bf Y$ is \cadlag we deduce that the set 
   $\{{\bf Y}_s\colon s\in [0,1], n\in \mathbb N\}\cup\{\widehat{\bf Y}^n_{\lambda_n(s)}\colon s\in [0,1], n\in \mathbb N\}$
is bounded.
By Remark \ref{rem: G2metric} we thus get 
			\begin{align*}
   \lim_{n\rightarrow\infty}\sup_{s\in [0,1]}\rho(\widehat{\mathbf{Y}}^n_{\lambda_n(s)},\mathbf{Y}_s)=0,
			\end{align*}
			implying in particular that for all $|I|\leq 2$ 
			\begin{align*}
				\lim_{n\rightarrow\infty}\sup_{s\in [0,1]} |\langle \e_I,\widehat{\mathbf{Y}}^n_{\lambda_n(s)}-\mathbf{Y}_s\rangle |=0.
			\end{align*}
			Since for all $n\in\Nset$, $s\in [0,1]$, $\langle \e_{-1}, \widehat{\mathbf{Y}}^n_{\lambda_n(s)}\rangle=\lambda_n(s)$, by condition \eqref{eq:lambdan}, $\langle \e_{-1},\mathbf{Y}_s\rangle =s$ for all $s\in [0,1]$.
			Furthermore, since continuity points of ${\widehat{\bf Y}}^n$ have Lebesgue measure 1, for all $i\in \{-1,1,\dots,d+1\}$ we have 
			\begin{align*}
				\langle \e_i\otimes \e_{-1},\widehat{\mathbf{Y}}^n_{\lambda_n(s)}\rangle= &\int_0^{\lambda_n(s)}  \langle \e_i,\widehat{\mathbf{Y}}^n_{u^-} \rangle du=\int_0^{\lambda_n(s)} \langle \e_i,\widehat{\mathbf{Y}}^n_{u} \rangle du
				=\int_0^1 \langle \e_i,\widehat{\mathbf{Y}}^n_{u} \rangle 1_{\{u\leq \lambda_n(s)\}}du.
			\end{align*}
   Finally, observe that for fixed $s\in [0,1]$ and for all continuity points $u\in [0,1]\setminus\{s\}$ of $\mathbf{Y}$,
			\begin{align*}
				\lim_{n\rightarrow\infty} \langle \e_i,\widehat{\mathbf{Y}}^n_{u} \rangle 1_{\{u\leq \lambda_n(s)\}}= \langle \e_i,\mathbf{Y}_{u} \rangle 1_{\{u\leq s\}}.
			\end{align*}
			Thus, the dominated convergence theorem yields
			\begin{align*}
				\lim_{n\rightarrow\infty} \langle \e_i\otimes \e_{-1},\widehat{\mathbf{Y}}^n_{\lambda_n(s)}\rangle = \langle \e_i\otimes \e_{-1},\mathbf{Y}_s\rangle =\int_0^s\langle \e_i,\mathbf{Y}_u\rangle du,
			\end{align*}
	and the claim follows.

       \paragraph{Step 3:} Set $\widehat F:\widehat{D}^p([0,1],G^{2}(\Rset^{d+2}))\rightarrow \Rset$ as $\widehat F( \widehat {\bf Z}):=F( {{\bf Z}})$ where ${{\bf Z}} \in {D}^p([0,1],G^{2}(\Rset^{d+1}))$ denotes the w.g.~rough path obtained by removing the time-extension of $\widehat {\bf Z}$ (see Remark \ref{rem: time-extended paths} \ref{rem: time-extended paths ii}). Observe that by Notation \ref{notation proof UAT non anticipativ}, if in particular $\widehat{\bf Z}\in \widehat{H}$, then $\widehat{\bf Z}=\widecheck{{\mathbf{X}}}^t$ and ${\bf Z}=\widehat {\bf X}^t$ for some $t\in[0,1]$ and $\widehat{\mathbf{X}}\in  K$.
       Since  $\widehat{F}$ is a
continuous path functional,
 an application of Proposition~\ref{UAT1} to the compact set $cl(\widehat{H})$ yields the existence of a linear functional $\tilde l\in\mathscr L_{\Rset^{d+2}}$ such that
	\[
 \sup_{t\in[0,1]}\sup_{{\widehat{\mathbf{X}}}\in  K}| \widehat F({\widecheck{{\mathbf{X}}}^t})-\tilde l({\widecheck{\mathbb{X}}}^t_1)|
 \leq
	\sup_{\widehat{\mathbf{Y}}\in cl(\widehat{H} )}|\widehat F(\widehat{\mathbf{Y}})-\tilde l(\mathbb{\widehat{Y}}_1)|\leq\varepsilon.
	\]

By definition of the path functional $\widehat{F}$ and Notation \ref{notation proof UAT non anticipativ}, we then conclude that
\begin{align*}
    \sup_{t\in[0,1]}\sup_{{\widehat{\mathbf{X}}}\in  K} | F( \widehat{\bf X}^t)-\tilde l({\widecheck{\mathbb{X}}}^t_1)|\leq\varepsilon.
\end{align*}
		\paragraph{Step 4:} We show that for every linear functional $\tilde l\in \mathscr{L}_{\Rset^{d+2}}$, there exists another linear functional ${l}\in \mathscr{L}_{\Rset^{d+1}}$ such that for all $t\in [0,1]$ and $\widehat{\mathbf{X}}\in K$  $$\tilde l({\widecheck{\mathbb{X}}}^t_1)={l}(\widehat{\mathbb{X}}_t),$$
  where $\widehat
  {\mathbb{X}}_t$ denotes the signature at time $t$ of the path $\widehat{\bf X}$.
Observe first that since each component of $\widecheck {\mathbb X}_{t,1}^t$ is a polynomial in $t$,
for each multi-index $I$ there is a map $l_I\in \mathscr L_{\Rset^{d+2}}$ such that 
$$\langle \e_I,\widecheck{\mathbb X}_{t,1}^t\rangle=l_I(\widecheck {\mathbb X}^1_t),$$ where $\widecheck {\mathbb X}^1_t$ denotes the signature at time $t$ of the path $\widecheck{\mathbf{X}}^1\in H^1$.
This in particular implies that for each multi-index $I$ 
$$\langle \e_I,\widecheck {\mathbb X}^t_1\rangle
=\langle\e_I,\widecheck {\mathbb X}^1_t\otimes \widecheck {\mathbb X}_{t,1}^t\rangle
=\sum_{\e_{I_1}\otimes \e_{I_2}=\e_I}
\langle\e_{I_1},\widecheck {\mathbb X}^1_t\rangle
\langle\e_{I_2},\widecheck {\mathbb X}_{t,1}^t\rangle
=\sum_{\e_{I_1}\otimes \e_{I_2}=\e_I}
\langle\e_{I_1},\widecheck {\mathbb X}_t^1\rangle
l_{I_2}(\widecheck {\mathbb X}^1_t)=:h_I(\widecheck {\mathbb X}^1_t)$$
for some $h_I\in \mathscr L_{\Rset^{d+2}}$. Since each $\widehat{\bf X}$ in $K$ is already time-extended,  the claim follows. \qed

\begin{remark}
    Observe that Step 4 would not hold true if in the statement of Theorem \ref{th: UAT non antici} we considered $K$ to be a compact subset of the set of w.g.~rough paths which are not necessarily time-extended, i.e.~$K\subset {D}^p([0,1],G^{2}(\Rset^{d+1}))$.
\end{remark}

\section{Proofs of Section \ref{sec:Levy}} 
\subsection{Proof of Proposition~\ref{expsiglevy} }\label{appendixB:expsiglevy}

Fix $d\in \Nset$ and let $L$ be a $\Rset^d$-valued Lévy process with triplet $(b^L,C^L,F^L)$ with respect to the ``truncation function'' $\chi(x)=x$. For notational convenience we write $(b,C,F)=(b^L,C^L,F^L)$. We prove that \ref{it2i} the truncated signature $\mathbb L^N$ of $L$ is a polynomial process in the truncated tensor algebra $T^N(\R^d)$ and \ref{it2ii} the expected signature satisfies \eqref{eqn4}.
\begin{enumerate}
\item\label{it2i} First of all, recall that $\mathbb{L}$ is the solution of a  Marcus SDE \eqref{MarcusSDE1}. Therefore, for each multi-index $I$ with entries in $\{1,\dots,d\}$ we have
\begin{align*}
		\langle\epsilon_{I},\mathbb{L}_t\rangle =&\int_{0}^t \langle\epsilon_{I'},\mathbb{L}_{s^-}\rangle dL_s^{i_{|I|}}
		+C_{i_{{|I|}-1},i_{|I|}}\frac{1}{2}\int_{0}^t\langle\epsilon_{I''},\mathbb{L}_{s^-}\rangle ds \\
		&\qquad+\sum_{0<s\leq t}\sum_{\e_{I_1}\otimes \e_{I_2}=\e_I}\frac{1}{|I_2|!}1_{\{|I_2|> 1\}}\langle\epsilon_{I_1},\mathbb{L}_{s^-}\rangle\langle \e_{I_2}, (\Delta L_s)^{\otimes |I_2|}\rangle.
	\end{align*}
Thus, $\mathbb{L}^N$ is a $T^N(\Rset^d)$-valued semimartingale. Moreover, by condition \eqref{eqn1}, $a.s.$
\[
\sum_{\e_{I_1}\otimes \e_{I_2}=\e_I}\frac{1}{|I_2|!}1_{\{|I_2|> 1\}}|\langle\epsilon_{I_1},\mathbb{L}_{s^-}\rangle|\int_{\Rset^d}|\langle \e_{I_2}, x^{\otimes |I_2|}\rangle |F(dx)<\infty.
\]
Therefore, the process $\int_0^\cdot\int_{\Rset^d}\sum_{\e_{I_1}\otimes \e_{I_2}=\e_I}\frac{1}{|I_2|!}1_{\{|I_2|> 1\}}\langle\epsilon_{I_1},\mathbb{L}_{s^-}\rangle\langle \e_{I_2}, x^{\otimes |I_2|}\rangle F(dx)ds$ is of locally integrable variation and by Theorem~II.1.8 in \cite{JS:87}, 
\[
\int_0^\cdot\int_{\Rset^d}\sum_{\e_{I_1}\otimes \e_{I_2}=\e_I}\frac{1}{|I_2|!}1_{\{|I_2|> 1\}}\langle\epsilon_{I_1},\mathbb{L}_{s^-}\rangle\langle \e_{I_2}, x^{\otimes |I_2|}\rangle (\mu(dx,ds)-F(dx)ds)
\]
is a local martingale.
We now prove that $\mathbb{L}^N$ is a jump-diffusion in the sense of Definition~III.2.18 in \cite{JS:87}. In order to simplify the notation, for each $\mathbf{y}\in T^N(\R^d)$ and each multi-index $I$, we set $\mathbf{y}_I:=\langle \epsilon_I,\mathbf{y}\rangle$ and
$$\delta(\mathbf{y},x)_I:=\sum_{\e_{I_1}\otimes \e_{I_2}=\e_I}\frac {1}{|I_2|!}\mathbf{y}_{I_1}\langle \e_{I_2},x^{\otimes |I_2|}\rangle1_{\{|I_2|>0\}},\quad\quad \gamma(\mathbf{y},x)_I:=\delta(\mathbf{y},x)_I-\mathbf{y}_{I'}x_{i_{|I|}}.$$ An application of It\^o's formula yields that for any bounded $C^2$ function on $T^N(\Rset^d)$
\[
f(\mathbb{L}^N_t)-f(\mathbb{L}^N_0)-\int_0^t\mathcal{A}f(\mathbb{L}^N_s)ds
\]
is a local martingale, where $\mathcal{A}:C^2_b(T^N(\Rset^d))\rightarrow C(T^N(\Rset^d))$ is such that for any $f\in C^2_b(T^N(\Rset^d))$ and $\mathbf{y}\in T^N(\R^d)$
\begin{align}\label{eq:generatorpoly}
    	\mathcal{A}f(\mathbf{y})=&\sum_{1\leq|I|\leq N}\frac{\partial f(\mathbf{y})}{\partial \mathbf{y}_{I}}\mathbf{y}_{I'}b_{i_{|I|}}
	+\sum_{2\leq|I|\leq N}\frac{\partial f(\mathbf{y})}{\partial \mathbf{y}_{I}}\bigg(\frac{1}{2}\mathbf{y}_{I''}C_{i_{|I|-1},i_{|I|}}
	+\int_{\R^d} \gamma(\mathbf{y},x)_IF(dx)\bigg)\nonumber \\
	&+\frac{1}{2} \sum_{1\leq|I|,|J|\leq  N}\frac{\partial^2f(\mathbf{y})}{\partial \mathbf{y}_{I}\partial \mathbf{y}_{J}}\mathbf{y}_{I'}\mathbf{y}_{J'}C_{i_{|I|},j_{|J|}}\\
	&+\int_{\Rset^d} \nonumber  \bigg(f(\mathbf{y}+\delta(\mathbf{y},x))-f(\mathbf{y}) -\sum_{1\leq |I|\leq N}\frac{\partial f(\mathbf{y})}{\partial \mathbf{y}_{I}}\delta(\mathbf{y},x)_I\bigg)F(dx).
\end{align}
The operator $\mathcal{A}$ is a jump-diffusion operator whose coefficients $$(\beta,\zeta,K):=((\beta_I)_{|I|\leq N},(\zeta_{I,J})_{|I|,|J|\leq N}, K)$$ are given by
\begin{align*}
\beta_I(\mathbf{y})&:=\mathbf{y}_{I'}b_{i_{|I|}}1_{|I|\geq1}
	+\bigg(\frac{1}{2}\mathbf{y}_{I''}C_{i_{|I|-1},i_{|I|}}+\int_{\R^d}\gamma(\mathbf{y},x)_I F(dx)\bigg)1_{|I|\geq2},\\
	\zeta_{I,J}(\mathbf{y})&:=\mathbf{y}_{I'}\mathbf{y}_{J'}C_{i_{|I|},j_{|J|}}1_{|I|,|J|\geq1},\\
	K(\mathbf{y},A)&:=\int_{\R^d}1_{A\setminus\{0\}}(\delta(\mathbf{y},x))F(dx).
\end{align*}
Therefore, by Theorem~II.2.42 in \cite{JS:87}, $\mathbb{L}^N$ is a jump-diffusion on $T^N(\Rset^d)$ whose
characteristics, with respect to the ``truncation function" $\chi(\mathbf{y})=\mathbf{y}$, are given by
\[
\int_0^t\beta(\mathbb{L}^N_s)ds, \quad \int_0^t\zeta(\mathbb{L}^N_s)ds \quad \text{and}\quad K(\mathbb{L}^N_{t^-},dx)\times dt.
\]
Notice also that 
$
K(\mathbf{y},\{0\})=F(\{0\})=0$ for all $\mathbf{y}\in T^N(\Rset^d)$,
 and, by condition \eqref{eqn1},
\[
\int_{T^N(\Rset^d)}\|\boldsymbol{\xi}\|\wedge \|\boldsymbol{\xi}\|^2K(\mathbf{y},d\boldsymbol{\xi})=\int_{\Rset^d}\|\delta(\mathbf{y},x)\|\wedge\| \delta(\mathbf{y},x)\|^2F(dx)<\infty.
\]          
Finally, we show that $\mathbb{L}^N$ is a polynomial jump diffusion on $T^N(\Rset^d)$. Once again by condition \eqref{eqn1},
\begin{align*}
	&\int_{T^N(\Rset^d)}\|\boldsymbol{\xi}\|^kK(\mathbf{y},d\boldsymbol{\xi})=\int_{\Rset^d}\|\delta(\mathbf{y},x)\|^kF(dx)<\infty,
\end{align*}
for all $\mathbf{y}\in T^N(\Rset^d)$ and $k\geq 2$.
Moreover, for every set of multi-indices $I_1,\ldots, I_k$ it holds
\begin{align*}
	\int_{T^N(\Rset^d)}\boldsymbol{\xi}_{I_1}\cdots \boldsymbol{\xi}_{I_k}K(\mathbf{y},d\boldsymbol{\xi})=\int_{\Rset^d}\delta(\mathbf{y},x)_{I_1}\cdots \delta(\mathbf{y},x)_{I_k} F(dx)
\end{align*}
is a polynomial of degree at most $k$ on $T^N(\Rset^d)$. Since each component of $\beta$ and $\zeta$ is polynomial on $T^N(\Rset^d)$ respectively of degree $1$ and $2$, by Lemma~2.2 in \cite{FL:20} we can conclude that $\mathbb{L}^N$ is a polynomial jump-diffusion on $T^N(\Rset^d)$.

\item\label{it2ii} We compute now  $\E[\langle \e_I,\L^N_t\rangle]$ for a fixed multi-index $I$ with $|I|\leq N$. To this end we apply the so-called moment formula (see for instance Theorem~2.5 in \cite{FL:20}).
Observe that the projection function 
\begin{align*}
	\langle \epsilon_{I},\cdot\rangle : \ &T^N(\Rset^d)\longrightarrow \Rset \\
	& \mathbf{y}\longmapsto\langle \epsilon_{I},\mathbf{y}\rangle
\end{align*}
is a polynomial of degree 1 on $T^N(\Rset^d)$. Since we already proved in \ref{it2i} that $\L^N$ is a polynomial process, we also know that the corresponding generator $\Acal$ is mapping the space of polynomials of degree at most one on $T^N(\Rset^d)$ to itself. More precisely, it follows from \eqref{eq:generatorpoly} that for each multi-index $J$ we have 
\[
\mathcal{A}\langle\epsilon_{J},\cdot\rangle=\sum_{\e_{J_1}\otimes\e_{J_2}=\e_J}\langle \e_{J_1},Q\rangle
\langle\epsilon_{J_{2}},\cdot\rangle,
\]
where $Q:=\left(0,b,\frac{1}{2} \left(C+\int_{\Rset^d}x^{\otimes 2}F(dx)\right),
	\dots,\frac{1}{k!}\int_{\Rset^d}x^{\otimes k}F(dx),\dots \right)$.
This can equivalently be written as
\[
\mathcal{A}\langle\epsilon_{J},\cdot\rangle(\mathbf{y})=\langle \e_{J},Q\otimes \mathbf{y}\rangle
\]

By Theorem~2.5 in \cite{FL:20} we can then conclude that
\begin{align*}\Ex[\langle \epsilon_{I},\mathbb{L}_t^N\rangle ]
	= \sum_{k=0}^{\infty}\frac{t^k}{k!}\mathcal A^k\langle\e_I,\cdot\rangle(\e_\emptyset)
	=\sum_{k=0}^{\infty}\frac{t^k}{k!}\langle\e_I,Q^{\otimes k}\otimes\e_\emptyset\rangle
	=\langle\e_I,\exp(tQ)\otimes\e_\emptyset\rangle.
\end{align*}

\end{enumerate}

\subsection{Proof of Proposition~\ref{index-transf}}\label{appendixB:index-transf}

The claim is clear for $I=\emptyset$. If $I\neq \emptyset$, recall from equation \eqref{component-wise} that
	\begin{align*}
		\langle\epsilon_{I}\otimes \e_j,\mathbb{X}_t\rangle =&\int_{0}^t \langle\epsilon_{I},\mathbb{X}_{s^-}\rangle dX_s^{j}
		+1_{\{i_{|I|}=j=0\}}\frac{1}{2}\int_{0}^t\langle\epsilon_{I'},\mathbb{X}_{s^-}\rangle ds \\
		&\quad+\sum_{0<s\leq t}\sum_{\e_{I_1}\otimes \e_{I_2}=\e_I}\frac{1}{(|I_2|+1)!}1_{\{|I_2|>0\}}\langle\epsilon_{I_1},\mathbb{X}_{s^-}\rangle \Delta X_s^{S(I_2)+j}1_{\{j>0\}}1_{\{I_2\in\{1,\ldots,d\}^{|I_2|}\}}.
	\end{align*}
Let $\mathcal I(r,k)$ be the set of all multi-indices $J\in \Nset^k$ such that $S(J)=r$. Set 
$$\alpha(r,k):=\sum_{J\in \mathcal I(r,k)}\prod_{i=1}^k\frac 1 {(j_i+1)!}$$ 
and note that
\begin{align*}
	\int_{0}^t \langle\epsilon_{I},\mathbb{X}_{s^-}\rangle dX_s^{j}
	 =&\langle\epsilon_{I}\otimes \epsilon_{j},\mathbb{X}_t\rangle-1_{\{i_{|I|}=j=0\}}\frac{1}{2}\int_{0}^t\langle\epsilon_{I'},\mathbb{X}_{s^-}\rangle ds \\
	&\qquad-\sum_{0<s\leq t}\sum_{\e_{I_1}\otimes \e_{I_2}=\e_I}\alpha(|I_2|,1)\langle\epsilon_{I_1},\mathbb{X}_{s^-}\rangle \Delta X_s^{S(I_2)+j}1_{\{j>0\}}1_{\{I_2\in\{1,\ldots,d\}^{|I_2|}\}}\\
	=&\langle\epsilon_{I}\otimes \epsilon_{j},\mathbb{X}_t\rangle-1_{\{i_{|I|}=j=0\}}\frac{1}{2}\int_{0}^t\langle\epsilon_{I'},\mathbb{X}_{s^-}\rangle ds \\
	&\qquad-\sum_{\e_{I_1}\otimes \e_{I_2}=\e_I}\alpha(|I_2|,1)\int_0^t\langle\epsilon_{I_1},\mathbb{X}_{s^-}\rangle dX_s^{S(I_2)+j}1_{\{j>0\}}1_{\{I_2\in\{1,\ldots,d\}^{|I_2|}\}}.
\end{align*}

In particular, observe that for each $j>0$ it holds
\begin{align*}
	\int_0^t\langle\epsilon_{I},\mathbb{X}_{s^-}\rangle \ dX_s^{j}
	=&\langle\epsilon_{I}\otimes \epsilon_{j},\mathbb{X}_t\rangle 
	-\sum_{\e_{I_1}\otimes \e_{I_2}=\e_I}\alpha(|I_2|,1)\int_0^t\langle\epsilon_{I_1},\mathbb{X}_{s^-}\rangle dX_s^{S(I_2)+j}1_{\{j>0\}}1_{\{I_2\in\{1,\ldots,d\}^{|I_2|}\}}.
	\end{align*}
	We now prove by induction that for each $r\in \mathbb N_0$ it holds
		\begin{align*}
	\int_{0}^t \langle\epsilon_{I},\mathbb{X}_{s^-}\rangle dX_s^{j}
	-&\langle\epsilon_{I}\otimes \epsilon_{j},\mathbb{X}_t\rangle
	+1_{\{i_{|I|}=j=0\}}\frac{1}{2}\int_{0}^t\langle\epsilon_{I'},\mathbb{X}_{s^-}\rangle ds \\
	&=\sum_{k=1}^r(-1)^k\sum_{\e_{I_1}\otimes \e_{I_2}=\e_I}\alpha(|I_2|,k)\langle\epsilon_{I_1}\otimes \e_{S(I_2)+j},\mathbb{X}_{t}\rangle1_{\{j>0\}}1_{\{I_2\in\{1,\ldots,d\}^{|I_2|}\}}\\
	&\qquad (-1)^{r+1}
	\sum_{\e_{I_1}\otimes \e_{I_2}=\e_I}\alpha(|I_2|,r+1)\int_0^t\langle\epsilon_{I_1},\mathbb{X}_{s^-}\rangle dX_s^{S(I_{2})+j}1_{\{j>0\}}1_{\{I_2\in\{1,\ldots,d\}^{|I_2|}\}}.
\end{align*}
We already showed that the claim holds for $r=0$. Supposing that it is true for $r-1$, we get
	\begin{align*}
	    \sum_{\e_{I_1}\otimes \e_{I_2}=\e_I}&\alpha(|I_2|,r)\int_0^t\langle\epsilon_{I_1},\mathbb{X}_{s^-}\rangle dX_s^{S(I_{2})+j}1_{\{j>0\}}1_{\{I_2\in\{1,\ldots,d\}^{|I_2|}\}}\\
	    	&= \sum_{\e_{I_1}\otimes \e_{I_2}=\e_I}\alpha(|I_2|,r)\langle\epsilon_{I_1}\otimes \e_{S(I_2)+j},\mathbb{X}_{t}\rangle1_{\{j>0\}}1_{\{I_2\in\{1,\ldots,d\}^{|I_2|}\}}\\
	    	&\qquad-\sum_{\e_{I_1}\otimes \e_{I_2}=\e_I}\alpha(|I_2|,r)\sum_{\e_{I_{11}}\otimes \e_{I_{12}}=\e_{I_1}}\alpha(|I_{12}|,1)\\
	    	&\qquad\times\int_0^t\langle\epsilon_{I_{11}},\mathbb{X}_{s^-}\rangle dX_s^{S(I_{12})+S(I_2)+j}1_{\{I_{12}\in\{1,\ldots,d\}^{|I_{12}|}\}}
	1_{\{j>0\}}
	1_{\{I_2\in\{1,\ldots,d\}^{|I_2|}\}}\\
		&\stackrel{(\star)}=
		 \sum_{\e_{I_1}\otimes \e_{I_2}=\e_I}\alpha(|I_2|,r)\langle\epsilon_{I_1}\otimes \e_{S(I_2)+j},\mathbb{X}_{t}\rangle1_{\{j>0\}}1_{\{I_2\in\{1,\ldots,d\}^{|I_2|}\}}\\
	    	&\qquad-\sum_{\e_{I_{11}}\otimes \e_{I_{122}}=\e_I}\alpha(|I_{122}|,r+1)
		\int_0^t\langle\epsilon_{I_{11}},\mathbb{X}_{s^-}\rangle dX_s^{S(I_{122})+j}
	1_{\{j>0\}}1_{\{I_{122}\in\{1,\ldots,d\}^{|I_{122}|}\}},
	\end{align*}
	where in $(\star)$ we set $\e_{I_{122}}=\e_{I_{12}}\otimes \e_{I_2}$.
	Since $\alpha(r,k)=0$ for $k>r$, we can conclude that
		\begin{align*}
	\int_{0}^t \langle\epsilon_{I},\mathbb{X}_{s^-}\rangle dX_s^{j}
	-&\langle\epsilon_{I}\otimes \epsilon_{j},\mathbb{X}_t\rangle
	+1_{\{i_{|I|}=j=0\}}\frac{1}{2}\int_{0}^t\langle\epsilon_{I'},\mathbb{X}_{s^-}\rangle ds \\
	&=\sum_{\e_{I_1}\otimes \e_{I_2}=\e_I}\sum_{k=1}^{|I_2|}(-1)^k\alpha(|I_2|,k)\langle\epsilon_{I_1}\otimes \e_{S(I_2)+j},\mathbb{X}_{t}\rangle1_{\{j>0\}}1_{\{I_2\in\{1,\ldots,d\}^{|I_2|}\}}.
\end{align*}

\subsection{Proof of Theorem~\ref{fromStoX}}\label{appendixB:fromStoX}

Let $k\in \Nset$ and observe that 
$
\big(\sum_{h=1}^M\ell_h\langle \epsilon_{J^h},\mathbb{X}_{s^-}\rangle \big)^k
$
is a polynomial of degree $k$ in $(\ell_h)_{1\leq h\leq M}$. Precisely by the multinomial theorem,
for each
$
\alpha=(\alpha^1,\alpha^2,\dots,\alpha^M)\in \Nset^M
$ it holds 
\begin{equation}\label{notationpolynom}
	\Big(\sum_{h=1}^M\ell_h\langle \epsilon_{J^h},\mathbb{X}_{s^-}\rangle \Big)^k=k!\sum_{\alpha\in \Nset^M, \ S(\alpha)=k} (\ell_{1,\dots,M})^{\alpha}\langle \epsilon^{\shuffle \alpha},\mathbb{X}_{s^-}\rangle.
\end{equation}
Setting $U_I(\boldsymbol{\ell})$ as in the statement of the proposition, we prove 
	$	\langle\e_I,\widehat{\mathbb{S}}(\boldsymbol{\ell})_t\rangle=\langle \epsilon_{U_I(\boldsymbol{\ell})},\mathbb{X}_t\rangle$
	proceeding by induction on $m=|I|$. Since $	\langle\e_\emptyset,\widehat{\mathbb{S}}(\boldsymbol{\ell})_t\rangle=1=\langle \e_\emptyset,\mathbb{X}_t\rangle$ the claim follow for $m=0$.
	Next, suppose that the claim holds for $|I|=m-1$ and fix $I$ with $I=(i_1,\ldots,i_m)$. By equation \eqref{sigofS}, $\langle \epsilon_{I}, \widehat{\mathbb{S}}(\boldsymbol{\ell})_t\rangle$ can be written as
	\begin{align*}
	    &1_{\{i_{|I|}=-1\}}\int_0^t\langle \epsilon_{I'}, \widehat{\mathbb{S}}(\boldsymbol{\ell})_{s^-}\rangle \ ds +1_{\{i_{|I|}=1\}}\int_0^t\langle \epsilon_{I'}, \widehat{\mathbb{S}}(\boldsymbol{\ell})_{s^-}\rangle \left(\sum_{h=1}^M\ell^{W}_h\langle \epsilon_{J^h},\mathbb{X}_{s^-}\rangle \right)dW_s\\
		&\quad+ 1_{\{i_{|I|}=1\}}\int_0^t\int_\Rset\langle \epsilon_{I'}, \widehat{\mathbb{S}}(\boldsymbol{\ell})_{s^-}\rangle \left(\sum_{h=1}^M\ell^{\mu}_h\langle \epsilon_{J^h},\mathbb{X}_{s^-}\rangle \right)x\  ({\mu}-\nu)(ds,dx)\\		&\quad+1_{\{(i_{|I|-1},i_{|I|})=(1,1)\}}\frac{1}{2}\int_0^t\langle \epsilon_{I''}, \widehat{S}(\boldsymbol{\ell})_{s^-}\rangle \left(\sum_{h=1}^M\ell^{W}_h\langle \epsilon_{J^h},\mathbb{X}_{s^-}\rangle \right)^2 ds\\
		&\quad+\sum_{\e_{I_1}\otimes\e_{I_2}=\e_I}\frac 1 {|I_2|!}1_{\{|I_2|>1\}}1_{\{I_2=(1,\ldots,1)\}}
		\int_0^t\int_\Rset \langle\epsilon_{I_1},\widehat{\mathbb{S}}(\boldsymbol{\ell})_{s^-}\rangle \left(\sum_{h=1}^M\ell^{\mu}_h\langle \epsilon_{J^h},\mathbb{X}_{s^-}\rangle \right)^{|I_2|} x^{|I_2|} \mu(ds,dx).
	\end{align*}
 Therefore, as a direct application of Proposition~\ref{index-transf} and the induction hypothesis,  with  the notation  of \eqref{notationpolynom} we can write
\begin{align*}
    \langle \epsilon_{I}, \widehat{\mathbb{S}}(\boldsymbol{\ell})_t\rangle=&
    \langle (\e_{{U_{I'}(\boldsymbol{\ell})}};{\e_{-1}})^{\thicksim},\mathbb{X}_t\rangle1_{\{i_{|I|}=-1\}}
    +\!\!\sum_{\alpha\in \Nset^M , \ S(\alpha)=1}\!(\ell^W_{{1,\dots,M}})^{\alpha} \langle(\e_{U_{I'}(\boldsymbol{\ell})}\shuffle {\epsilon^{\shuffle \alpha}};\e_0)^{\thicksim},\mathbb{X}_t\rangle1_{\{i_{|I|}=1\}}\\ 
		&+\!\!\sum_{\alpha\in \Nset^M, \ S(\alpha)=1}\!(\ell^\nu_{{1,\dots,M}})^{\alpha} \langle(\e_{U_{I'}(\boldsymbol{\ell})}\shuffle {\epsilon^{\shuffle \alpha}};\e_1)^{\thicksim},\mathbb{X}_t\rangle1_{\{i_{|I|}=1\}}\\
		&+\!\!\sum_{\alpha\in \Nset^M, \ S(\alpha)=2}\!(\ell^W_{{1,\dots,M}})^{\alpha} \langle(\e_{U_{I''}(\boldsymbol{\ell})}\shuffle \epsilon^{\shuffle \alpha};\e_{-1})^{\thicksim},\mathbb{X}_t\rangle1_{\{|I|>1\}}1_{\{(i_{|I|-1},i_{|I|})=(1,1)\}} \\ 
		&+\!\!\sum_{\e_{I_1}\otimes \e_{I_2}=\e_I}\! \sum_{\alpha\in \Nset^M, \ S(\alpha)=|I_2|}(\ell^\nu_{{1,\dots,M}})^{\alpha} \langle (\e_{U_{I_1}(\boldsymbol{\ell})}\shuffle \epsilon^{\shuffle \alpha};\e_{|I_2|})^{\thicksim},\mathbb{X}_t\rangle1_{\{|I_2|>1\}}1_{\{I_2=(1,\dots,1)\}}\\
		=&\langle \epsilon_{U_I(\boldsymbol{\ell})},\mathbb{X}_t\rangle.
\end{align*}
Note that $\langle \epsilon_{U_I(\boldsymbol{\ell})},\mathbb{X}_t\rangle$ is a polynomial of degree $|I|$ in $(\ell^W_{{1,\dots,M}},\ell^\nu_{{1,\dots,M}})\in \Rset^{2M}$ by construction.

\begin{remark}
	The necessity of the condition  $N\geq |I|(nd+1)$ follows from the fact that  
		\begin{equation}\label{eqn7}
	\langle\epsilon_{I},\widehat{\mathbb{S}}(\boldsymbol{\ell})_t\rangle=\langle \epsilon_{U_I(\boldsymbol{\ell})},\mathbb{X}_t\rangle
	\end{equation}
	is well defined if and only if the entries of every multi-index appearing in $U_I(\boldsymbol{\ell})$ is bounded by $N$. More precisely, the claim is clear for $|I|=0$. For $|I|=1$ it follows by noticing that by 	Remark~\ref{gooddefinitiontilda}
the expression $\langle( \epsilon_{J};\epsilon_1)^\thicksim,\mathbb X_t\rangle$ is well defined for each $|J|\leq n$ if and only if $nd+1\leq N$. Finally, observe that well definiteness of
	$\langle(\epsilon^{\shuffle \alpha};\e_{|I|})^{\thicksim}, \mathbb X\rangle$
	 for each  $|I|\geq 2$ and $\alpha$ with $S(\alpha)=|I|$ is a necessary condition for well definiteness of \eqref{eqn7}. This includes in particular the case $\epsilon^{\shuffle \alpha}=(\e_d^{\otimes n})^{\shuffle |I|}$ for which it holds $S(\epsilon^{\shuffle \alpha})=nd|I|$. Proceeding as in Remark~\ref{gooddefinitiontilda} we can conclude that this is the case if and only if $nd|I|+|I|\leq N$. 
	
\end{remark}

\subsection{Proof of Theorem~\ref{localriskminimization}}\label{appendixB:localriskminimization}

Since
\begin{equation}\label{eqn8}
    S(\boldsymbol{\ell})_t=\int_0^t\Big(\sum_{|H|\leq n}\ell^H_W\langle \epsilon_{H},\mathbb{X}_{s^-}\rangle \Big)dW_s+\int_0^t\int_\Rset\Big(\sum_{|H|\leq n}\ell^H_\nu\langle \epsilon_{H},\mathbb{X}_{s^-}\rangle x\Big) (\mu-\nu)(ds,dx)
\end{equation}
we already know that
$$d [S(\boldsymbol{\ell})]^{pred}_t
=\int_0^t\Big(\sum_{|H|\leq n}\ell^J_W\langle \epsilon_{H},\mathbb{X}_{s^-}\rangle \Big)^2+\Big(\sum_{|H|\leq n}\ell^H_\nu\langle \epsilon_{H},\mathbb{X}_{s^-}\rangle \Big)^2\int_\Rset x^2F(dx)ds.$$
Moreover using that $X$ is a L\'evy process by Proposition~\ref{Xislevy} and Remark~\ref{rem3} we also get that
\begin{align*}
    \Ex[\langle \e_I,\mathbb{\widehat{S}(\boldsymbol{\ell})}_T\rangle|\mathcal F_t]
    &=\Ex[\langle \e_{U_I(\boldsymbol \ell)},\mathbb{{X}}_T\rangle|\mathcal F_t]
    =\sum_{\e_J\otimes \e_H=\e_{U_I(\boldsymbol \ell)}}
    \langle \e_J,\mathbb X_t\rangle\Ex[\langle \e_H,\mathbb X_{T-t}\rangle].
\end{align*}
Observe that in the first step we used the result of Theorem~\ref{fromStoX}.
Since for each $|H|\geq 1$ by Remark~\ref{rem2} it holds
$$
\Ex[\langle \e_H,\mathbb X_{T-t}\rangle]=\sum_{k=1}^{|H|}\frac{(T-t)^k}{k! }\sum_{\e_{H_1}\otimes\cdots\otimes\e_{H_k}=\e_H}\langle \e_{H_1},Q\rangle\cdots\langle \e_{H_k},Q\rangle,
$$
 which is a continuous process of finite variation, by the product rule we  know that $t\mapsto\langle \e_J,\mathbb X_t\rangle\Ex[\langle \e_H,\mathbb X_{T-t}\rangle]$ can be decomposed as the sum of a continuous process of finite variation and 
$\int_0^\cdot \Ex[\langle \e_H,\mathbb X_{T-s}\rangle]d\langle \e_J,\mathbb X_s\rangle.$
Observe also that for each $J_1$ and $J_2$ such that $S(J_2)\geq2$, it holds
	\begin{align*}
	    \sum_{0<s\leq t}\langle\epsilon_{J_1},\mathbb{X}_{s^-}\rangle \Delta X_s^{S(J_2)}
	=&\int_0^t\int_{\mathbb R}\langle\epsilon_{J_1},\mathbb{X}_{s^-}\rangle  x^{S(J_2)}  (\mu-\nu)(ds,dx)\\
	&+
	\int_0^t\langle\epsilon_{J_1},\mathbb{X}_{s^-}\rangle ds\int_{\mathbb R}x^{S(J_2)}  F(dx).
	\end{align*}
Representation \eqref{component-wise}  thus yields that $\langle \e_J,\mathbb X\rangle$ can be decomposed into a continuous process of finite variation  and a square integrable martingale $M^J$ given by
\begin{align*}
    M^J_t=&
    \int_{0}^t \langle\epsilon_{J'},\mathbb{X}_{s^-}\rangle dW_s1_{\{j_{|J|}= 0\}} 
   + \int_{0}^t\int_{\mathbb R} \langle\epsilon_{J'},\mathbb{X}_{s^-}\rangle x^{j_{|J|}}(\mu-\nu)(ds,dx)1_{\{j_{|J|}> 0\}} 
    \\
		&+\sum_{\e_{J_1}\otimes \e_{J_2}=\e_J}\frac{1}{|J_2|!}1_{\{|J_2|> 1\}}1_{\{J_2\in\{1,\dots,N\}^{|J_2|}\}}\int_0^t\int_{\mathbb R}\langle\epsilon_{J_1},\mathbb{X}_{s^-}\rangle  x^{S(J_2)}  (\mu-\nu)(ds,dx)\\
		=&\sum_{\e_{J_1}\otimes \e_{J_2}=\e_J}
    \gamma_W(J_2)\int_{0}^t \langle\epsilon_{J_1},\mathbb{X}_{s^-}\rangle dW_s +\gamma_\nu(J_2)\int_0^t\int_{\mathbb R}\langle\epsilon_{J_1},\mathbb{X}_{s^-}\rangle  x^{S(J_2)}  (\mu-\nu)(ds,dx),
\end{align*}
for $\gamma_W(J_2)=1_{\{J_2=(0)\}}$ and
$\gamma_\nu(J_2)=\frac{1}{|J_2|!}1_{\{|J_2|\geq 1\}}1_{\{J_2\in\{1,\dots,N\}^{|J_2|}\}}.$
Define $V_t:=\mathbb E\big[\langle \epsilon_{I},\widehat{\mathbb{S}}(\boldsymbol{\ell})_T\rangle\big|\mathcal{F}_t\big]$ for each  $t\in[0, T]$. Writing $J_3$ in place of $H$ and using that $V$ is a true martingale by definition, we can thus conclude that
\begin{align*}
    V_t
    &=V_0
    +
    \sum_{\e_J\otimes \e_{J_3}=\e_{U_I(\boldsymbol \ell)}}
    \int_0^t \Ex[\langle \e_{J_3},\mathbb X_{T-s}\rangle]dM_s^J\\
    &=
    V_0
    +
    \sum_{\e_{J_1}\otimes\e_{J_2} \otimes \e_{J_3}=\e_{U_I(\boldsymbol \ell)}}\\
    &\qquad\qquad \int_0^t \langle\epsilon_{J_1},\mathbb{X}_{s^-}\rangle \Ex[\langle \e_{J_3},\mathbb X_{T-s}\rangle](\gamma_W(J_2)dW_s+\gamma_\nu(J_2)x^{S(J_2)}  (\mu-\nu)(ds,dx)),
\end{align*}
for $V_0=\Ex[\langle \epsilon_{I},\widehat{\mathbb{S}}(\boldsymbol{\ell})_T\rangle]$.
 Using \eqref{eqn8} we can then compute 
\begin{align*}
    d[V,{S(\boldsymbol{\ell})}]^{pred}_t
    &=  
    \sum_{\e_{J_1}\otimes \e_{J_2}\otimes \e_{J_3}=\e_{U_I(\boldsymbol \ell)}}
    \sum_{|H|\leq n}
      \langle\epsilon_{J_1},\mathbb{X}_{t^-}\rangle
      \gamma_2(J_2,H,\boldsymbol{\ell})
      \Ex[\langle \e_{J_3},\mathbb X_{T-t}\rangle]\langle \epsilon_{H},\mathbb{X}_{t^-}\rangle dt,
\end{align*}
where $\gamma_2(J_2,H,\boldsymbol\ell)=\gamma_W(J_2)\ell^H_W
      +
    \gamma_\nu(J_2)\ell^H_\nu\int_{\mathbb R} x^{S(J_2)+1} F(dx)$.
The claim follows.

\subsection{Proof of Proposition \ref{prop1}}\label{proof: prop1}
First of all, observe that by Theorem~9 in \cite{PS:06}, the pair $(f(t),g(x))$ is such that the solution of the SDE \eqref{exponentialmartingale} is a true martingale.  

Since $g$ is deterministic and independent of time, $\mu$ is a homogeneous $\P$-Poisson random measure (see e.g.~Theorem~II.4.8 in \cite{JS:87}). Moreover, $e^{g(x)}F(dx)$ is a Lévy measure such that $\int_{|x|>1}|x|^k e^{g(x)}F(dx)<\infty$, for all $k\in\Nset$. Thus, proceeding as in the proof of Proposition~\ref{Xislevy}, we conclude that the process $Y$ defined in \eqref{eq:Y} is a Lévy process with finite moments of all orders.

By equation \eqref{SunderP}, we know that the $\mathbb{P}$-characteristics of $S(\boldsymbol{\ell})$ with respect to the ``truncation function" $\chi(x)=x$ are given by 
\begin{align*}
		   &B^S_t=\int_0^t\bigg(\sum_{|J|\leq n}\ell^J_W\langle \epsilon_{J},\mathbb{X}_{s^-}\rangle f(s)+\sum_{|J|\leq n}\ell^J_\nu\langle \epsilon_{J},\mathbb{X}_{s^-}\rangle  \int_\Rset x(e^{g(x)}-1)F(dx)\bigg)ds,\\
		   &C^S_t=\int_0^t\bigg(\sum_{|J|\leq n}\ell^J_W\langle \epsilon_{J},\mathbb{X}_{s^-}\rangle \bigg)^2ds, \\
		   &\nu^S(dt\times dx)=K^{\mathbb{P}}_t(dx)\times dt,
		\end{align*}
		where
$K^{\mathbb{P}}_t(A)=\int_\Rset 1_{A\setminus\{0\}}\big(\sum_{|J|\leq n}\ell^J_\nu\langle \epsilon_{J},\mathbb{X}_{t^-}\rangle \ x\big)e^{g(x)}F(dx).
$
Since by condition \eqref{conditionmart} we know that
$
    \int_\Rset y(e^{g(x)}-1)F(dx)<\infty,
$
and 
$$ f(t)=\sum_{|I|\leq p}f^{I}\langle \epsilon_{I},\mathbb{X}_{t^-}\rangle=\sum_{|I|\leq p}f^{I}\langle \epsilon_{I},\mathbb{Y}_{t^-}\rangle$$
it suffices to show that for each $|J|\leq n$ it holds
$
    \langle \epsilon_{J},\mathbb{X}_{t}\rangle=\langle \epsilon_{J}^\P,\mathbb{Y}_{t}\rangle
$
 for some linear combination of multi-indices $\e^\P_J$.  We proceed by induction. For $I=\emptyset$ the result is clear. Fix $k\leq n$, $J$ such that $|J|=k$, and suppose that the claim holds for each $|J|< k$. Then equation \eqref{component-wise} yields
$$
\begin{aligned}
		\langle\epsilon_{J},\mathbb{X}_t\rangle 
		=&\int_{0}^t \langle\epsilon_{J'}^\P,\mathbb{Y}_{s^-}\rangle dX_s^{j_k}
		+1_{\{j_{k-1}=j_k=0\}}\frac{1}{2}\int_{0}^t\langle\epsilon_{J''}^\P,\mathbb{Y}_{s^-}\rangle ds \\
		&\qquad+\sum_{0<s\leq t}\sum_{\e_{J_1}\otimes \e_{J_2}=\e_J}\frac{1}{|J_2|!}1_{\{|J_2|> 1\}}\langle\epsilon_{J_1}^\P,\mathbb{Y}_{s^-}\rangle \Delta X_s^{S(J_2)}1_{\{J_2\in\{1,\dots,N\}^{|J_2|}\}}.
	\end{aligned}
$$
Note that since $S(J_2)>2$ we have that $\Delta X_s^{S(J_2)}=\Delta Y_s^{S(J_2)}$.
Moreover, 
\begin{align*}
	\int_{0}^t \langle\epsilon_{J'}^\P,\mathbb{Y}_{s^-}\rangle dW^\mathbb{Q}_s
	=&\int_{0}^t \langle\epsilon_{J'}^\P,\mathbb{Y}_{s^-}\rangle dW^\mathbb{P}_s
	+ f(s) \langle\epsilon_{J'}^\P,\mathbb{Y}_{s^-}\rangle\ ds\\
	=&\int_{0}^t \langle\epsilon_{J'}^\P,\mathbb{Y}_{s^-}\rangle dW^\mathbb{P}_s
	+\sum_{|I|\leq p}f^{I}\langle \epsilon_{J'}^\P\shuffle \epsilon_{I},\mathbb{Y}_{s^-}\rangle ds
\end{align*}
and $\int_{0}^t\int_\Rset \langle\epsilon_{J'}^\P,\mathbb{Y}_{s^-}\rangle  x   (\mu-\nu^\mathbb{Q})(ds,dx)$ can be written as
\begin{align*}
\int_0^t \int_\Rset\langle\epsilon_{J'}^\P,\mathbb{Y}_{s^-}\rangle x   (\mu-\nu^\mathbb{P})(ds,dx)
	+  \langle\epsilon_{J'}^\P,\mathbb{Y}_{s^-}\rangle \int_\Rset x(e^{g(x)}-1)F(dx)  ds.
\end{align*}
Thus, by Proposition~\ref{index-transf} we can write
$$
\begin{aligned}
		\langle\epsilon_{J},\mathbb{X}_t^\P\rangle 
		=& \langle(\epsilon_{J'}^\P;\e_{j_k})^{\thicksim},\mathbb{Y}_{t}\rangle
		+1_{\{j_{k-1}=j_k=0\}}\frac{1}{2}\langle\epsilon_{J''}^{\mathbb{P}}\otimes \e_{-1},\mathbb{Y}_{t}\rangle  \\
		&\qquad+\sum_{\e_{J_1}\otimes \e_{J_2}=\e_J}\frac{1}{|J_2|!}1_{\{|J_2|> 1\}}\langle(\epsilon_{J_1}^\P; \e_{S(J_2)})^{\thicksim},\mathbb{Y}_{t}\rangle 1_{\{J_2\in\{1,\dots,N\}^{|J_2|}\}}\\
		&+1_{\{j_k=0\}}\sum_{|I|\leq p}f^{I}\langle (\epsilon_{J'}^\P\shuffle \epsilon_{I})\otimes \e_0,\mathbb{Y}_{t}\rangle
		+1_{\{j_k=1\}}\langle\epsilon_{J'}^\P\otimes  \e_0,\mathbb{Y}_{t}\rangle \int_\Rset x (e^{g(x)}-1)F(dx).
	\end{aligned}
$$
concluding the proof.

\begin{remark}\phantomsection
 For notational convenience, the function $f$ in \eqref{eqn9} has been defined via multi-indices $I$ with indices in the set $\{-1,0,1,2,\dots,N\}\setminus\{0,1\}$. However, observe that the result of Proposition~\ref{prop1} still holds if we consider the larger set of multi-indices with indices in  $\{-1,0,1,2,\dots,N\}\setminus\{0\}$.

\end{remark}

\section{Auxiliary results on the signature of \cadlag rough paths}\label{sec:auxiliary}
In this section, we present some auxiliary results on the signature of w.g.~\cadlag rough paths. The main purpose is to give a more comprehensive account of the Marcus rough differential equation \eqref{Rdemarcus} solved by the minimal jump extension (truncated signature) of a w.g.~\cadlag rough path. To this end, we first elaborate on the different definitions of w.g.~rough paths often used in literature, then after some brief reminders about the rough integral, we conclude 
by providing an explicit form of equation \eqref{Rdemarcus} along with some properties of the signature of time-extended w.g.~\cadlag rough paths. 

\subsection{Equivalent definitions of weakly geometric \cadlag rough paths}\label{sec:equi}
In Definition \ref{def: roughpathsGROUP}, we introduced the notion of w.g.~\cadlag $p$-rough path as a group-valued path with a predetermined regularity. More often, however, especially when the main interest is in rough integration and differential equations, it is presented following the line below.  

\begin{definition}\label{rp1}
\begin{enumerate}
    \item Let $p\in[1,2)$. A \cadlag path $X:[0,1]\rightarrow\Rset^d$ is called \emph{w.g.~\cadlag $p$-rough path over $\Rset^d$} if $\|X\|_{p-var}:=\sup_{\mathcal{D}\subset[0,1]}\left(\sum_{t_i\in\mathcal{D}}\|X_{t_i,t_{i+1}}\|^p\right)^{\frac{1}{p}}<\infty$\footnote{For a partition $\mathcal{D}=\{0=t_0<t_1<\dots<t_k=1\}$ of $[0,1]$ we write $\sum_{t_i\in \mathcal{D}}$ for the summation over all points in $\mathcal{D}$. Moreover, for $X:[0,1]\rightarrow\Rset^d$, we set $X_{s,t}:=X_t-X_s$ for all $(s,t)\in[0,1]^2, \ s\leq t$.}.
    \item 		 Let $p\in [2,3)$ and $\Delta_1:=\{(s,t)\in[0,1]^2 \ | \ s\leq t\}$. A pair $\mathbf{X}=(X,\mathbb{X}^{(2)})$ is called \emph{w.g.~\cadlag $p$-rough path over $\Rset^d$} if  $X:[0,1]\rightarrow\Rset^d$ and $\mathbb{X}^{(2)}:\Delta_1\rightarrow(\Rset^d )^{\otimes 2}$ satisfy:
		 	\begin{enumerate}
		 		\item The map $[0,1]\ni t\longmapsto (X_{0,t},\mathbb{X}^{(2)}_{0,t})\in \Rset^d \oplus(\Rset^{d})^{\otimes2}$ is \cadlag.
		 		\item Chen's relation holds:
$\X_{s,t}=\X_{s,u}+\X_{u,t}+X_{s,u}\otimes X_{u,t}, \   0\leq s<u< t \leq 1.$
\item \label{analiii}$		\vvvert\mathbf{X}\vvvert_{p-var}:=\|X\|_{p- var}+\|\mathbb{X}^{(2)}\|_{p/2 - var}^{1/2}$\footnote{For $\mathbb{X}^{(2)}:\Delta_1\rightarrow(\Rset^d )^{\otimes 2}$, we set $\|\X\|_{p/2-var}:=\sup_{\mathcal{D}\subset[0,1]}\left(\sum_{t_i,t_{i+1}\in\mathcal{D}}\|\X_{t_i,t_{i+1}}\|^{p/2}\right)^{\frac{2}{p}}$.}$<\infty.$

     \item\label{symmetric} For all $0\leq s< t\leq 1$
 		$
Sym(\X_{s,t})\footnote{$Sym(\mathbb{X}^{(2)})$ and $Anti(\mathbb{X}^{(2)})$ denote the symmetric and antysimmetric part of $\mathbb{X}^{(2)}$, respectively.}=\frac{1}{2}X_{s,t}\otimes X_{s,t}.
 		$
		 	\end{enumerate}
    If moreover $\lim_{s \nnearrow t} Anti(\X_{s,t})=0, \ \text{for all }t\in[0,1]$, then $\mathbf{X}$ is said to be \emph{Marcus-like}.
	 
\end{enumerate}
	
\end{definition}

The information provided by Definition~\ref{def: roughpathsGROUP} and Definition~\ref{rp1}  is exactly the same, except for the value of the of $X$ (with values in $\Rset^d)$ at $0$, i.e.~$X_0$, which is omitted when adopting the Lie-group valued path point of view. Indeed, for $p\in[1,2)$, $[p]=1$, $G^1(\Rset^d)=\{1\}\oplus\Rset^d$ and the $d_{CC}$ distance on $G^1(\Rset^d)$ coincides with the Euclidean one. If $p\in[2,3)$ and $(X,\X)$ is  a w.g.~\cadlag $p$-rough path as introduced in Definition \ref{rp1}, consider the \cadlag path $[0,1]\ni t\mapsto \mathbf{X}_t:=(1,X_{0,t},\mathbb{X}_{0,t}^{(2)})\in T^{2}_1(\Rset^d)$. Define as in \eqref{eq:inv} the path increments via
    \begin{align}\label{eq:pathincrements}
        \mathbf{X}_{s,t}:=\mathbf{X}_s^{-1}\otimes \mathbf{X}_t=(1,X_{s,t},\mathbb{X}_{s,t}^{(2)}),
    \end{align}
   and notice that  $\mathbb{X}_{s,t}^{(2)}=\mathbb{X}_{0,t}^{(2)}-\mathbb{X}_{0,s}^{(2)}-X_{0,s}\otimes X_{s,t}$, which is in line with the Chen's relation. 
   Next, by property \ref{symmetric} of Definition~\ref{rp1}, it can be deduced that the path $[0,1]\ni t\mapsto \mathbf{X}_t$ takes values in fact in $G^{2}(\Rset^d)\subset T^{2}_1(\Rset^d)$.
Indeed, recall that $G^{2}(\Rset^d)=\exp^{(2)}(\mathfrak{g}^2(\Rset^d))$, and  $\mathfrak{g}^2(\Rset^d)=\Rset^d\otimes [\Rset^d,\Rset^d]$. Since $[\Rset^d,\Rset^d]$ is nothing but the space of antisymmetric $d\times d$ matrices, we have that
\begin{equation*}
    \mathbf{X}_t=(1,X_{0,t},\frac{1}{2}X_{0,t}^{\otimes 2}+Anti(\X_{0,t}))=\exp^{(2)}(X_{0,t},Anti(\X_{0,t}))\in G^{2}(\Rset^d) .
\end{equation*}
Finally, by Theorem~7.44 in \cite{FV:10} (see also the discussion in \cite{FZ:17}) the analytic condition  \ref{analiii} is equivalent to the one expressed in Definition \ref{def: roughpathsGROUP} by means of the distance $d_{CC}$ on $G^2(\Rset^d)$.

\subsection{Rough integration}

We briefly recall the notion of rough integral with respect to a w.g.~\cadlag $p-$rough path, for $p\in[2,3)$, as originally introduced in \cite{FS:17}. Let $V,W$ be two normed vector spaces and denote by $\mathcal{L}(V,W)$ the space of linear maps from $V$ into $W$. Remember that for $\mathbf{X}\in D^p([0,1],G^2(\Rset^d))$ the path increments are given in equation \eqref{eq:pathincrements}.
\begin{definition}\label{roughintegraldefinition}
Let $p\in[2,3)$, $\mathbf{X}\in D^p([0,1],G^2(\Rset^d))$ and $m\in\Nset$.
\begin{enumerate}
    \item  A pair of \cadlag paths $(Y,Y')$ is called a \textit{controlled rough path} (with respect to $X$) if $Y\in D^{p}([0,1],\mathcal{L}(\Rset^d,\Rset^m)$, $Y'\in D^{p}([0,1],\mathcal{L}(\Rset^d,\mathcal{L}(\Rset^d,\Rset^m))$ and $R:\Delta_1\rightarrow\mathcal{L}(\Rset^d,\Rset^m)$, defined by \[
    R_{s,t}:=Y_{s,t}-Y_s'X_{s,t}, \quad \text{for }\ (s,t)\in \Delta_1,
    \]
    has finite $\frac{p}{2}$-variation. 
    \item The \textit{rough integral} of a controlled rough path $(Y,Y')$ with respect to ${X}$ is defined as follows: for $t\in[0,1]$
    \begin{align}\label{RI}
        \int_0^tY_{s^-}d\mathbf{X}_s:=\lim_{(RRS)\ |\mathcal{D}|\rightarrow0}\sum_{t_{i}\in \mathcal{D}}Y_{t_{i}}X_{t_{i},t_{i+1}}+Y'_{t_{i}}\X_{t_{i},t_{i+1}},
    \end{align}
    where $\mathcal{D}$ denotes a partition of $[0,t]$ and the limit is understood in Refinement Riemann-Stieltjes (RRS) sense, as introduced in Definition~1 in \cite{FS:17}.
\end{enumerate}
\end{definition}
\begin{remark}
\begin{enumerate}
    \item The convergence in equation \eqref{RI} has been proved to hold also in Mesh Riemann-Stieltjes sense (see Definition~1.1 and Proposition~2.6 in \cite{FZ:17}).
    \item Remember that $\mathcal{L}(\Rset^d,\mathcal{L}(\Rset^d,\Rset^m))$ can be identify with the space $\mathcal{L}(\Rset^{d \times d},\Rset^m)$. This allows us to make
sense of the product $Y'\X$ in equation \eqref{RI}, which then takes values in $\Rset^m$.
\end{enumerate}
  \end{remark}
\subsection{Marcus rough differential equation}
We are now ready to analyze in detail the explicit form of equation \eqref{Rdemarcus}. 

Let $\mathbf{X}\in D^p([0,1],G^{[p]}(\Rset^d))$ and assume first that $p\in[2,3)$. Set $\Delta\mathbf{X}_{s}:= \lim_{u \nnearrow s} \mathbf{X}_{u,s}=(1,\Delta X_s,\Delta \mathbb{X}^{(2)}_s)$ for $s\in [0,1]$, and $\Nset \ni N\geq 3$. Then,  using \eqref{eqn3}, equation \eqref{Rdemarcusexplicit} reads
		\begin{align*} 
		\mathbb{X}^{N}_t=1+\int_0^t\mathbb{X}^{N}_{s^-}\otimes d\mathbf{X}_s+\sum_{0<s\leq t}\mathbb{X}^{N}_{s^-}\otimes\left(-\frac 1 2 (\Delta X_s)^{\otimes 2} +\sum_{k=2}^N \frac{1}{k!}(\Delta X_s,  Anti(\Delta \mathbb{X}_s))^{\otimes k}\right).
		\end{align*}
	Moreover, let $Z$ denote the rough integral 
$
	Z=\int_0^\cdot \mathbb{X}^{N}_{s^-}\otimes d\mathbf{X}_s,
	$
	and let $Z^{(k)}:=\pi_k(Z)$, $k=3,\dots,N$. Then, $Z^{(k)}$ is precisely the rough integral of the controlled rough path $(\pi_{k-1}(\mathbb{X}^{N}),\pi_{k-2}(\mathbb{X}^{N})$ with respect to $\mathbf{X}$.
		This means that 
		\begin{align}\label{eq:explicitMarcus}
		  \langle \epsilon_I,X^N_t\rangle =&\lim_{(RRS)\ |\mathcal{D}|\rightarrow0}\sum_{t_{i}\in \mathcal{D}}\langle \epsilon_{I'},\mathbb{X}^{N}_{t_i} \rangle  \langle \e_{i_{|I|}},X_{t_{i},t_{i+1}}\rangle +\langle \epsilon_{I''},\mathbb{X}^{N}_{t_i} \rangle  \langle \e_{i_{|I|-1},i_{|I|}},\mathbb{X}^{(2)}_{t_{i},t_{i+1}}\rangle \nonumber\\
		  \ -&\sum_{0<s\leq t}\sum_{\epsilon_{I_1}\otimes \epsilon_{I_2}=\epsilon_I}\langle\epsilon_{I_1},\mathbb{X}^N_{s^-}\rangle \langle \epsilon_{I_2},  \frac{1}{2}(\Delta X_s)^{\otimes 2}\rangle \\\
		   \ +&\sum_{0<s\leq t}\sum_{\epsilon_{I_1}\otimes \epsilon_{I_2}=\epsilon_I}\sum_{k=2}^N\frac{1}{|I_2|!}\langle\epsilon_{I_1},\mathbb{X}^N_{s^-}\rangle \langle \epsilon_{I_2}, (\Delta X_s, Anti(\Delta\X_s))^{\otimes k}\rangle \nonumber,
		\end{align}
		for each multi-index $I=(i_1,\dots,i_{|I|})\in\{1,\dots,d\}^{|I|}$ with $3\leq|I|\leq N$. Finally, notice that, if  $\mathbf{X}$ is  additionally Marcus-like, using \eqref{eqn3} equation \eqref{Rdemarcusexplicit} simplifies to
		\begin{align*}
			\mathbb{X}^{N}_t=1+\int_0^t\mathbb{X}^{N}_{s^-}\otimes d\mathbf{X}_s+\sum_{0<s\leq t}\mathbb{X}^{N}_{s^-}\otimes  (\exp^{(N)}(\Delta X_s)-1-\Delta X_s-\frac 1 2 (\Delta X_s)^{\otimes 2}),
		\end{align*}
		and the summation term can be expressed only by means of jumps at the first level of $\mathbf{X}$.

	If $p\in[1,2)$, equation \eqref{Rdemarcusexplicit} reads
\begin{equation}\label{eq:MarcusYoung}
\mathbb{X}^N_t=1+\int_0^t\mathbb{X}^N_{s^-}\otimes d X_s+\sum_{0<s\leq t}\mathbb{X}^N_{s^-}\otimes \{\exp^{(N)}(\Delta X_s)-1-\Delta X_s\}. 
		\end{equation} 
		Here, the integral is understood as Young integral. Specifically, $\langle \epsilon_{I},\mathbb{X}^N_t\rangle$ can be written as
		\begin{align*}
		    \lim_{(RRS)\ |\mathcal{D}|\rightarrow0}\sum_{t_{i}\in \mathcal{D}}\langle \epsilon_{I'},\mathbb{X}^{N}_{t_i}\rangle X^{i_{|I|}}_{t_{i},t_{i+1}}+\sum_{0<s\leq t}\sum_{\epsilon_{I_1}\otimes \epsilon_{I_2}=\epsilon_I}\frac{1}{|I_2|!}\langle\epsilon_{I_1},\mathbb{X}^N_{s^-}\rangle \langle \e_{I_2},(\Delta X_s)^{\otimes  |I_2|}\rangle1_{\{|I_2|>1\}},
		\end{align*}
		for each multi-index $I=(i_1,\dots,i_{|I|})\in\{1,\dots,d\}^{|I|}$ with $2\leq|I|\leq N$.

Finally, we conclude this section by highlighting a fundamental property of the signature of time-extended w.g.~rough paths. It is based on the following key lemma concerning the consistency of the Young and the rough integral.
\begin{lemma}\label{Rough=Young}
	Let $p\in[2,3)$ and $\mathbf{X}\in D^p([0,1],G^2(\Rset^d))$. Assume that for some $i\in \{1,\dots,d\}$ $\langle \e_i,\mathbf{X}\rangle $ has finite variation and $\langle \e_{(j,i)},\mathbf{X}\rangle$ is defined via the Young integral for all $j\in\{1,\dots,d\}$. Let $N\geq3$, $Y\in D^p([0,1],\mathcal{L}(\Rset^d,(\Rset^d)^{\otimes N}))$ and $Y'\in D^p([0,1],\mathcal{L}((\Rset^d),\mathcal{L}(\Rset^d,(\Rset^d)^{\otimes N})))$ be such that $(Y,Y')$ is a controlled rough path, and I:=$(j_1,\dots,j_{n-1},i)\in\{1,\dots,d\}^n$ a multi-index with length $n\leq N$. Then, the rough integral
	$
	Z_t:=\int_0^tY_{s^-}\otimes d\mathbf{X}_s, 
	$
	satisfies
$$
\langle \epsilon_I,Z_t\rangle = \int_0^t\langle \epsilon_{I'}, Y_{s^-}\rangle   d\langle \e_i,\mathbf{X}_s\rangle ,
$$
where the integral on the right hand side is a Young integral. 
\end{lemma}
\begin{proof}
Recall from Definition~\ref{roughintegraldefinition} that
\[
\langle \epsilon_I,Z_t\rangle=\lim_{(RRS)\ |\mathcal{D}|\rightarrow0}\sum_{t_{i}\in \mathcal{D}}\langle\epsilon_{I'}, Y_{t_{i}}\rangle \langle \e_i,\mathbf{X}_{t_{i},t_{i+1}}\rangle+\langle\epsilon_{I''}, Y'_{t_{i}}\rangle \langle \e_{(j_{n-1},i)},\mathbb{X}^{(2)}_{t_{i},t_{i+1}}\rangle.
\]
By Proposition~25 in \cite{FS:17} 
\[
\lim_{(RRS)\ |\mathcal{D}|\rightarrow0}\sum_{t_{i}\in \mathcal{D}}\langle\epsilon_{I'}, Y_{t_{i}}\rangle 
\langle \e_i,\mathbf{X}_{t_{i},t_{i+1}}\rangle=\int_0^t\langle \epsilon_{I'}, Y_{s^-}\rangle \  d\langle \e_i,\mathbf{X}_
s\rangle.
\]
Therefore, the claim follows if we show that
\[
\lim_{(RRS)\ |\mathcal{D}|\rightarrow0}\sum_{t_{i}\in \mathcal{D}}\langle\epsilon_{I''}, Y'_{t_{i}}\rangle \langle \e_{(
j_{n-1},1)},\mathbb{X}^{(2)}_{t_{i},t_{i+1}}\rangle=0.
\]
This is obtained by a slight modification of the proof of Theorem~35 in \cite{FS:17}. Using their notation, observe that in this case $[2,3)\ni p\neq q=1 $. Therefore, one needs to choose $2\leq p<p'< 3$ and $1< q'$ such that $\frac{1}{p'}+\frac{1}{q'}>1$, and the Young integral is still well defined. Then, the result follows by applying the same reasoning.
\end{proof}

\begin{proposition}
    \label{prop:rough=Young}
        Let $p\in[2,3)$ and $\mathbf{X}\in \widehat{D}^p([0,1],G^2(\Rset^{d+1}))$. Fix $N\geq 3$ and let $\widehat{\mathbb{X}}^N$ be the solution of the Marcus-type RDE \eqref{Rdemarcus}. Let $I=(i_1,\dots,i_{|I|-1},-1)\in\{-1,1,\dots,d\}^n$ be a multi-index with $|I|\leq N$. Then, 
        \[
        \langle \epsilon_I,\widehat{\mathbb{X}}^N_t\rangle=\int_0^t\langle \epsilon_{I'},\widehat{\mathbb{X}}^N_{s^-}\rangle \ ds.
        \]
\end{proposition}
\begin{proof}
    The proof follows by the definition of time-extended w.g.~\cadlag rough path, equation \eqref{eq:explicitMarcus} and Lemma ~\ref{Rough=Young}.
\end{proof}

\end{document}